\documentclass[11pt]{article}

\usepackage{graphicx}
\usepackage{natbib}
\usepackage{hyperref}
\usepackage{url}

\usepackage{amstext}
\usepackage{amssymb}
\usepackage{amsmath}
\usepackage{bbm}
\usepackage{booktabs}
\usepackage{caption}
\usepackage{subcaption}
\usepackage{multirow}

\usepackage{comment}

\usepackage{algorithm}
\usepackage{algorithmic}

\DeclareMathOperator*{\argmax}{arg\,max}

\newcommand{\eg}{\textit{e.g.}}
\newcommand{\ie}{\textit{i.e.}}

\newcommand{\Halmos}{$\hfill \square$}

\newtheorem{lemma}{Lemma}
\newtheorem{theorem}{Theorem}
\newtheorem{definition}{Definition}
\newtheorem{proposition}{Proposition}
\newtheorem{corollary}{Corollary}

\newcommand{\boxxx}[1]
 {\begin{center}\fbox{\begin{minipage}{13cm}#1\smallskip\end{minipage}}\end{center}}

\usepackage{tikz}

\allowdisplaybreaks

\title{{\bf Dynamic Facility Location under Cumulative Customer Demand}}

\author{}

\author{{\bf Warley Almeida Silva}\\
CIRRELT, Université de Montréal\\
warley.almeida.silva@umontreal.ca
\and
{\bf Margarida Carvalho} \\
    CIRRELT, Université de Montréal \\
carvalho@iro.umontreal.ca\and
{\bf Sanjay Dominik Jena}\\
CIRRELT, Université du Quebec à Montréal\\
jena.sanjay-dominik@uqam.ca}

\date{}

\newcommand{\mathset}[1]{\{{#1}\}}

\newcommand{\ourproblem}[1]{DFLP-CCD}

\newcommand{\ceil}[1]{\left\lceil #1  \right\rceil}

\begin{document}
\bibliographystyle{apalike}

\pagenumbering{gobble}

\maketitle

\begin{abstract}
    Dynamic facility location problems aim at placing one or more valuable resources over a planning horizon to meet customer demand.
    Existing literature commonly assumes that customer demand quantities are defined independently for each time period. 
    In many planning contexts, however, unmet demand carries over to future time periods. Unmet demand at some time periods may therefore affect decisions of subsequent time periods.
    This work studies a novel location problem, where the decision maker places facilities over time to capture cumulative customer demand. 
    We propose two mixed-integer programming formulations for this problem, and show that one of them has a tighter continuous relaxation and allows the representation of more general customer demand behaviour.
    We characterize the computational complexity for this problem, and analyze which problem characteristics result in NP-hardness. We then propose an exact branch-and-Benders-cut method, and show that this method is approximately five times faster, on average, than solving the tighter formulation directly in our computational experiments. 
    Our results also quantify the benefit of accounting for cumulative customer demand within the optimization framework, since the corresponding planning solutions perform much better than those obtained by ignoring cumulative demand or employing myopic heuristics.
    We also draw managerial insights on the quality of service perceived by customers when the provider places facilities under cumulative customer demand.
    
    \,

    \textbf{Keywords:} Facility Location; Multi-period Planning; Cumulative Customer Demand.

\end{abstract}

\newpage

\section{Introduction}

Dynamic facility location problems are a classical family of combinatorial problems that aim to place one or more valuable resources over a planning horizon \citep{nickel2019multi}.
Examples of applications include energy suppliers that need to place electric vehicle charging stations over time while facing different levels of customer demand \citep{lamontagneOptimisingElectricVehicle2023} and humanitarian organizations that must locate relief facilities over time while handling constant demand shifts
\citep{alizadehMultiPeriodMaximalCovering2021}.
In broad terms, the literature often designates the valuable resource to be located as a \textit{facility}, the entity seeking services at one of these facilities as a \textit{customer}, the quantity of service sought by a customer at a facility as \textit{demand}, and key moments in time where location decisions are assumed to be made as \textit{time periods} \citep{laporteLocationScience2019}.

In most dynamic location problems, customer demand is fixed for each time period and must be served while optimizing a performance measure \citep[see, \eg,][]{vanroyDualBasedProcedureDynamic1982, sahyouniFacilityLocationModel2007, jenaDynamicFacilityLocation2015}.
In the case where customer demand cannot be served completely (\eg, due to the lack of sufficient resources or technical restrictions), most works explicitly or implicitly maximize captured demand \citep[see, \eg,][]{gunawardaneDynamicVersionsSet1982, marinMultiperiodStochasticCovering2018a,vatsaCapacitatedMultiperiodMaximal2021}.
Unmet demand is typically assumed to vanish, thus not impacting location decisions of subsequent time periods.
In numerous planning contexts, however, unmet demand carries over to future time periods, until eventually being served.
For example, in temporary retail \citep{rosenbaumBenefitsPitfallsContemporary2021}, pop-up stores attempt to  capture customer demand for seasonal or luxury goods, which typically accumulate until it is met. In the healthcare sector, mobile clinics \citep{qiMobileFacilityRouting2017} have to serve patients' demand for treatment, which persists until it is met. In humanitarian logistics \citep{daneshvarTwostageStochasticPostdisaster2023}, temporary relief chains alleviate the need for shelter, food and medication. Unmet demand for such commodities accumulates over time and may even spread if not addressed in time, for example, in the case of epidemic outbreaks.
Despite unmet demand accumulating over time being a natural phenomenon, the literature on this type of demand behaviour is quite sparse.

\textit{Cumulative customer demand} may also appear in other multi-period planning problems such as vehicle routing (\eg, if a customer is visited at a later stage of the route rather than at the beginning, the demand of that customer may have already  increased), general network design and distribution networks (\eg, customer demand may increase as long as said customer is not connected to the distribution network), and production scheduling (\eg, products scheduled for later time periods may have a higher demand than initially expected).
However, to the best of our knowledge, the literature on cumulative customer demand for other multi-period planning problems is similarly quite sparse. 

This paper aims at filling this gap in the context of a dynamic facility location problem faced by a service provider where unmet demand carries over to future time periods.
The provider may be locating, for example, temporary retail stores \citep{rosenbaumBenefitsPitfallsContemporary2021, retail1}, mobile healthcare units \citep{busingRobustStrategicPlanning2021,healthcare2}, or teams providing consular services \citep{nica2020diaspora, consular1}.
We focus on real-world applications where the provider does not control how each customer patronizes available facilities.
Customers may decide that available facilities are sufficiently convenient to obtain service (\eg, retail goods,  medical supplies, and official documents) and, in this case, they satisfy their entire accumulated demand at the most preferred facility.
In contrast, customers may find available facilities not attractive enough and, in this case, they retain their unmet demand and wait for the next opportunity to obtain service.
Explicitly considering the accumulation of unmet demand and independent customer agency may result not only in a different location policy, but also in 
a higher collection of profit.
The provider therefore needs to be strategic about when and where to install facilities over the planning horizon.

We investigate a novel location problem named the \textit{Dynamic Facility Location Problem under Cumulative Customer Demand (\ourproblem{})}.
We concentrate on the deterministic case, where parameters modelling cumulative customer demand are either known or can be sufficiently well-estimated in advance.
More specifically, we contribute to the literature on location problems as follows:
\begin{enumerate}
    \item We introduce a novel multi-period deterministic location problem, referred to as \ourproblem{}, where the decision maker places facilities over time to capture cumulative customer demand. We characterize the computational complexity of some special cases of the \ourproblem{}, and provide insights on which problem characteristics render it NP-hard or even inapproximable.
    \item We propose two mixed-integer programming formulations, and show that one of them \textit{(i)} provides a tighter continuous relaxation and \textit{(ii)} allows the representation of more general customer demand behaviour. We also devise an exact Benders Decomposition, which (i) is five times faster, on average, than solving the tightest formulation directly and \textit{(ii)} proves optimality gaps about two times smaller in our computational experiments.
    \item We highlight the benefits of explicitly accounting for cumulative customer demand within the optimization framework, since myopic heuristics yield suboptimal solutions whenever customer demand is, in fact, cumulative -- in our computational experiments, on average, within $6\%$ of the optimal solution. We also draw insights on how different instance attributes seem to influence the structure of the optimal solution.
\end{enumerate}

The remainder of this paper is organized as follows.
Section~\ref{sec:literature} discusses the literature on dynamic location problems, paying close attention to how customer demand and customer preferences are modelled.
Section~\ref{sec:formulations} defines the \ourproblem{} and presents its two-mixed-integer programming formulations.
Section~\ref{sec:complexity} summarizes theoretical results concerning the computational complexity of the \ourproblem{}, and Section~\ref{sec:methods} describes the proposed exact and heuristic solution methods.
Section~\ref{sec:experiments} presents computational results and managerial insights, and Section~\ref{sec:conclusion} concludes with final remarks and future work.

\section{Literature Review}
\label{sec:literature}

In this section, we discuss literature relevant to our work. Given that demand accumulation is critical in multi-period planning contexts, we focus on multi-period demand modelling in Section~\ref{sub:demand}. We then discuss how to model customer preferences in Section~\ref{sub:preferences}.

\subsection{Demand Modelling}
\label{sub:demand}

Location problems seek to locate one or several facilities to satisfy customer demand (\ie, the \textit{quantity} of commodities or service sought by a customer). 
Literature on location problems dates back several decades and generally aims to strike a balance between the cost of locating facilities and the cost of serving customers. In multi-period contexts, customers tend to demand different quantities at each of the time periods. 
Since early works on multi-period location problems for a single facility \citep[see, \eg,][]{ballouDynamicWarehouseLocation1968,wesolowskyDynamicFacilityLocation1973} and for multiple facilities \citep[see, \eg,][]{wesolowskyMultiperiodLocationAllocationProblem1975, sweeneyImprovedLongRunModel1976}, it has been a classical assumption that the decision maker must satisfy the entire demand at each of the time periods.
Such an assumption remains predominant in the literature on location problems \citep[see, \eg,][]{jenaDynamicFacilityLocation2015,nickel2019multi}.

In certain planning contexts, demand cannot be satisfied completely. For example, healthcare clinics may not be within reach for some patients 
\citep{vatsaCapacitatedMultiperiodMaximal2021}.
In these contexts, the decision maker typically implicitly or explicitly maximizes captured demand to ensure the maximization of total profit, coverage, or welfare 
\citep[see, \eg,][]{ marinMultiperiodStochasticCovering2018a, alizadehMultiPeriodMaximalCovering2021, vatsaCapacitatedMultiperiodMaximal2021}.

A large body of the literature assumes that demand is deterministic, \ie, it can be sufficiently well forecast in advance. When demand is considered uncertain, the decision maker's ability to satisfy all demand may be particularly compromised. To address this challenge, the literature has extensively contributed to the development of stochastic optimization models, where facilities aim at satisfying demand represented by probability distributions \citep[see, \eg,][]{marinMultiperiodStochasticCovering2018a}, as well as robust optimization models \citep[see, \eg,][]{baron2011facility}, seeking to provide sufficient supply even when faced with worst-case demand.

Existing works propose diverse demand representations and requirements in terms of complete demand satisfaction. Yet, they generally assume that the unmet demand occurring at a time period is ignored in subsequent time periods, thus not impacting location decisions for the remainder of the planning horizon. 
In many planning contexts, this may be, however, a severe oversimplification of the underlying demand behaviour, as unmet demand of one time period may still be important in subsequent time periods. 
For example, in the context of humanitarian supply chains, \cite{daneshvarTwostageStochasticPostdisaster2023} explicitly model the accumulation of unmet demand for first-aid response resources such as medication, shelter and food. Demand that has not been satisfied at some time period remains critical. It may carry over to future time periods, and, in fact, lead to a spread of disease and therefore increase future demand.

One may naturally associate demand accumulation with the concept of inventory, where customer demand accumulates as inventories need more stocking over time \citep[see, \eg,][]{archettiBranchCutAlgorithmVendorManaged2007,zhangBendersDecompositionApproach2021}.
We highlight, however, that planning problems involving inventories differ considerably from the \ourproblem{}.
First, the former supposes that the decision maker decides how to conduct inventory restocking (\eg, which clients to replenish and by how many units), whereas the latter assumes that customers are the ones deciding how to satisfy their accumulated demand.
Second, the former has hard constraints to guarantee that client inventories are above a minimal level, whereas the latter does not have similar requirements.
In this sense, existing formulations for these planning problems involving inventories cannot be used to devise solutions for cumulative customer demand.

While cumulative demand can be seen as a specific case of decision-dependent demand, where the demand level depends on other model decisions, cumulative demand has not yet been addressed in the facility location literature.
Indeed, works that consider decision-dependent demand within facility location problems assume that the demand level  depends on the number of available facilities \citep[see, \eg,][]{basciftciDistributionallyRobustFacility2021}.
Furthermore, given the inherent model complexity, most works are restricted to single-period problems \citep[see, \eg,][]{huRobustFacilityLocation2025}, where the issue of accumulation of unmet demand over time does not occur.

To the best of our knowledge, cumulative customer demand has not yet been considered in the context of facility location or within other multi-period planning problems with a similar structure. Except for \cite{daneshvarTwostageStochasticPostdisaster2023}, it has not been studied in logistics or transportation problems either.

\subsection{Customer Preferences}
\label{sub:preferences}

The vast majority of the literature on facility location problems \citep[see, \eg,][]{vanroyDualBasedProcedureDynamic1982, jenaDynamicFacilityLocation2015, malladiFacilityLocationProblem2024} assumes that the decision maker decides through which facility each customer demand is satisfied.
However, in many planning contexts, customers are independent agents that decide how to patronize facilities based on their individual preferences.
For example, individuals may not buy an electric vehicle if charging stations are not sufficiently convenient, so modelling their preferences is key when developing the associated infrastructure \citep{lamontagneOptimisingElectricVehicle2023}.
This independent customer agency can be modelled by discrete choice models \citep[][]{bierlaireDiscreteChoiceModels1998}.

Early works on location problems with customer preferences suppose that customer behaviour is deterministic and can be expressed through rankings \citep[see, \eg,][]{hanjoulFacilityLocationProblem1987,canovasStrengthenedFormulationSimple2007}.
More precisely, each customer has a ranking over the locations and chooses to patronize the facility with the highest rank.
This choice model has the advantage of being easily integrated into the constraints of mixed-integer programming formulations, but the disadvantage of requiring customers to always make decisions based on the same criterion, thus not allowing the representation of more complex choice behaviour.

Recent works on location problems with customer preferences have employed parametric choice models, where each customer is assumed to patronize available facilities according to a probability distribution  \citep[\eg][]{linLocatingFacilitiesCompetition2022, qiSequentialCompetitiveFacility2022}.
These choice models rely on customer parameters estimated from historical data to quantify the utility of each location to each customer, and set patronizing probabilities proportional to these utility values.
On the one hand, parametric choice models can describe choice behaviours that are more complex than those represented through customer rankings (\eg, some level of uncertainty).
On the other hand, these choice models induce harder formulations from a methodological perspective (\eg, nonlinear objective functions), requiring some additional work to be solved in practice.

One may instead employ nonparametric ones, among which the rank-based choice model is quite popular \citep[see, \eg,][]{fariasNonparametricApproachModeling2013, vanryzinMarketDiscoveryAlgorithm2015, jenaPartiallyRankedChoice2020}.
This choice model assumes that each customer can be represented by multiple customer profiles, each one with a realization probability and ranking over locations.
We highlight that rank-based choice models perform better than parametric choice models when it comes to  describing real customer behaviour \citep{berbegliaComparativeEmpiricalStudy2022}.
In addition, their implementation shares many similarities with customer rankings, including its easy integration through constraints of mixed-integer programming formulations.
In this sense, we present the \ourproblem{} with customer rankings (\ie, a simple deterministic choice model), and explain in Online Appendix~A how our solution methods can be readily applied to a rank-based choice model.

\section{Mathematical Models}
\label{sec:formulations}

We define the \ourproblem{} in Section~\ref{sub:definition}, propose two formulations for it in Sections \ref{sub:formulation-DI} and \ref{sub:formulation-SI}, and highlight some of their properties in Section~\ref{sub:properties}.
We also discuss relevant problem extensions in Section~\ref{sec:ProblemExtensions}.

\subsection{Problem Definition}
\label{sub:definition}

We consider a service provider that locates facilities over time and profits from capturing customer demand.
Facilities can be located on a discrete set of candidate locations.
Demand can be captured from a discrete set of targeted customers. Unmet demand is assumed to accumulate over time.
More specifically, at each time period, a specific demand quantity, referred to as \textit{spawning demand}, adds to the already accumulated unmet demand. If not captured by a facility, the accumulated demand carries over to the next time period.
We assume that customers satisfy their entire accumulated demand when they go to a facility.
The provider does not control how customers attend available facilities. Instead, customers visit facilities according to their individual preferences. 
In other words, at each time period, each customer patronizes the facility with the highest rank among those made available by the provider, and may patronize no facility at all if available facilities are considered to be worse than not obtaining service.
Customers may, for example, prioritize facilities that are close, or those that have shorter waiting times to provide service. Only facilities meeting those criteria appear in the consideration set of the customer.
Given that the provider typically does not have sufficient facilities available to capture all customers at each of the time periods, the provider aims to determine the location of the facilities over the planning horizon that maximizes the total profit.

\subsubsection{Mathematical Notation.}

Throughout the rest of this paper, we employ bold letters and sets to denote vectors. Let $\mathcal{I}$ be the set of candidate locations, $\mathcal{J}$ be the set of targeted customers, and $\mathcal{T} = \{1, \ldots, T\}$ be the set of time periods.
Let also $\mathcal{T}^{S} = \{0\} \cup \mathcal{T}$ and $\mathcal{T}^{F} = \mathcal{T} \cup \{T+1\}$ be the set of time periods plus the start period $0$ and the final period $T + 1$, respectively. 
These artificial time periods become relevant later when writing mathematical formulations in Sections~\ref{sub:formulation-DI} and \ref{sub:formulation-SI}.

\subsubsection{Location Decisions.}
Since we are interested in the location of facilities over time, we define a ubiquitous set of decision variables $y^{t}_{i} \in \mathset{0,1}$ that equal $1$ if the provider opens a facility in location $i$ at time period $t$, $0$ otherwise.
Location decisions are constrained by a feasible set $\mathcal{Y}$, which may contain cardinality or budget constraints, phase-in or phase-out constraints, among others.
A location policy $\boldsymbol{y}$ incurs costs denoted by $f(\boldsymbol{y})$, which are related to assembling, disassembling, and transportation costs, among others.
While the literature considers a plethora of different feasible sets $\mathcal{Y}$  and cost functions $f(\boldsymbol{y})$, we here focus on a simple cardinality constraint to focus on understanding how cumulative customer demand influences location decisions over time. 
Specifically, we consider a simple feasible set of the form $\mathcal{Y} = \mathset{\sum_{i \in \mathcal{I}} y^{t}_{i} \leq h, \forall t \in \mathcal{T}}$, where $h \in \mathbb{Z}^{+}$ is the number of facilities available at each time period, and a trivial cost function $f(\boldsymbol{y}) = 0, \forall \boldsymbol{y} \in \mathcal{Y}$.
For example, in the context of mobile consular services, such a feasible set allows the consulate to allocate $h$ teams at each week, no matter their locations in the previous week.
We highlight, however, that our solution methods can be applied to different feasible sets $\mathcal{Y}$ and different cost functions $f(\boldsymbol{y})$ as long as they can be represented linearly.

\subsubsection{Customer Preferences.}
Each customer $j$ is only willing to attend a subset of candidate locations called the consideration set. 
Facilities in this set are ranked by preference, and a customer is assumed to attend the location within the consideration set that has the highest preference rank. 
More specifically, we denote $\succ_{j}$ the ranking of customer $j$ over the set of candidate locations $\mathcal{I}$, where the ranking $i \succ_{j} k$ indicates that customer $j$ prefers location $i$ over location $k$.
We employ an artificial location $0$ to allow customers to rank the choice of no service at all (\ie, the ranking $0 \succ_{j} i$ stipulates that customer $j$ prefers no service at all over patronizing location $i$).
All locations ranked after $0$ are not part of the consideration set.
We define an auxiliary parameter $a_{ij} \in \mathset{0,1}$ that equals $1$ if the ranking $i \succ_{j} 0$ holds for customer $j$ 
 in location $i$, $0$ otherwise.
Even though we here assume that each customer has a single ranking, our solution methods can readily account for representing situations where the choice behaviour of a customer is represented by a rank-based choice model \citep[][]{fariasNonparametricApproachModeling2013}.
This is achieved by creating one customer duplicate for each customer profile in the choice model, and scaling down parameters of these duplicates according to the underlying probability distribution (for more details, see Online Appendix~A). 

\subsubsection{Demand Accumulation.}
Let $d^{t}_{j} \in \mathbb{R}^{+}$ be the spawning demand of customer $j$ at time period $t$.
Inspired by \cite{daneshvarTwostageStochasticPostdisaster2023}, we formally define the accumulated demand of customer $j$ at the beginning of time period $t$ as $c^{t}_{j} (\boldsymbol{y}) = e_{j} u^{t-1}_{j} (\boldsymbol{y}) + d^{t}_{j}$, where $e_{j}\in \mathbb{R}^{+}$ is the spread factor of customer $j$ and $u^{t}_{j} (\boldsymbol{y})$ is the unmet demand of customer $j$ at the end of time period $t$.
Intuitively, the spread factor controls the percentage of the unmet demand from the previous time period that carries over to the current time period.
If $e_{j} > 1$ (respectively, $e_{j} < 1$), customer $j$ carries over strictly more (respectively, less) than the unmet demand $u^{t-1}_{j} (\boldsymbol{y})$.
Throughout the remainder of the paper, we assume that the unmet demand carries over entirely to the future (\ie, $e_{j} = 1, \forall j \in \mathcal{J}$).
We remark that our solution methods can still be adapted to account for different spread factors, but some positive theoretical results might not hold (namely, the approximation guarantees and the polynomial case further discussed in Section~\ref{sec:complexity}).
We define the unmet demand $u^{t}_{j} (\boldsymbol{y})$ of customer $j$ at time period $t$ as
    $u^{t}_{j} (\boldsymbol{y}) = \begin{cases}
        \left(1 - \max_{i \in \mathcal{I}} \mathset{a_{ij} y^{t}_{i}} \right) c^{t}_{j} (\boldsymbol{y}), & \text{ if } t \in \mathcal{T}, \\
        0, & \text{ if } t = 0.
    \end{cases}$
Note that the term $\max_{i \in \mathcal{I}} \mathset{a_{ij} y^{t}_{i}}$ equals $1$ if and only if customer $j$ has been captured by some location $i$ at time period $t$.

\begin{figure}[!ht]
     \centering
     \caption{Examples of spawning demands $d^{t}_{j}$ and their accumulated demand functions $c^{t}_{j} (\boldsymbol{y})$ for customer $j$ under an empty location policy $\boldsymbol{y} = \boldsymbol{0}$. Note that we interpolate the discrete points to facilitate the visualization of cumulative demand behaviours, but the planning horizon remains discrete.}
    \vspace{0.3cm}
     \begin{subfigure}[b]{0.49\textwidth}
        \centering
        \caption{Spawning demands $d^{t}_{j}$.}
        \begin{tikzpicture}[scale=.7, every node/.style={scale=.7}]
\draw[line width=0.5mm,thick,->] (0,0) -- (9.5,0);
\draw[line width=0.5mm,thick,->] (0,0) -- (0,5.5);
\draw (10, 0) node[anchor=mid] {$t$};
\draw (0, 6.0) node[anchor=mid] {$d^t_j$};
\draw (0,-0.5) node[anchor=mid] {$0$};
\draw (1,-0.5) node[anchor=mid] {$1$};
\draw (2,-0.5) node[anchor=mid] {$2$};
\draw (3,-0.5) node[anchor=mid] {$3$};
\draw (4,-0.5) node[anchor=mid] {$4$};
\draw (5,-0.5) node[anchor=mid] {$5$};
\draw (6,-0.5) node[anchor=mid] {$6$};
\draw (7,-0.5) node[anchor=mid] {$7$};
\draw (8,-0.5) node[anchor=mid] {$8$};
\draw (9,-0.5) node[anchor=mid] {$9$};

\draw (-0.5, 0) node[anchor=mid] {$0$};
\draw (-0.5, 1) node[anchor=mid] {$2$};
\draw (-0.5, 2) node[anchor=mid] {$4$};
\draw (-0.5, 3) node[anchor=mid] {$6$};
\draw (-0.5, 4) node[anchor=mid] {$8$};
\draw (-0.5, 5) node[anchor=mid] {$10$};

\draw[dotted, teal]  (1.5, 6)--(2.0, 6);
\draw[teal]  (2.0, 6) node[anchor=west] {$d^{t}_{j} = 5$};

\draw[dashed, red]  (1.5, 5.5)--(2.0, 5.5);
\draw[red]  (2.0, 5.5) node[anchor=west] {$d^{t}_{j} = 10 - t$};

\draw[dashdotted, blue]  (5.5, 6)--(6.0, 6);
\draw[blue]  (6.0, 6) node[anchor=west] {$d^{t}_{j} = t$};

\draw[solid, orange]  (5.5, 5.5)--(6.0, 5.5);
\draw[orange]  (6.0, 5.5) node[anchor=west] {$d^{t}_{j} = 5 \sin{t} + 5$};

\draw[dotted, teal] (0.0,2.5)--(1.0,2.5);
\draw[teal,fill=teal] (1.0,2.5) circle (2pt);
\draw[dotted, teal] (1.0,2.5)--(2.0,2.5);
\draw[teal,fill=teal] (2.0,2.5) circle (2pt);
\draw[dotted, teal] (2.0,2.5)--(3.0,2.5);
\draw[teal,fill=teal] (3.0,2.5) circle (2pt);
\draw[dotted, teal] (3.0,2.5)--(4.0,2.5);
\draw[teal,fill=teal] (4.0,2.5) circle (2pt);
\draw[dotted, teal] (4.0,2.5)--(5.0,2.5);
\draw[teal,fill=teal] (5.0,2.5) circle (2pt);
\draw[dotted, teal] (5.0,2.5)--(6.0,2.5);
\draw[teal,fill=teal] (6.0,2.5) circle (2pt);
\draw[dotted, teal] (6.0,2.5)--(7.0,2.5);
\draw[teal,fill=teal] (7.0,2.5) circle (2pt);
\draw[dotted, teal] (7.0,2.5)--(8.0,2.5);
\draw[teal,fill=teal] (8.0,2.5) circle (2pt);
\draw[dotted, teal] (8.0,2.5)--(9.0,2.5);
\draw[teal,fill=teal] (9.0,2.5) circle (2pt);
\draw[dashdotted, blue] (0.0,0.0)--(1.0,0.5);
\draw[blue,fill=blue] (1.0,0.5) circle (2pt);
\draw[dashdotted, blue] (1.0,0.5)--(2.0,1.0);
\draw[blue,fill=blue] (2.0,1.0) circle (2pt);
\draw[dashdotted, blue] (2.0,1.0)--(3.0,1.5);
\draw[blue,fill=blue] (3.0,1.5) circle (2pt);
\draw[dashdotted, blue] (3.0,1.5)--(4.0,2.0);
\draw[blue,fill=blue] (4.0,2.0) circle (2pt);
\draw[dashdotted, blue] (4.0,2.0)--(5.0,2.5);
\draw[blue,fill=blue] (5.0,2.5) circle (2pt);
\draw[dashdotted, blue] (5.0,2.5)--(6.0,3.0);
\draw[blue,fill=blue] (6.0,3.0) circle (2pt);
\draw[dashdotted, blue] (6.0,3.0)--(7.0,3.5);
\draw[blue,fill=blue] (7.0,3.5) circle (2pt);
\draw[dashdotted, blue] (7.0,3.5)--(8.0,4.0);
\draw[blue,fill=blue] (8.0,4.0) circle (2pt);
\draw[dashdotted, blue] (8.0,4.0)--(9.0,4.5);
\draw[blue,fill=blue] (9.0,4.5) circle (2pt);
\draw[dashed, red] (0.0,5.0)--(1.0,4.5);
\draw[red,fill=red] (1.0,4.5) circle (2pt);
\draw[dashed, red] (1.0,4.5)--(2.0,4.0);
\draw[red,fill=red] (2.0,4.0) circle (2pt);
\draw[dashed, red] (2.0,4.0)--(3.0,3.5);
\draw[red,fill=red] (3.0,3.5) circle (2pt);
\draw[dashed, red] (3.0,3.5)--(4.0,3.0);
\draw[red,fill=red] (4.0,3.0) circle (2pt);
\draw[dashed, red] (4.0,3.0)--(5.0,2.5);
\draw[red,fill=red] (5.0,2.5) circle (2pt);
\draw[dashed, red] (5.0,2.5)--(6.0,2.0);
\draw[red,fill=red] (6.0,2.0) circle (2pt);
\draw[dashed, red] (6.0,2.0)--(7.0,1.5);
\draw[red,fill=red] (7.0,1.5) circle (2pt);
\draw[dashed, red] (7.0,1.5)--(8.0,1.0);
\draw[red,fill=red] (8.0,1.0) circle (2pt);
\draw[dashed, red] (8.0,1.0)--(9.0,0.5);
\draw[red,fill=red] (9.0,0.5) circle (2pt);

\draw[solid, orange] (0.0,2.5)--(0.1,2.75);
\draw[solid, orange] (0.1,2.75)--(0.2,3.0);
\draw[solid, orange] (0.2,3.0)--(0.3,3.24);
\draw[solid, orange] (0.3,3.24)--(0.4,3.47);
\draw[solid, orange] (0.4,3.47)--(0.5,3.7);
\draw[solid, orange] (0.5,3.7)--(0.6,3.91);
\draw[solid, orange] (0.6,3.91)--(0.7,4.11);
\draw[solid, orange] (0.7,4.11)--(0.8,4.29);
\draw[solid, orange] (0.8,4.29)--(0.9,4.46);
\draw[solid, orange] (0.9,4.46)--(1.0,4.6);
\draw[orange,fill=orange] (1.0,4.6) circle (2pt);
\draw[solid, orange] (1.0,4.6)--(1.1,4.73);
\draw[solid, orange] (1.1,4.73)--(1.2,4.83);
\draw[solid, orange] (1.2,4.83)--(1.3,4.91);
\draw[solid, orange] (1.3,4.91)--(1.4,4.96);
\draw[solid, orange] (1.4,4.96)--(1.5,4.99);
\draw[solid, orange] (1.5,4.99)--(1.6,5.0);
\draw[solid, orange] (1.6,5.0)--(1.7,4.98);
\draw[solid, orange] (1.7,4.98)--(1.8,4.93);
\draw[solid, orange] (1.8,4.93)--(1.9,4.87);
\draw[solid, orange] (1.9,4.87)--(2.0,4.77);
\draw[orange,fill=orange] (2.0,4.77) circle (2pt);
\draw[solid, orange] (2.0,4.77)--(2.1,4.66);
\draw[solid, orange] (2.1,4.66)--(2.2,4.52);
\draw[solid, orange] (2.2,4.52)--(2.3,4.36);
\draw[solid, orange] (2.3,4.36)--(2.4,4.19);
\draw[solid, orange] (2.4,4.19)--(2.5,4.0);
\draw[solid, orange] (2.5,4.0)--(2.6,3.79);
\draw[solid, orange] (2.6,3.79)--(2.7,3.57);
\draw[solid, orange] (2.7,3.57)--(2.8,3.34);
\draw[solid, orange] (2.8,3.34)--(2.9,3.1);
\draw[solid, orange] (2.9,3.1)--(3.0,2.85);
\draw[orange,fill=orange] (3.0,2.85) circle (2pt);
\draw[solid, orange] (3.0,2.85)--(3.1,2.6);
\draw[solid, orange] (3.1,2.6)--(3.2,2.35);
\draw[solid, orange] (3.2,2.35)--(3.3,2.11);
\draw[solid, orange] (3.3,2.11)--(3.4,1.86);
\draw[solid, orange] (3.4,1.86)--(3.5,1.62);
\draw[solid, orange] (3.5,1.62)--(3.6,1.39);
\draw[solid, orange] (3.6,1.39)--(3.7,1.18);
\draw[solid, orange] (3.7,1.18)--(3.8,0.97);
\draw[solid, orange] (3.8,0.97)--(3.9,0.78);
\draw[solid, orange] (3.9,0.78)--(4.0,0.61);
\draw[orange,fill=orange] (4.0,0.61) circle (2pt);
\draw[solid, orange] (4.0,0.61)--(4.1,0.45);
\draw[solid, orange] (4.1,0.45)--(4.2,0.32);
\draw[solid, orange] (4.2,0.32)--(4.3,0.21);
\draw[solid, orange] (4.3,0.21)--(4.4,0.12);
\draw[solid, orange] (4.4,0.12)--(4.5,0.06);
\draw[solid, orange] (4.5,0.06)--(4.6,0.02);
\draw[solid, orange] (4.6,0.02)--(4.7,0.0);
\draw[solid, orange] (4.7,0.0)--(4.8,0.01);
\draw[solid, orange] (4.8,0.01)--(4.9,0.04);
\draw[solid, orange] (4.9,0.04)--(5.0,0.1);
\draw[orange,fill=orange] (5.0,0.1) circle (2pt);
\draw[solid, orange] (5.0,0.1)--(5.1,0.19);
\draw[solid, orange] (5.1,0.19)--(5.2,0.29);
\draw[solid, orange] (5.2,0.29)--(5.3,0.42);
\draw[solid, orange] (5.3,0.42)--(5.4,0.57);
\draw[solid, orange] (5.4,0.57)--(5.5,0.74);
\draw[solid, orange] (5.5,0.74)--(5.6,0.92);
\draw[solid, orange] (5.6,0.92)--(5.7,1.12);
\draw[solid, orange] (5.7,1.12)--(5.8,1.34);
\draw[solid, orange] (5.8,1.34)--(5.9,1.57);
\draw[solid, orange] (5.9,1.57)--(6.0,1.8);
\draw[orange,fill=orange] (6.0,1.8) circle (2pt);
\draw[solid, orange] (6.0,1.8)--(6.1,2.04);
\draw[solid, orange] (6.1,2.04)--(6.2,2.29);
\draw[solid, orange] (6.2,2.29)--(6.3,2.54);
\draw[solid, orange] (6.3,2.54)--(6.4,2.79);
\draw[solid, orange] (6.4,2.79)--(6.5,3.04);
\draw[solid, orange] (6.5,3.04)--(6.6,3.28);
\draw[solid, orange] (6.6,3.28)--(6.7,3.51);
\draw[solid, orange] (6.7,3.51)--(6.8,3.74);
\draw[solid, orange] (6.8,3.74)--(6.9,3.95);
\draw[solid, orange] (6.9,3.95)--(7.0,4.14);
\draw[orange,fill=orange] (7.0,4.14) circle (2pt);
\draw[solid, orange] (7.0,4.14)--(7.1,4.32);
\draw[solid, orange] (7.1,4.32)--(7.2,4.48);
\draw[solid, orange] (7.2,4.48)--(7.3,4.63);
\draw[solid, orange] (7.3,4.63)--(7.4,4.75);
\draw[solid, orange] (7.4,4.75)--(7.5,4.84);
\draw[solid, orange] (7.5,4.84)--(7.6,4.92);
\draw[solid, orange] (7.6,4.92)--(7.7,4.97);
\draw[solid, orange] (7.7,4.97)--(7.8,5.0);
\draw[solid, orange] (7.8,5.0)--(7.9,5.0);
\draw[solid, orange] (7.9,5.0)--(8.0,4.97);
\draw[orange,fill=orange] (8.0,4.97) circle (2pt);
\draw[solid, orange] (8.0,4.97)--(8.1,4.92);
\draw[solid, orange] (8.1,4.92)--(8.2,4.85);
\draw[solid, orange] (8.2,4.85)--(8.3,4.76);
\draw[solid, orange] (8.3,4.76)--(8.4,4.64);
\draw[solid, orange] (8.4,4.64)--(8.5,4.5);
\draw[solid, orange] (8.5,4.5)--(8.6,4.34);
\draw[solid, orange] (8.6,4.34)--(8.7,4.16);
\draw[solid, orange] (8.7,4.16)--(8.8,3.96);
\draw[solid, orange] (8.8,3.96)--(8.9,3.75);
\draw[solid, orange] (8.9,3.75)--(9.0,3.53);
\draw[orange,fill=orange] (9.0,3.53) circle (2pt);
\draw[solid, orange] (9.0,3.53)--(9.1,3.3);
\draw[solid, orange] (9.1,3.3)--(9.2,3.06);
\draw[solid, orange] (9.2,3.06)--(9.3,2.81);
\draw[solid, orange] (9.3,2.81)--(9.4,2.56);
\draw[solid, orange] (9.4,2.56)--(9.5,2.31);

\end{tikzpicture}
        \label{fig:spawning}
     \end{subfigure}
    \hfill
     \begin{subfigure}[b]{0.50\textwidth}
        \centering
        \caption{Accumulated demand functions $c^{t}_{j} (\boldsymbol{y} = \boldsymbol{0}) = \sum^{t}_{k = 1} d^{k}_{j}$.}
        \begin{tikzpicture}[scale=.7, every node/.style={scale=.7}]
\draw[line width=0.5mm,thick,->] (0,0) -- (9.5,0);
\draw[line width=0.5mm,thick,->] (0,0) -- (0,5.5);
\draw (10, 0) node[anchor=mid] {$t$};
\draw (0,6.0) node[anchor=mid] {$c^t_j(\boldsymbol{y} = \boldsymbol{0})$};
\draw (0,-0.5) node[anchor=mid] {$0$};
\draw (1,-0.5) node[anchor=mid] {$1$};
\draw (2,-0.5) node[anchor=mid] {$2$};
\draw (3,-0.5) node[anchor=mid] {$3$};
\draw (4,-0.5) node[anchor=mid] {$4$};
\draw (5,-0.5) node[anchor=mid] {$5$};
\draw (6,-0.5) node[anchor=mid] {$6$};
\draw (7,-0.5) node[anchor=mid] {$7$};
\draw (8,-0.5) node[anchor=mid] {$8$};
\draw (9,-0.5) node[anchor=mid] {$9$};

\draw (-0.5, 0) node[anchor=mid] {$0$};
\draw (-0.5, 1) node[anchor=mid] {$10$};
\draw (-0.5, 2) node[anchor=mid] {$20$};
\draw (-0.5, 3) node[anchor=mid] {$30$};
\draw (-0.5, 4) node[anchor=mid] {$40$};
\draw (-0.5, 5) node[anchor=mid] {$50$};


\draw[dotted, teal] (0.0,0.0)--(1.0,0.5);
\draw[teal,fill=teal] (1.0,0.5) circle (2pt);
\draw[dotted, teal] (1.0,0.5)--(2.0,1.0);
\draw[teal,fill=teal] (2.0,1.0) circle (2pt);
\draw[dotted, teal] (2.0,1.0)--(3.0,1.5);
\draw[teal,fill=teal] (3.0,1.5) circle (2pt);
\draw[dotted, teal] (3.0,1.5)--(4.0,2.0);
\draw[teal,fill=teal] (4.0,2.0) circle (2pt);
\draw[dotted, teal] (4.0,2.0)--(5.0,2.5);
\draw[teal,fill=teal] (5.0,2.5) circle (2pt);
\draw[dotted, teal] (5.0,2.5)--(6.0,3.0);
\draw[teal,fill=teal] (6.0,3.0) circle (2pt);
\draw[dotted, teal] (6.0,3.0)--(7.0,3.5);
\draw[teal,fill=teal] (7.0,3.5) circle (2pt);
\draw[dotted, teal] (7.0,3.5)--(8.0,4.0);
\draw[teal,fill=teal] (8.0,4.0) circle (2pt);
\draw[dotted, teal] (8.0,4.0)--(9.0,4.5);
\draw[teal,fill=teal] (9.0,4.5) circle (2pt);
\draw[dashdotted, blue] (0.0,0.0)--(1.0,0.1);
\draw[blue,fill=blue] (1.0,0.1) circle (2pt);
\draw[dashdotted, blue] (1.0,0.1)--(2.0,0.3);
\draw[blue,fill=blue] (2.0,0.3) circle (2pt);
\draw[dashdotted, blue] (2.0,0.3)--(3.0,0.6);
\draw[blue,fill=blue] (3.0,0.6) circle (2pt);
\draw[dashdotted, blue] (3.0,0.6)--(4.0,1.0);
\draw[blue,fill=blue] (4.0,1.0) circle (2pt);
\draw[dashdotted, blue] (4.0,1.0)--(5.0,1.5);
\draw[blue,fill=blue] (5.0,1.5) circle (2pt);
\draw[dashdotted, blue] (5.0,1.5)--(6.0,2.1);
\draw[blue,fill=blue] (6.0,2.1) circle (2pt);
\draw[dashdotted, blue] (6.0,2.1)--(7.0,2.8);
\draw[blue,fill=blue] (7.0,2.8) circle (2pt);
\draw[dashdotted, blue] (7.0,2.8)--(8.0,3.6);
\draw[blue,fill=blue] (8.0,3.6) circle (2pt);
\draw[dashdotted, blue] (8.0,3.6)--(9.0,4.5);
\draw[blue,fill=blue] (9.0,4.5) circle (2pt);
\draw[dashed, red] (0.0,0.0)--(1.0,0.9);
\draw[red,fill=red] (1.0,0.9) circle (2pt);
\draw[dashed, red] (1.0,0.9)--(2.0,1.7);
\draw[red,fill=red] (2.0,1.7) circle (2pt);
\draw[dashed, red] (2.0,1.7)--(3.0,2.4);
\draw[red,fill=red] (3.0,2.4) circle (2pt);
\draw[dashed, red] (3.0,2.4)--(4.0,3.0);
\draw[red,fill=red] (4.0,3.0) circle (2pt);
\draw[dashed, red] (4.0,3.0)--(5.0,3.5);
\draw[red,fill=red] (5.0,3.5) circle (2pt);
\draw[dashed, red] (5.0,3.5)--(6.0,3.9);
\draw[red,fill=red] (6.0,3.9) circle (2pt);
\draw[dashed, red] (6.0,3.9)--(7.0,4.2);
\draw[red,fill=red] (7.0,4.2) circle (2pt);
\draw[dashed, red] (7.0,4.2)--(8.0,4.4);
\draw[red,fill=red] (8.0,4.4) circle (2pt);
\draw[dashed, red] (8.0,4.4)--(9.0,4.5);
\draw[red,fill=red] (9.0,4.5) circle (2pt);
\draw[solid, orange] (0.0,0.0)--(1.0,0.92);
\draw[orange,fill=orange] (1.0,0.92) circle (2pt);
\draw[solid, orange] (1.0,0.92)--(2.0,1.88);
\draw[orange,fill=orange] (2.0,1.88) circle (2pt);
\draw[solid, orange] (2.0,1.88)--(3.0,2.45);
\draw[orange,fill=orange] (3.0,2.45) circle (2pt);
\draw[solid, orange] (3.0,2.45)--(4.0,2.57);
\draw[orange,fill=orange] (4.0,2.57) circle (2pt);
\draw[solid, orange] (4.0,2.57)--(5.0,2.59);
\draw[orange,fill=orange] (5.0,2.59) circle (2pt);
\draw[solid, orange] (5.0,2.59)--(6.0,2.95);
\draw[orange,fill=orange] (6.0,2.95) circle (2pt);
\draw[solid, orange] (6.0,2.95)--(7.0,3.78);
\draw[orange,fill=orange] (7.0,3.78) circle (2pt);
\draw[solid, orange] (7.0,3.78)--(8.0,4.77);
\draw[orange,fill=orange] (8.0,4.77) circle (2pt);
\draw[solid, orange] (8.0,4.77)--(9.0,5.48);
\draw[orange,fill=orange] (9.0,5.48) circle (2pt);

\end{tikzpicture}
        \label{fig:accumulated}
     \end{subfigure}
    \label{fig:demand-functions}
\end{figure}

Figure~\ref{fig:demand-functions} provides examples of spawning demands $d^{t}_{j}$ and their accumulated demands $c^{t}_{j} (\boldsymbol{y})$ under an empty location policy $\boldsymbol{y} = \boldsymbol{0}$.
If the provider had to choose between these four customers at time period $t = 6$ based on spawning demands, customer blue would be the most profitable. 
However, under cumulative customer demand, customer blue is rather the least profitable at time period $t = 6$, and customer red should be chosen.
Cumulative customer demand introduces a layer of complexity when placing facilities, as location decisions from time periods $1,2, \ldots, t-1$ directly affect how much customer demand is available at time period $t$, and embedding the accumulated demand function $c^{t}_{j} (\boldsymbol{y})$ into mixed-integer programming formulations is not trivial.

\subsubsection{Reward Computation.}
Let $r_{i} \in \mathbb{R}^{+}$ be the reward per demand unit captured through location $i$.
For a location policy $\boldsymbol{y}$, if the provider opens a location $i$ at time period $t$ (\ie, $y^{t}_{i} = 1$) and customer $j$ decides to patronize it (\ie, $i \succ_{j} 0$ and $i \succ_{j} k, \forall k \in \mathcal{I}$ such that $y^{t}_{k} = 1$), we compute its marginal reward as the product $r_{i} c^{t}_{j} (\boldsymbol{y})$.
Note that we can seamlessly account for customer-dependent rewards by employing $r_{ij}$ without affecting the theoretical results presented in Section~\ref{sec:complexity}.
We highlight that the facility through which the provider captures customer $j$ directly influences the marginal reward obtained by the provider.

\subsection{Double-Index Formulation}
\label{sub:formulation-DI}

We propose first a general formulation that is capable of representing complex accumulation patterns.
While a formulation based on classical formulation techniques would track the accumulated customer demand at each time period (see Single-Index Formulation in Section \ref{sub:formulation-SI}), the proposed general formulation relies on pairs $(\ell, t)$ of time periods, $\ell{} < t$, to represent the accumulated demand function $c^{t}_{j} (\boldsymbol{y})$.
We therefore refer to this formulation as the \textit{Double-Index Formulation (DIF)}.

Let decision variables $x^{\ell{}t}_{ij} \in \mathset{0,1}$ be $1$ if customer $j$ patronizes location $i$ at time period $t$ after being last captured at time period $\ell{}$, and $0$ otherwise.
If customer $j$ was last captured at time period $\ell{}$, we can easily compute its accumulated demand at time period $t$ as $D^{\ell{}t}_{j} = \sum_{\substack{s \in \mathcal{T} \\ s > \ell{} \\ s \leq t}} d^{s}_{j}$, and its marginal reward through location $i$ as $G^{\ell{}t}_{ij} = r_{i} D^{\ell{}t}_{j}$.
The DIF can be written as follows:
\begin{subequations}
    \label{pgm:network-dflp-ccd}
    \begin{align}
        \max_{\boldsymbol{x}, \boldsymbol{y}} \quad
        &  \sum_{t \in \mathcal{T}} \sum_{\substack{\ell{} \in \mathcal{T}^{S} : \\ \ell{} < t}} \sum_{i \in \mathcal{I}} \sum_{j \in \mathcal{J}} G^{\ell{}t}_{ij}  x^{\ell{}t}_{ij} - f(\boldsymbol{y})
        && \label{eq:network-dflp-ccd-obj}\\
        \text{s.t.} \quad
        & \boldsymbol{y} \in \mathcal{Y}
        && \label{eq:network-dflp-ccd-ct1} \\
        & \sum_{\substack{\ell{} \in \mathcal{T}^{S} \\ \ell{} < t}} x^{\ell{}t}_{ij} \leq a_{ij} y^{t}_{i} 
        && \substack{\forall i \in \mathcal{I}, \forall j \in \mathcal{J}, \\ \forall t \in \mathcal{T}} \label{eq:network-dflp-ccd-ct2} \\
        & \sum_{\substack{\ell{} \in \mathcal{T}^{S}: \\ \ell{} < t}} \sum_{k \in \mathcal{I}} x^{\ell{}t}_{kj} \geq a_{ij} y^{t}_{i}
        && \substack{\forall i \in \mathcal{I}, \forall j \in \mathcal{J}, \\ \forall t \in \mathcal{T}} \label{eq:network-dflp-ccd-ct3} \\
        & a_{ij} y^{t}_{i} + \sum_{\substack{\ell{} \in \mathcal{T}: \\ \ell{} < t}} \sum_{\substack{k \in \mathcal{I}: \\ i \succ_{j} k}} x^{\ell{}t}_{kj} \leq a_{ij} 
        &&  \substack{\forall i \in \mathcal{I}, \forall j \in \mathcal{J}, \\ \forall t \in \mathcal{T}} \label{eq:network-dflp-ccd-ct4} \\
        & \sum_{\substack{s \in \mathcal{T}^{F}: \\ s > t}} \sum_{i \in \mathcal{I}} x^{ts}_{ij} - \sum_{\substack{s \in \mathcal{T}^{S}: \\ s < t}} \sum_{i \in \mathcal{I}} x^{st}_{ij} = 0
        && \substack{\forall j \in \mathcal{J}, \\ \forall t \in \mathcal{T}} \label{eq:network-dflp-ccd-ct5} \\
        & \sum_{s \in \mathcal{T}^{F}} \sum_{i \in \mathcal{I}} x^{0s}_{ij} = 1
        && \forall j \in \mathcal{J} \label{eq:network-dflp-ccd-ct6} \\
        & x^{\ell{}t}_{ij} \in \mathset{0,1}
        && \substack{\forall i \in \mathcal{I}, \forall j \in \mathcal{J}, \\ \forall \ell{} \in \mathcal{T}^{S}, \forall t \in \mathcal{T}^{F} : \\ \ell{} < t} \label{eq:network-dflp-ccd-dm1} \\
        & y^{t}_{i} \in \mathset{0,1}
        && \substack{\forall i \in \mathcal{I}, \\ \forall t \in \mathcal{T}}. \label{eq:network-dflp-ccd-dm2}
    \end{align}
\end{subequations}

Objective Function \eqref{eq:network-dflp-ccd-obj} maximizes the total profit obtained by the provider (\ie, total reward minus total costs).
Constraints \eqref{eq:network-dflp-ccd-ct1} guarantee that location decisions are feasible.
Constraints \eqref{eq:network-dflp-ccd-ct2} prevent customers from patronizing facilities that are not open by the provider, while Constraints \eqref{eq:network-dflp-ccd-ct3} force customers to patronize one of the available facilities if at least one of them is better than no service at all.
Constraints \eqref{eq:network-dflp-ccd-ct4} enforce customer preferences (\ie, customer $j$ cannot patronize a location $k$ less preferred than location $i$ if there is a facility in location $i$ at time period $t$).
Constraints \eqref{eq:network-dflp-ccd-ct5}--\eqref{eq:network-dflp-ccd-ct6} preserve the capture flow for customers throughout the planning horizon (\ie, if customer $j$ has been captured at time period $t$ after lastly being captured at some previous time period $s_{1}$, it must also be captured  at some future time period $s_{2}$ after lastly being captured at time period $t$). 
Finally, Constraints \eqref{eq:network-dflp-ccd-dm1}--\eqref{eq:network-dflp-ccd-dm2} define feasible variable domains.
Note that, once variables $y^{t}_{i}$ are fixed to binary values, Constraints \eqref{eq:network-dflp-ccd-ct2}--\eqref{eq:network-dflp-ccd-ct6} define a totally unimodular matrix. As the resulting formulation has integral extreme points, we can relax binary constraints on variables $x^{\ell{}t}_{ij}$.

\subsection{Single-Index Formulation}
\label{sub:formulation-SI}

We now introduce an alternative formulation 
based on classical formulation techniques \citep[\eg, ][]{daneshvarTwostageStochasticPostdisaster2023}.
This formulation represents the accumulated demand function $c^{t}_{j} (\boldsymbol{y})$ and the unmet demand function $u^{t}_{j} (\boldsymbol{y})$ through two sets of continuous decision variables, each with a single time index $t$.
As such, we refer to this formulation as the \textit{Single-Index Formulation (SIF)}.

Let $c^{t}_{j} \in \mathbb{R}^{+}$ and $u^{t}_{j} \in \mathbb{R}^{+}$ be two sets of decision variables that store the accumulated demand of customer $j$ at the beginning of time period $t$ and the unmet demand of customer $j$ at the end of time period $t$, respectively.
Let also $x^{t}_{ij} \in \mathset{0,1}$ be a set of decision variables that equals $1$ if customer $j$ patronizes a facility in location $i$ at time period $t$, $0$ otherwise.
The SIF can be written as follows:
\begin{subequations}
    \label{pgm:nonlinear-dflp-ccd}
    \begin{align}
        \max_{\boldsymbol{u}, \boldsymbol{c}, \boldsymbol{x}, \boldsymbol{y}}  \quad
        &  \sum_{t \in \mathcal{T}} \sum_{i \in \mathcal{I}} \sum_{j \in \mathcal{J}} r_{i} c^{t}_{j} x^{t}_{ij} - f(\boldsymbol{y}) 
        && \label{eq:nonlinear-dflp-ccd-obj} \\
        \text{s.t.} \quad
        & \boldsymbol{y} \in \mathcal{Y} 
        && \label{eq:nonlinear-dflp-ccd-ct1} \\
        & x^{t}_{ij} \leq a_{ij} y^{t}_{i} 
        && \substack{\forall i \in \mathcal{I}, \forall j \in \mathcal{J}, \\ \forall t \in \mathcal{T}} \label{eq:nonlinear-dflp-ccd-ct2} \\
        & \sum_{k \in \mathcal{I}} x^{t}_{kj} \geq a_{ij} y^{t}_{i}
        && \substack{\forall i \in \mathcal{I}, \forall j \in \mathcal{J}, \\ \forall t \in \mathcal{T}} \label{eq:nonlinear-dflp-ccd-ct3} \\
        & a_{ij} y^{t}_{i} + \sum_{\substack{k \in \mathcal{I}: \\ i \succ_{j} k}} x^{t}_{kj} \leq a_{ij} &&  \substack{\forall i \in \mathcal{I}, \forall j \in \mathcal{J}, \\ \forall t \in \mathcal{T}} \label{eq:nonlinear-dflp-ccd-ct4} \\
        & u^{0}_{j} = 0
        && \forall j \in \mathcal{J}
        \label{eq:nonlinear-dflp-ccd-ct5} \\
        & c^{t}_{j} = u^{t-1}_{j} + d^{t}_{j}
        && \substack{\forall j \in \mathcal{J}, \\ \forall t \in \mathcal{T}} \label{eq:nonlinear-dflp-ccd-ct6} \\
        & u^{t}_{j} = c^{t}_{j} - \sum_{i \in \mathcal{I}} c^{t}_{j} x^{t}_{ij}
        && \substack{\forall j \in \mathcal{J}, \\ \forall t \in \mathcal{T}} \label{eq:nonlinear-dflp-ccd-ct7} \\
        & u^{t}_{j} \in \mathbb{R}^{+}
        && \substack{\forall j \in \mathcal{J}, \\ \forall t \in \mathcal{T}}^{S} \label{eq:nonlinear-dflp-ccd-dm1} \\
        & c^{t}_{j}  \in \mathbb{R}^{+}
        && \substack{\forall j \in \mathcal{J}, \\ \forall t \in \mathcal{T}} \label{eq:nonlinear-dflp-ccd-dm2} \\
        & x^{t}_{ij} \in \mathset{0,1}
        && \substack{\forall i \in \mathcal{I}, \forall j \in \mathcal{J}, \\ \forall t \in \mathcal{T}} \label{eq:nonlinear-dflp-ccd-dm3} \\
        & y^{t}_{i} \in \mathset{0,1}
        && \substack{\forall i \in \mathcal{I}, \\ \forall t \in \mathcal{T}}. \label{eq:nonlinear-dflp-ccd-dm4}
    \end{align}
\end{subequations}

Objective Function \eqref{eq:nonlinear-dflp-ccd-obj} and Constraints \eqref{eq:nonlinear-dflp-ccd-ct1} --\eqref{eq:nonlinear-dflp-ccd-ct4} have the same meaning as their counterparts in the DIF, whereas Constraints \eqref{eq:nonlinear-dflp-ccd-ct5}--\eqref{eq:nonlinear-dflp-ccd-ct7} ensure proper cumulative demand behaviour over time and Constraints \eqref{eq:nonlinear-dflp-ccd-dm1}--\eqref{eq:nonlinear-dflp-ccd-dm4} define feasible variable domains.
Although the SIF is nonlinear, we can easily linearize it through standard techniques (for more details, see Online Appendix~B).
Even though the linearization requires the use of big-$M$ constraints, preliminary results show that off-the-shelf solvers perform better on the linearized version rather than on the nonlinear one.
Therefore, we employ the linearized version of the SIF.

\subsection{Formulation Comparison}
\label{sub:properties}

Once a location policy $\boldsymbol{y}$ is chosen by the provider, the values of the remaining components within the DI and SI Formulations can be unambiguously computed, as stated in Proposition~\ref{prp:representation}.
We thus refer to feasible solutions of the \ourproblem{} solely by a location policy $\boldsymbol{y}$, and employ $z(\boldsymbol{y})$ to denote its total profit.

\begin{proposition}
    \label{prp:representation}
    Feasible solutions of the DIF and the SIF can be solely represented by a location policy $\boldsymbol{y}$, since the remaining components have fixed values computable in polynomial time.
\end{proposition}

These formulations are equivalent (\ie, they have the same space of feasible and optimal integer solutions) and only differ in how they represent the cumulative customer demand.
They are considerably different, however, when it comes to their continuous relaxation bounds.
Intuitively, the linearization of bilinear terms $c^{t}_{j} x^{t}_{ij}$ within the SIF requires the addition of big-$M$ constraints, which tend to yield loose continuous relaxation bounds.
In this sense, we can show that the DIF provides a tighter continuous relaxation than the SIF, as stated in Theorem~\ref{thm:tightness}.

\begin{theorem}
    \label{thm:tightness}
    The DIF provides a tighter continuous relaxation than the SIF.
\end{theorem}

\textbf{Proof.}
    See Online Appendix~C.

We highlight that the DIF not only provides tighter continuous relaxation bounds, but also allows the representation of more general customer demand behaviour.
In fact, it can represent any customer demand behaviour $G^{\ell{}t}_{ij}$ that depends solely on time periods $\ell{}$ and $t$ for each customer $j$ and each location $i$.
Although some positive theoretical results may not hold for the problem variant with more general customer demand behaviour (namely, the approximation guarantees and the polynomial case further discussed in Section~\ref{sec:complexity}), we highlight that our solution methods, further presented in Section~\ref{sec:methods}, could be seamlessly adapted.
The main drawback of the DIF is the number of variables $x^{\ell{}t}_{ij}$, which grows quadratically with the size of the planning horizon $T$ and requires more computational resources to be solved.

\subsection{Problem Extensions}
\label{sec:ProblemExtensions}

We now briefly discuss some relevant extensions of the \ourproblem{} and its formulations in practice.

\subsubsection{Time-Dependent Parameters.}
The provider may need to account for time-dependent parameters.
For example, we may consider time-dependent customer rankings $\succ^{t}_{j}$ (\eg, when customer preferences vary with the season), time-dependent number of facilities $h^{t}$ (\eg, when some facilities become available or unavailable throughout time), and time-dependent rewards per captured demand unit $r^{t}_{i}$ (\eg, when profit margins fluctuate throughout the planning horizon).
We highlight that our solution methods can be easily adapted for time-dependent parameters, but they prevent us from holding onto some positive theoretical results (namely, the approximation guarantees and the polynomial case further discussed in Section~\ref{sec:complexity}).

\subsubsection{Penalties for Unmet Demand.}
The provider may need to account for penalties when the spawning demand $d^{t}_{j}$ of customer $j$ is not satisfied at time period $t$, but rather in future time periods $t + 1, \ldots, T$ (\eg, late-service fines charged by the government institutions or price discounts offered as customer service) or not at all.
In this sense, we remark that the DIF can readily account for several penalty structures by changing the formula to compute coefficients $G^{\ell{}t}_{ij}$. This includes complex nonlinear penalty functions, as long as they can be precomputed based on arguments $i$, $j$, $\ell{}$, $t$.
The SIF, however, requires structural changes of the formulation tailored to each specific penalty function (\eg, additional variables and constraints) to achieve the same goal.
For example, consider the case where the provider pays $p_{j} \in \mathbb{R}^{+}$ per demand of customer $j$ that is not served as soon as it appears.
In the DIF, we can compute coefficients $G^{\ell{}t}_{ij}$ as $G^{\ell{}t}_{ij} = r_{i} D^{\ell{}t}_{j} - \sum_{\substack{s \in \mathcal{T}: \\ s > \ell{} \\ \ell{} < t}} p_{j} d^{s}_{j}$ to account for penalties of customer $j$ between time periods $\ell{}$ and $t$.
In the SIF, we need to add the sum $- \sum_{t \in \mathcal{T}} \sum_{j \in \mathcal{J}} p_{j} d^{t}_{j} (1 - \sum_{i \in \mathcal{I}} x^{t}_{ij})$ to Objective Function~\eqref{eq:nonlinear-dflp-ccd-obj}.
Note that these changes do not affect the strength of the formulations, so Theorem~\ref{thm:tightness} still holds.
On the other hand, some positive theoretical results might not hold under penalties (namely, the approximation guarantees and the polynomial case further discussed in Section~\ref{sec:complexity}).

\section{Computational Complexity}
\label{sec:complexity}

In this section, we study the computational complexity of the \ourproblem{}.
We characterize special cases that are NP-hard with or without approximation guarantees, establishing the complexity of the general problem.
In addition, we discuss special cases solvable in polynomial time. 
For the sake of brevity, we present the proofs of Theorems~\ref{thm:nphard-multiple}--\ref{thm:polynomial} in Online Appendix~D.
We fix $\mathcal{Y} = \mathset{\sum_{i \in \mathcal{I}} y^{t}_{i} \leq h, \forall t \in \mathcal{T}}$ and $f(\boldsymbol{y}) = 0, \forall \boldsymbol{y} \in \mathcal{Y}$.

\subsection{General Case}
\label{sub:general}

We first consider the general case, where the provider might be able to locate multiple facilities over the planning horizon.
Let us formally define the decision version of the \ourproblem{} as follows:
\vspace{0.3cm}
\boxxx{
{\bf Decision version of the \ourproblem{}}: \\
{\sc instance}: Finite sets $\mathcal{T}=\{1, \ldots,T\}$, $\mathcal{I}=\{1, \ldots,I\}$ and $\mathcal{J}=\{1, \ldots,J\}$, positive rational numbers $\{r_{i}\}_{i \in \mathcal{I}}$ and $\{d_{j}^t\}_{j \in \mathcal{J}, t\in \mathcal{T}}$, positive integer number $h$, rankings $\{\succ_{j}\}_{j \in \mathcal{J}}$, binary values $\{a_{ij}\}_{i\in \mathcal{I}, j \in \mathcal{J}}$ and a positive rational number $Z$.\\
{\sc question}: Is there a feasible location policy $\boldsymbol{y}$ with an objective value of at least $Z$?
}
\vspace{0.3cm}

The \ourproblem{} is NP-hard because it generalizes the classical single-period facility location problem, which is known to be NP-hard \citep[][]{cornuejols1983uncapicitated}, as stated in Theorem~\ref{thm:nphard-multiple}.

\begin{theorem}
    \label{thm:nphard-multiple}
    The decision version of the \ourproblem{} is NP-complete.
\end{theorem}

Theorem~\ref{thm:nphard-multiple} implies that the \ourproblem{} is not an easy problem computationally (unless P = NP).
Naturally, single-period problem variants cannot account for cumulative demand behaviour.
As such, these problem variants do not provide insights on how cumulative customer demand might influence the theoretical intractability, which, as thoroughly discussed, has been rather ignored in the literature.
In the following, we therefore focus on the case with a single facility, which is considered an ``easy'' problem and based on which we can isolate the impact of introducing cumulative customer demand.

\subsection{Single Facility}
\label{sub:single}

We now consider a special case where the provider can only install a single facility throughout the planning horizon (\ie, $h = 1$), which we refer to as 1-\ourproblem{}.
This special case is not only relevant in practice, particularly in applications where the provider can deploy only one facility over the planning horizon,  but also allows us to underline in more detail how the cumulative demand behaviour (versus a noncumulative demand behaviour) contributes to the theoretical intractability of the \ourproblem{}.

In the 1-\ourproblem{}, the provider only opens a single facility per time period, so we can fully describe customer behaviour with parameters $a_{ij}$ (\ie, customer $j$ either patronizes the single facility at location $i$ or not).
Let us formally define the decision version of the 1-\ourproblem{} as follows:

\vspace{0.3cm}
\boxxx{
{\bf Decision version of the 1-\ourproblem{}}: \\
{\sc instance}: Finite sets $\mathcal{T}=\{1, \ldots,T\}$, $\mathcal{I}=\{1, \ldots,I\}$ and $\mathcal{J}=\{1, \ldots,J\}$, positive rational numbers $\{r_{i}\}_{i \in \mathcal{I}}$ and $\{d_{j}^t\}_{j \in \mathcal{J}, t\in \mathcal{T}}$, binary values $\{a_{ij}\}_{i\in \mathcal{I}, j \in \mathcal{J}}$ and a positive rational number $Z$.\\
{\sc question}: Is there a feasible location policy $\boldsymbol{y}$ with an objective value of at least $Z$?
}
\vspace{0.3cm}

For the sake of comparison, if we ignore the assumption that customer demand is cumulative in the 1-\ourproblem{}, we obtain the noncumulative variant of the 1-\ourproblem{}, which we refer to as 1-DFLP:
\begin{subequations}
    \label{pgm:linear-1-dflp}
    \begin{align}
        \max_{\boldsymbol{x}, \boldsymbol{y}} \quad
        &  \sum_{t \in \mathcal{T}} \sum_{i \in \mathcal{I}}  \sum_{j \in \mathcal{J}} r_{i} d^{t}_{j} x^{t}_{ij}
        && \label{eq:linear-1-dflp-obj}\\
        \text{s.t.} \quad
        & \sum_{i \in \mathcal{I}} y^{t}_{i} \leq 1
        && \forall t \in \mathcal{T} \label{eq:linear-1-dflp-ct1} \\
        & x^{t}_{ij} \leq a_{ij} y^{t}_{i} 
        && \forall i \in \mathcal{I}, \forall j \in \mathcal{J}, \forall t \in \mathcal{T} \label{eq:linear-1-dflp-ct2} \\
        & \sum_{k \in \mathcal{I}} x^{t}_{kj} \geq a_{ij} y^{t}_{i}
        && \forall i \in \mathcal{I}, \forall j \in \mathcal{J}, \forall t \in \mathcal{T} \label{eq:linear-1-dflp-ct3} \\
        & a_{ij} y^{t}_{i} + \sum_{\substack{k \in \mathcal{I}: \\ i \succ_{j} k}} x^{t}_{kj} \leq a_{ij} &&  \forall i \in \mathcal{I}, \forall j \in \mathcal{J}, \forall t \in \mathcal{T} \label{eq:linear-1-dflp-ct4} \\
        & x^{t}_{ij} \in \mathset{0,1}
        && \forall i \in \mathcal{I}, \forall j \in \mathcal{J}, \forall t \in \mathcal{T} \label{eq:linear-1-dflp-dm1} \\
        & y^{t}_{i} \in \mathset{0,1}
        && \forall i \in \mathcal{I}, \forall t \in \mathcal{T}. \label{eq:linear-1-dflp-dm2}
    \end{align}
\end{subequations}

Objective Function \eqref{eq:linear-1-dflp-obj} and Constraints~\eqref{eq:linear-1-dflp-ct1}--\eqref{eq:linear-1-dflp-dm2} have the same meaning as their counterparts in the SIF.
Note that we can easily build the optimal location policy $\boldsymbol{y}$ for the 1-DFLP by selecting the location with the largest marginal reward at each time period, so the 1-DFLP is polynomially solvable.
We study the theoretical intractability of the 1-\ourproblem{} to understand which problem characteristics, when interacting with cumulative customer demand, may then render it NP-hard.
To better categorize 1-\ourproblem{} instances in the upcoming discussion, we define the following instance descriptors:

\begin{definition}[Loyal and flexible customers] \label{def:customers}
    A 1-\ourproblem{} instance is said to have loyal customers if every customer is willing to attend only one location (\ie, $\sum_{i \in \mathcal{I}} a_{ij} = 1, \forall j \in \mathcal{J}$). A 1-\ourproblem{} instance is said to have flexible customers if it does not have loyal customers.
\end{definition}

\begin{definition}[Identical and different rewards] \label{def:identical-rewards}
    A 1-\ourproblem{} instance is said to have identical rewards if they are the same throughout the planning horizon (\ie, $r_{i}  = R, \forall i \in \mathcal{I},  R \in \mathbb{R^+}$).  A 1-\ourproblem{} instance is said to have different rewards if it does not have identical rewards.
\end{definition}

We show that the 1-\ourproblem{} is NP-hard through a reduction from the Set Packing Problem (SPP), which is known to be NP-hard \citep{karp1972reducibility} and cannot be approximated within a constant factor \citep{hazanComplexityApproximatingKset2006}.

\begin{theorem}
    \label{thm:nphard-single}
    The decision version of the 1-\ourproblem{} is NP-complete, and the 1-\ourproblem{} cannot be approximated within a factor $T^{1-\alpha}$ for any $\alpha > 0$, unless P = NP. 
\end{theorem}

The proof of Theorem \ref{thm:nphard-single} heavily relies on parameters $a_{ij}$ and rewards $r_{i}$ to obtain the NP-hardness results by assuming flexible customers and different rewards.
Note, however, that some instances may have loyal customers or identical rewards, and the 1-\ourproblem{} may become theoretically tractable for these special cases.
In this sense, we first investigate the 1-\ourproblem{} with identical rewards, and show that it remains NP-hard through a reduction from the Satisfiability Problem with exactly three variables per clause (3SAT), which is known to be strongly NP-hard \citep{karp1972reducibility, garey1979computers}.

\begin{theorem}
    \label{thm:nphard-identical}
    The decision version of the 1-\ourproblem{} with identical rewards remains NP-complete, and the 1-\ourproblem{} is strongly NP-hard.
\end{theorem}

Moreover, we can show that the 1-\ourproblem{} with identical rewards is approximable through a greedy algorithm that builds a location policy in reverse order.
More specifically, the iteration linked to time period $t$ finds which locations to open at said time period, assuming that no locations have been installed from time period $1$ to $t-1$, and fixes them in the solution.
Note that this heuristic takes the cumulative demand behaviour into consideration, but neglects effects of location decisions for earlier time periods.
We refer to this algorithm as the Backward Greedy Heuristic (BGH), and present its pseudocode in Algorithm~\ref{alg:bcw-greedy}.
We highlight that solving the formulation at hand for each time period can be done by inspection for the 1-\ourproblem{}, so this algorithm runs in polynomial time.
Note also that Algorithm~\ref{alg:bcw-greedy} is a heuristic for the \ourproblem{}.

\begin{algorithm}
    \caption{Backward Greedy Heuristic.}
    \label{alg:bcw-greedy}
    \begin{algorithmic}
        \REQUIRE Double-Index Formulation DIF, $\mathcal{I}$, $\mathcal{T} = \mathset{1, \ldots, T}$
        \STATE Add constraints $y^{t}_{i} = 0, \forall i \in \mathcal{I}, \forall t \in \mathcal{T}$ to the DIF
        \FORALL{$t = T, \ldots, 1$}
            \STATE Remove constraints $y^{t}_{i} = 0, \forall i \in \mathcal{I}$ from the DIF
            \STATE Solve the DIF to find optimal solution $\boldsymbol{y}^{\star}$
            \STATE Add constraints $y^{t}_{i} = {y^{t}_{i}}^{\star}, \forall i \in \mathcal{I}$ to the DIF
        \ENDFOR
        \STATE Solve the DIF to find optimal solution $\boldsymbol{y}^{\star}$
        \STATE \RETURN Location policy $\boldsymbol{y}^{\star}$.
    \end{algorithmic}
\end{algorithm}

\begin{theorem}
    \label{thm:approximation}
Algorithm~\ref{alg:bcw-greedy} is a 2-approximation algorithm for the 1-\ourproblem{} with identical rewards.
\end{theorem}

Theorems~\ref{thm:nphard-single}--\ref{thm:approximation} show that having different rewards seems to \textit{strengthen} the NP-hardness of the 1-\ourproblem{} (\ie, it turns a problem with potential approximation guarantees into a problem without them).
We now turn to the special case with loyal customers.
Surprisingly, we are able to show that this special case is solvable in polynomial time without further assumptions about rewards.

\begin{theorem}
    \label{thm:polynomial}
    The 1-\ourproblem{} with loyal customers is polynomially solvable.
\end{theorem}

Theorem~\ref{thm:polynomial}, along with Theorem~\ref{thm:nphard-identical}, implies that cumulative customer demand by itself does not generate the NP-hardness, but rather its intrinsic interaction with customer preferences over the planning horizon.
We highlight that the strongly NP-hard nature of the 1-\ourproblem{} (see Theorem~\ref{thm:nphard-identical}), as well as the inapproximability result (see Theorem~\ref{thm:nphard-single}), extend to the \ourproblem{}, since the former is a special case of the latter.
As previously mentioned, the approximation guarantees (see Theorem~\ref{thm:approximation}) as well as the polynomial case (see Theorem~\ref{thm:polynomial}) do not necessarily hold for problem extensions discussed in Section~\ref{sec:ProblemExtensions}.

\section{Solution Methods}
\label{sec:methods}

We now propose several solutions methods for the \ourproblem{}. 
We first present an exact Benders Decomposition \citep[][]{bendersPartitioningProceduresSolving1962, rahmanianiBendersDecompositionAlgorithm2017} of the DIF in Section~\ref{sub:benders}, and propose an analytical procedure to compute optimality cuts for the 1-\ourproblem{} (\ie, the special case studied in Section~\ref{sec:complexity}) in Section \ref{sub:analytical}.
We then provide a myopic heuristic to derive reasonable solutions when the provider ignores (or, equivalently, is unaware of) cumulative customer demand in Section \ref{sub:heuristics}.

\subsection{Benders Decomposition}
\label{sub:benders}

We can overcome the main drawback of the DIF, which is the large number of variables $x^{\ell{}t}_{ij}$, through an exact Benders Decomposition.
We keep variables $y^{t}_{i}$ in the master program and move variables $x^{\ell{}t}_{ij}$ to subproblems, one per customer $j$.
The master program of the DIF writes as follows:
\begin{subequations}
    \label{pgm:master-dflp-ccd}
    \begin{align}
        \max_{\boldsymbol{w}, \boldsymbol{y}} \quad
        &  \sum_{j \in \mathcal{J}} w_{j}
        \label{eq:master-dflp-ccd-obj}\\
        \text{s.t.} \quad
        & \boldsymbol{y} \in \mathcal{Y}
         \label{eq:master-dflp-ccd-ct1} \\     
        & w_{j} \leq  \sum_{t \in \mathcal{T}} \sum_{i \in  \mathcal{I}} \scriptstyle  a_{ij} \left[ ({\beta^{t}_{i}}^{\star} - {\zeta^{t}_{i}}^{\star} - {\delta^{t}_{i}}^{\star}) y^{t}_{i} + {\zeta^{t}_{i}}^{\star} \right] + {\theta^{0}}^{\star} \nonumber \\
        & \qquad \qquad\forall j \in \mathcal{J}, \forall (\boldsymbol{\beta}^{\star}, \boldsymbol{\delta}^{\star}, \boldsymbol{\zeta}^{\star}, \boldsymbol{\theta}^{\star}) \in \mathcal{O}_{j} \label{eq:master-dflp-ccd-ct2} \\
        & w_{j} \in \mathbb{R}^{+}
        \qquad \forall j \in \mathcal{J} \label{eq:master-dflp-ccd-dm1} \\
        & y^{t}_{i} \in \mathset{0,1}
        \qquad \forall i \in \mathcal{I}, \forall t \in \mathcal{T}, \label{eq:master-dflp-ccd-dm2}
    \end{align}
\end{subequations}
where $w_{j} \in \mathbb{R}^{+}$ is the objective value estimation of the subproblem of customer $j$ and $\mathcal{O}_{j}$ is the set of all optimality cuts for customer $j$ (\ie, the set of all optimal solutions of the dual subproblem of customer $j$).
The dual subproblem $w^{D}_{j} (\boldsymbol{y})$ of customer $j$ writes as follows:
\begin{subequations}
    \label{pgm:dual-dflp-ccd}
    \begin{align}
        w^{D}_{j} (\boldsymbol{y}): 
        \min_{\boldsymbol{\beta}, \boldsymbol{\delta}, \boldsymbol{\zeta}, \boldsymbol{\theta}} \quad
        &  \sum_{t \in \mathcal{T}} \sum_{i \in \mathcal{I}} \scriptstyle a_{ij} \left[ y^{t}_{i} (\beta^{t}_{i}  - \delta^{t}_{i} - \zeta^{t}_{i}) + \zeta^{t}_{i} \right] + \theta^{0} 
        \label{eq:dual-dflp-ccd-obj}\\
        \text{s.t.} \quad
        & \scriptstyle \beta^{t}_{i} - \sum_{k \in \mathcal{I}} \delta^{t}_{k} + \sum_{\substack{k \in \mathcal{I}: \\ k \succ_{j} i}} \zeta^{t}_{k} + \theta^{\ell{}} - \theta^{t}  \geq G^{\ell{}t}_{ij} 
        \nonumber \\ & \qquad \scriptstyle  \forall i \in \mathcal{I}, \forall \ell{} \in \mathcal{T}^{S}, \forall t \in \mathcal{T} : \ell{} < t \label{eq:dual-dflp-ccd-ct1} \\
        & \theta^{\ell{}} \geq G^{\ell{}(T+1)}_{ij}
        \quad \forall i \in \mathcal{I}, \forall \ell{} \in \mathcal{T}^{S} \label{eq:dual-dflp-ccd-ct2} \\
        & \beta^{t}_{i}, \delta^{t}_{i}, \zeta^{t}_{i} \in \mathbb{R}^{+}
        \quad \forall i \in \mathcal{I}, \forall t \in \mathcal{T}, \label{eq:dual-dflp-ccd-dm1}
    \end{align}
\end{subequations}
where dual variables $\beta^{t}_{i}$, $\delta^{t}_{i}$, $\zeta^{t}_{i}$, and $\theta^{t}$ are related to Constraints \eqref{eq:network-dflp-ccd-ct2}, \eqref{eq:network-dflp-ccd-ct3}, \eqref{eq:network-dflp-ccd-ct4}, and \eqref{eq:network-dflp-ccd-ct5}--\eqref{eq:network-dflp-ccd-ct6}, respectively.
Note that Constraints \eqref{eq:network-dflp-ccd-ct1} unveil an intrinsic dependency between time periods $\ell{}$ and $t$, thus preventing us from further decomposing subproblems by time period $t$.

We solve the Standard Benders Decomposition (SBD) in a branch-and-Benders-cut fashion \citep[][]{cordeauBendersDecompositionVery2019}.
More specifically, we start with a restricted set of optimality cuts $\tilde{\mathcal{O}_{j}}$ and, whenever we find a integer location policy $\boldsymbol{y}$ in the branch-and-bound tree, we solve the dual subproblem $w^{D}_{j} (\boldsymbol{y})$ and add the optimality cut $(\boldsymbol{\beta}^{\star}, \boldsymbol{\delta}^{\star}, \boldsymbol{\zeta}^{\star}, \boldsymbol{\theta}^{\star})$ to the restricted set $\tilde{\mathcal{O}_{j}}$.

\subsection{Optimality Cuts}
\label{sub:analytical}

The dual subproblem $w^{D}_{j} (\boldsymbol{y})$ tends to have multiple optimal solutions and, consequently, multiple optimality cuts for the same location policy $\boldsymbol{y}$.
In this sense, optimality cuts obtained by simply solving the dual subproblem $w^{D}_{j} (\boldsymbol{y})$ may be shallow or lack structure among themselves \citep[][]{magnantiAcceleratingBendersDecomposition1981}.
In addition, solving dual subproblems might become a bottleneck for larger number of customers, as observed in preliminary experiments and also in the literature \citep[][]{cordeauBendersDecompositionVery2019}.

One strategy to address these challenges is to employ an analytical procedure to compute optimality cuts, which can be done faster and guarantees that they have, at least, a similar structure among themselves.
Although the dual subproblem $w^{D}_{j} (\boldsymbol{y})$ cannot be easily solved by an analytical procedure in the general case, we devise such a method to compute optimality cuts for the 1-\ourproblem{}, presented in Algorithm~\ref{alg:analytical}.
We further refer to the branch-and-Benders-cut implementation with the analytical procedure to compute optimality cuts as the Analytical Benders Decomposition (ABD).
\begin{algorithm}
    \caption{Optimality Cuts for the 1-\ourproblem{}.}
    \label{alg:analytical}
    \begin{algorithmic}
        \REQUIRE $\mathcal{I}$, $\mathcal{T} = \mathset{1, \ldots, T}$, $\boldsymbol{y}$, $\boldsymbol{x}^{\star}$
        \STATE ${\theta^{\ell{}}}^{\star} \gets \max_{\substack{i \in \mathcal{I}, s \in \mathcal{T}^{S}, t \in \mathcal{T} :\\ s < t, {x^{st}_{i}}^{\star} = 1 \\ \ell{} < t, \ell{} \neq s}} \mathset{G^{\ell{}t}_{ij} - G^{st}_{ij} + {\theta^{s}}^{\star}, G^{\ell{}(T+1)}_{ij}}, \forall \ell{} \in \mathset{T, \ldots, 1, 0} : \sum_{i\in \mathcal{I}} a_{ij} y^{\ell{}}_{i} = 1$
        \STATE ${\theta^{\ell{}}}^{\star} \gets \max_{\substack{i \in \mathcal{I}, s \in \mathcal{T}^{S}, t \in \mathcal{T} :\\ s < t, {x^{st}_{i}}^{\star} = 1 \\ \ell{} < t, \ell{} \neq s}} \mathset{G^{\ell{}t}_{ij} - G^{st}_{ij} + {\theta^{s}}^{\star}, G^{\ell{}(T+1)}_{ij}}, \forall \ell{} \in \mathset{T, \ldots, 1, 0} : \sum_{i\in \mathcal{I}} a_{ij} y^{\ell{}}_{i} = 0$
        \STATE ${\lambda^{t}_{i}}^{\star} \gets \max_{\ell{} \in \mathcal{T} : \ell{} < t} \mathset{G^{\ell{}t}_{ij} - {\theta^{\ell{}}}^{\star} + {\theta^{t}}^{\star}}, \forall i \in \mathcal{I}, \forall t \in \mathcal{T}$
        \STATE \RETURN Optimality cut $w_{j} \leq \sum_{t \in \mathcal{T}} \sum_{i \in \mathcal{I}} a_{ij} {\lambda^{t}_{i}}^{\star} y^{t}_{o} + {\theta^{0}}^{\star}$.
    \end{algorithmic}
\end{algorithm}

\begin{theorem}
    \label{thm:analytical}
    Algorithm~\ref{alg:analytical} computes optimality cuts of the form $w_{j} \leq \sum_{t \in \mathcal{T}} \sum_{i \in \mathcal{I}} a_{ij} {\lambda^{t}_{i}}^{\star} y^{t}_{o} + {\theta^{0}}^{\star}$ for the 1-\ourproblem{}.
\end{theorem}

\textbf{Proof.}
    See Online Appendix~E.

\subsection{Myopic Heuristic}
\label{sub:heuristics}

The provider may explicitly or implicitly ignore cumulative customer demand when devising a location policy, thus behaving myopically.
In this sense, we present myopic heuristics to derive what seem to be natural solutions when cumulative demand behaviour is overlooked.
We further employ these heuristics to evaluate the economic benefit of modelling cumulative customer demand.

\subsubsection{DFLP-Based Heuristic.}
If we ignore the assumption that customer demand is cumulative, we can obtain the noncumulative variant of the \ourproblem{}, which we refer to as DFLP.
The DFLP has a structure similar to the 1-DFLP presented in Section~\ref{sub:single}, but with a general feasible set $\mathcal{Y}$ and a general cost function $f(\boldsymbol{y})$.
One can solve the DFLP to obtain a location policy $\boldsymbol{y}$ that completely ignores demand accumulation, and evaluate how it performs under the cumulative demand behaviour by computing the total profit $z(\boldsymbol{y})$.
Note that, if we additionally fix $\mathcal{Y} = \mathset{\sum_{i \in \mathcal{I}} y^{t}_{i} \leq h, \forall t \in \mathcal{T}}$ and $f(\boldsymbol{y}) = 0, \forall \boldsymbol{y} \in \mathcal{Y}$, the DFLP becomes separable by time periods. In this sense, this heuristic, referred to as the DFLP-Based Heuristic (DBH), chooses locations that capture large amounts of spawning demands at each time period. For example, if spawning demands are constant, a single location is chosen for the entire time horizon.

\subsubsection{Forward Greedy Heuristic.}
Another reasonable strategy is to perceive customer demand at each time period and then myopically choose the best location (\ie, the one that provides the highest marginal contribution) accordingly.
This heuristic, referred to as Forward Greedy Heuristic (FGH), takes the cumulative demand behaviour into consideration, but neglects future effects of current location decisions.
Although the Backward Greedy Heuristic presented in Section~\ref{sec:complexity} is similar to the Forward Greedy Heuristic in nature, we cannot trivially extend the approximation guarantees of the former to the latter (see Theorem~\ref{thm:approximation}).
Algorithm~\ref{alg:frw-greedy} presents the pseudocode of this heuristic.

\begin{algorithm}
    \caption{Forward Greedy Heuristic.}
    \label{alg:frw-greedy}
    \begin{algorithmic}
        \REQUIRE Double-Index Formulation DIF, $\mathcal{I}$, $\mathcal{T} = \mathset{1, \ldots, T}$
        \STATE Add constraints $y^{t}_{i} = 0, \forall i \in \mathcal{I}, \forall t \in \mathcal{T}$ to the DIF
        \FORALL{$t = 1, \ldots, T$}
            \STATE Remove constraints $y^{t}_{i} = 0, \forall i \in \mathcal{I}$ from the DIF
            \STATE Solve the DIF to find optimal solution $\boldsymbol{y}^{\star}$
            \STATE Add constraints $y^{t}_{i} = {y^{t}_{i}}^{\star}, \forall i \in \mathcal{I}$ to the DIF
        \ENDFOR
        \STATE Solve the DIF to find optimal solution $\boldsymbol{y}^{\star}$
        \STATE \RETURN Location policy $\boldsymbol{y}^{\star}$.
    \end{algorithmic}
\end{algorithm}

\section{Computational Experiments}
\label{sec:experiments}

In this section, we study the performance of the proposed solution methods and devise some managerial insights about the \ourproblem{}.
In Section~\ref{sub:instance-generation}, we describe the experimental setup, as well as the instance generation procedure.
Then, in Section~\ref{sub:computational-performance}, we evaluate the performance of the exact solution methods in terms of solution quality and computing times.
Lastly, in Section~\ref{sub:managerial-insights}, we investigate the impact of heuristic decisions for the provider, as well as the structure of optimal solutions for different instance attributes.

\subsection{Experimental Setup}

\label{sub:instance-generation}

We implement most of our solution methods in Python (version \texttt{3.10}), and solve the mixed-integer programming formulations with Gurobi (version \texttt{12.0}).
We write the analytical procedure to obtain optimality cuts for the 1-\ourproblem{}  in C instead of Python to ensure a fair comparison with Gurobi, which also runs on C.
Each solution method had a time limit of 1 hour and was limited to a single thread to avoid bias related to computational resources.
All jobs were processed on the Cedar cluster of the Digital Research Alliance of Canada with a maximum RAM of 30GB.
The code is available on \href{https://github.com/almeidawarley/1-dflp-ccd}{Github} under the MIT license.

Since the \ourproblem{} is a novel problem, benchmark instances are unavailable.
We therefore generate synthetic instances inspired by other related papers in the literature \citep[see, \eg,][]{marinMultiperiodStochasticCovering2018a}.
We create instances with $|\mathcal{T}| \in \mathset{5, 7, 9}$ time periods, $|\mathcal{I}| \in \mathset{50, 100, 150}$ candidate locations, $|\mathcal{J}| \in \mathset{1|\mathcal{I}|, 3|\mathcal{I}|, 5|\mathcal{I}|}$ targeted customers, and $h \in \mathset{1,3,5}$ facilities to understand the scalability of our solution methods.
We consider $C \in \mathset{0.05, 0.10}$ to build \textit{short} or \textit{long} rankings.
More precisely, for each customer $j \in \mathcal{J}$, we sample $\ceil{C \cdot I}$ candidate locations in order to build customer rankings $\succ_{j}$ and, consequently, parameters $a_{ij}$.
We consider \textit{identical} ($r_{i} = |\mathcal{I}|, \forall i \in \mathcal{I}$) and \textit{different} ($r_{i} = \ceil{\frac{|\mathcal{I}|}{\sum_{j \in \mathcal{J}}  a_{ij}}}, \forall i \in \mathcal{I}$) rewards per unit of 
captured demand.
Intuitively, the former describes applications where the reward is independent of location, whereas the latter describes applications where popular locations tend to have larger costs and, consequently, smaller rewards.
Lastly, we consider \textit{constant} ($d^{t}_{j} = 1, \forall j \in \mathcal{J}, \forall t \in \mathcal{T}$) and \textit{sparse} ($d^{t}_{j} \sim \mathset{0,1}, \forall j \in \mathcal{J}, \forall t \in \mathcal{T}$) spawning demands.
In simple terms, the former exemplifies scenarios where customers always have demand appearing throughout the planning horizon, whereas the latter exemplifies scenarios where customers may not have demand appearing at some time periods.
The combination of the aforementioned parameters yields a benchmark with $3^4 \cdot 2^3 = 648$ instances.

\subsection{Computational Performance}
\label{sub:computational-performance}

We first evaluate the performance of the three exact methods (\ie, DIF, SIF, and SBD) when it comes to finding an optimal solution (and proving its optimality) within the time limit of 1 hour.
Table~\ref{tab:optimality-summary} presents the percentage of the benchmark solved to optimality by these exact methods, grouped by dimensional instance attributes (\ie, number of locations $|\mathcal{I}|$, customers $|\mathcal{J}|$, and time periods $|\mathcal{T}|$).
As the instance size increases, these exact methods face greater difficulty in finding the optimal solution within the time limit, as they take longer to explore larger sets of feasible solutions.
Although the SBD struggles more than the DIF and the SIF to solve some small instances to optimality (\eg, instances with $|\mathcal{I}| = 50$, $|\mathcal{J}| = 1|\mathcal{I}|$, $|\mathcal{T}| = 9$), it handles larger instances much better overall (\eg, instances with $|\mathcal{I}| = 150$, $|\mathcal{J}| = 1|\mathcal{I}|$, $|\mathcal{T}| = 9$).
These results indicate a dominance relationship among the three exact methods in practice, where the SBD comes first, followed by the DIF then the SIF.
\begin{table}[!ht]
    \centering
    \caption{Percentage of the benchmark solved to optimality by the exact methods, grouped by dimensional instance attributes (\ie, number of locations $|\mathcal{I}|$, customers $|\mathcal{J}|$, and time periods $|\mathcal{T}|$).}
    \begin{tabular}{ccccc|ccccc|ccccc}
    \toprule
        \multicolumn{5}{c}{$|\mathcal{I}| = 50$} & \multicolumn{5}{|c|}{$|\mathcal{I}| = 100$} & \multicolumn{5}{c}{$|\mathcal{I}| = 150$}  \\ \midrule
        $|\mathcal{J}|$ & $|\mathcal{T}|$ & SIF & DIF & SBD & $|\mathcal{J}|$ & $|\mathcal{T}|$ & SIF & DIF & SBD & $|\mathcal{J}|$ & $|\mathcal{T}|$ & SIF & DIF & SBD \\ \midrule
        \multirow{3}*{$1|\mathcal{I}|$} & $5$ & $100\%$ & $100\%$ & $100 \%$ & \multirow{3}*{$1|\mathcal{I}|$} & $5$ & $88\%$ & $92\%$ & $92 \%$ & \multirow{3}*{$1|\mathcal{I}|$} & $5$ & $54\%$ & $58\%$ & $83\%$  \\
        ~ & $7$ & $96\%$ & $100\%$ & $92 \%$ & ~ & $7$ & $75\%$ & $88\%$ & $79 \%$ & ~ & $7$ & $42\%$ & $42\%$ & $67\%$  \\ 
        ~ & $9$ & $92\%$ & $100\%$ & $88 \%$ & ~ & $9$ & $62\%$ & $71\%$ & $75 \%$ & ~ & $9$ & $21\%$ & $17\%$ & $67\%$  \\ \midrule
        \multirow{3}*{$3|\mathcal{I}|$} & $5$ & $83\%$ & $88\%$ & $88 \%$ & \multirow{3}*{$3|\mathcal{I}|$} & $5$ & $38\%$ & $42\%$ & $50 \%$ & \multirow{3}*{$3|\mathcal{I}|$} & $5$ & $33\%$ & $29\%$ & $33\%$  \\
        ~ & $7$ & $62\%$ & $79\%$ & $79 \%$ & ~ & $7$ & $25\%$ & $25\%$ & $38 \%$ & ~ & $7$ & $17\%$ & $17\%$ & $29\%$  \\ 
        ~ & $9$ & $54\%$ & $62\%$ & $67 \%$ & ~ & $9$ & $17\%$ & $17\%$ & $29 \%$ & ~ & $9$ & $12\%$ & $8\%$ & $17\%$  \\ \midrule
        \multirow{3}*{$5|\mathcal{I}|$} & $5$ & $50\%$ & $62\%$ & $75 \%$ & \multirow{3}*{$5|\mathcal{I}|$} & $5$ & $33\%$ & $33\%$ & $38 \%$ & \multirow{3}*{$5|\mathcal{I}|$} & $5$ & $33\%$ & $21\%$ & $33\%$  \\ 
        ~ & $7$ & $33\%$ & $33\%$ & $50 \%$ & ~ & $7$ & $21\%$ & $21\%$ & $33 \%$ & ~ & $7$ & $12\%$ & $17\%$ & $25\%$  \\ 
        ~ & $9$ & $33\%$ & $25\%$ & $50 \%$ & ~ & $9$ & $17\%$ & $8\%$ & $21 \%$ & ~ & $9$ & $0\%$ & $0\%$ & $12\%$  \\ \bottomrule
    \end{tabular}
    \label{tab:optimality-summary}
\end{table}

We now investigate whether such a dominance relationship holds true in terms of solution quality and computing times.
To this end, we compute the \textit{objective ratio} of each exact method as $\frac{Z^{b}}{Z^{\prime}}$, where $Z^{b}$ is the highest objective value found among the exact methods and $Z^{\prime}$ is the objective value obtained by the exact method at hand.
Similarly, we compute the \textit{runtime ratio} of each exact method as $\frac{\Delta^{\prime}}{\Delta^{b}}$, where $\Delta^{b}$ is the lowest computing time among the exact methods and $\Delta^{\prime}$ is the computing time taken by the exact method at hand.
Small (respectively, large) objective ratios indicate that the exact method at hand finds a solution with an objective value closer to (respectively, farther from) the best one.
Similarly, small (respectively, large) runtime ratios indicate that the exact method at hand has a computing time closer to (respectively, farther from) the fastest one.
Figure~\ref{fig:performance-profile} presents the performance profile of the three exact methods, where the $y$ axis presents the number of instances with a ratio smaller than or equal to the reference value on the $x$ axis.
We can see that the SBD provides the best objective value and yields the fastest computing time when compared to the DIF and the SIF, validating the premise on the dominance relationship among them in practice.

\begin{figure}
    \begin{minipage}{0.99\textwidth}
         \centering
         \caption{Performance profiles of the three exact methods, where the $y$ axis presents the number of instances with a ratio smaller than or  equal to the reference value on the $x$ axis.}
         \vspace{0.3cm}
         \begin{subfigure}[b]{0.49\textwidth}
            \centering
            \caption{Objective ratio.}
            \begin{tikzpicture}[scale=.7, every node/.style={scale=.7}]
\draw[line width=0.5mm,thick,->] (0,0) -- (10.5,0);
\draw[line width=0.5mm,thick,->] (0,0) -- (0,5.5);
\draw (9.5,0.5) node[anchor=mid] {objective ratio};
\draw (0,6) node[anchor=mid] {\# instances};
\draw (0,-0.5) node[anchor=mid] {$1.00$};
\draw (1,-0.5) node[anchor=mid] {$1.05$};
\draw (2,-0.5) node[anchor=mid] {$1.10$};
\draw (3,-0.5) node[anchor=mid] {$1.15$};
\draw (4,-0.5) node[anchor=mid] {$1.20$};
\draw (5,-0.5) node[anchor=mid] {$1.25$};
\draw (6,-0.5) node[anchor=mid] {$1.30$};
\draw (7,-0.5) node[anchor=mid] {$1.35$};
\draw (8,-0.5) node[anchor=mid] {$1.40$};
\draw (9,-0.5) node[anchor=mid] {$1.45$};
\draw (10,-0.5) node[anchor=mid] {$1.50$};
\draw (-0.5,0) node[anchor=mid] {$259$};
\draw (-0.5,1) node[anchor=mid] {$337$};
\draw (-0.5,2) node[anchor=mid] {$415$};
\draw (-0.5,3) node[anchor=mid] {$492$};
\draw (-0.5,4) node[anchor=mid] {$570$};
\draw (-0.5,5) node[anchor=mid] {$648$};
\draw[line width=0.5mm,line width=0.5mm,blue,dashdotted] (0.00,3.17)--(1.00,4.33);\draw[line width=0.5mm,line width=0.5mm,blue,dashdotted] (1.00,4.33)--(2.00,4.58);\draw[line width=0.5mm,line width=0.5mm,blue,dashdotted] (2.00,4.58)--(3.00,4.83);\draw[line width=0.5mm,line width=0.5mm,blue,dashdotted] (3.00,4.83)--(4.00,4.83);\draw[line width=0.5mm,line width=0.5mm,blue,dashdotted] (4.00,4.83)--(5.00,4.92);\draw[line width=0.5mm,line width=0.5mm,blue,dashdotted] (5.00,4.92)--(6.00,4.92);\draw[line width=0.5mm,line width=0.5mm,blue,dashdotted] (6.00,4.92)--(7.00,4.92);\draw[line width=0.5mm,line width=0.5mm,blue,dashdotted] (7.00,4.92)--(8.00,4.92);\draw[line width=0.5mm,line width=0.5mm,blue,dashdotted] (8.00,4.92)--(9.00,4.92);\draw[line width=0.5mm,line width=0.5mm,blue,dashdotted] (9.00,4.92)--(10.00,4.92);\draw[line width=0.5mm,line width=0.5mm,magenta,dashed] (0.00,1.08)--(1.00,3.75);\draw[line width=0.5mm,line width=0.5mm,magenta,dashed] (1.00,3.75)--(2.00,4.00);\draw[line width=0.5mm,line width=0.5mm,magenta,dashed] (2.00,4.00)--(3.00,4.08);\draw[line width=0.5mm,line width=0.5mm,magenta,dashed] (3.00,4.08)--(4.00,4.17);\draw[line width=0.5mm,line width=0.5mm,magenta,dashed] (4.00,4.17)--(5.00,4.17);\draw[line width=0.5mm,line width=0.5mm,magenta,dashed] (5.00,4.17)--(6.00,4.25);\draw[line width=0.5mm,line width=0.5mm,magenta,dashed] (6.00,4.25)--(7.00,4.25);\draw[line width=0.5mm,line width=0.5mm,magenta,dashed] (7.00,4.25)--(8.00,4.25);\draw[line width=0.5mm,line width=0.5mm,magenta,dashed] (8.00,4.25)--(9.00,4.25);\draw[line width=0.5mm,line width=0.5mm,magenta,dashed] (9.00,4.25)--(10.00,4.25);\draw[line width=0.5mm,line width=0.5mm,teal,dotted] (0.00,2.25)--(1.00,4.33);\draw[line width=0.5mm,line width=0.5mm,teal,dotted] (1.00,4.33)--(2.00,4.58);\draw[line width=0.5mm,line width=0.5mm,teal,dotted] (2.00,4.58)--(3.00,4.58);\draw[line width=0.5mm,line width=0.5mm,teal,dotted] (3.00,4.58)--(4.00,4.58);\draw[line width=0.5mm,line width=0.5mm,teal,dotted] (4.00,4.58)--(5.00,4.58);\draw[line width=0.5mm,line width=0.5mm,teal,dotted] (5.00,4.58)--(6.00,4.58);\draw[line width=0.5mm,line width=0.5mm,teal,dotted] (6.00,4.58)--(7.00,4.58);\draw[line width=0.5mm,line width=0.5mm,teal,dotted] (7.00,4.58)--(8.00,4.58);\draw[line width=0.5mm,line width=0.5mm,teal,dotted] (8.00,4.58)--(9.00,4.58);\draw[line width=0.5mm,line width=0.5mm,teal,dotted] (9.00,4.58)--(10.00,4.58);
\draw[line width=0.5mm, blue, dashdotted] (8.5, 1.00)--(9.0, 1.00);
\draw[line width=0.5mm, blue] (9.0, 1.00) node[anchor=west] {SBD};
\draw[line width=0.5mm, magenta, dashed] (8.5, 1.50)--(9.0, 1.50);
\draw[line width=0.5mm, magenta] (9.0, 1.50) node[anchor=west] {SIF};
\draw[line width=0.5mm, teal, dotted] (8.5, 2.00)--(9.0, 2.00);
\draw[line width=0.5mm, teal] (9.0, 2.00) node[anchor=west] {DIF};
\end{tikzpicture}
            \label{fig:performance-objective}
         \end{subfigure}
        \hfill
         \begin{subfigure}[b]{0.49\textwidth}
            \centering
            \caption{Runtime ratio.}
            \begin{tikzpicture}[scale=.7, every node/.style={scale=.7}]
\draw[line width=0.5mm,thick,->] (0,0) -- (10.5,0);
\draw[line width=0.5mm,thick,->] (0,0) -- (0,5.5);
\draw (9.5,0.5) node[anchor=mid] {runtime ratio};
\draw (0,6) node[anchor=mid] {\# instances};
\draw (0,-0.5) node[anchor=mid] {$1.0$};
\draw (1,-0.5) node[anchor=mid] {$1.9$};
\draw (2,-0.5) node[anchor=mid] {$2.8$};
\draw (3,-0.5) node[anchor=mid] {$3.7$};
\draw (4,-0.5) node[anchor=mid] {$4.6$};
\draw (5,-0.5) node[anchor=mid] {$5.5$};
\draw (6,-0.5) node[anchor=mid] {$6.4$};
\draw (7,-0.5) node[anchor=mid] {$7.3$};
\draw (8,-0.5) node[anchor=mid] {$8.2$};
\draw (9,-0.5) node[anchor=mid] {$9.1$};
\draw (10,-0.5) node[anchor=mid] {$10.0$};
\draw (-0.5,0) node[anchor=mid] {$259$};
\draw (-0.5,1) node[anchor=mid] {$337$};
\draw (-0.5,2) node[anchor=mid] {$415$};
\draw (-0.5,3) node[anchor=mid] {$492$};
\draw (-0.5,4) node[anchor=mid] {$570$};
\draw (-0.5,5) node[anchor=mid] {$648$};
\draw[line width=0.5mm,line width=0.5mm,blue,dashdotted] (0.00,1.75)--(1.11,2.58);\draw[line width=0.5mm,line width=0.5mm,blue,dashdotted] (1.11,2.58)--(2.22,3.33);\draw[line width=0.5mm,line width=0.5mm,blue,dashdotted] (2.22,3.33)--(3.33,3.92);\draw[line width=0.5mm,line width=0.5mm,blue,dashdotted] (3.33,3.92)--(4.44,4.25);\draw[line width=0.5mm,line width=0.5mm,blue,dashdotted] (4.44,4.25)--(5.56,4.42);\draw[line width=0.5mm,line width=0.5mm,blue,dashdotted] (5.56,4.42)--(6.67,4.58);\draw[line width=0.5mm,line width=0.5mm,blue,dashdotted] (6.67,4.58)--(7.78,4.67);\draw[line width=0.5mm,line width=0.5mm,blue,dashdotted] (7.78,4.67)--(8.89,4.75);\draw[line width=0.5mm,line width=0.5mm,blue,dashdotted] (8.89,4.75)--(10.00,4.75);\draw[line width=0.5mm,line width=0.5mm,magenta,dashed] (0.00,0.25)--(1.11,0.75);\draw[line width=0.5mm,line width=0.5mm,magenta,dashed] (1.11,0.75)--(2.22,1.17);\draw[line width=0.5mm,line width=0.5mm,magenta,dashed] (2.22,1.17)--(3.33,1.50);\draw[line width=0.5mm,line width=0.5mm,magenta,dashed] (3.33,1.50)--(4.44,1.83);\draw[line width=0.5mm,line width=0.5mm,magenta,dashed] (4.44,1.83)--(5.56,2.00);\draw[line width=0.5mm,line width=0.5mm,magenta,dashed] (5.56,2.00)--(6.67,2.17);\draw[line width=0.5mm,line width=0.5mm,magenta,dashed] (6.67,2.17)--(7.78,2.50);\draw[line width=0.5mm,line width=0.5mm,magenta,dashed] (7.78,2.50)--(8.89,2.75);\draw[line width=0.5mm,line width=0.5mm,magenta,dashed] (8.89,2.75)--(10.00,2.92);\draw[line width=0.5mm,line width=0.5mm,teal,dotted] (0.00,0.67)--(1.11,1.50);\draw[line width=0.5mm,line width=0.5mm,teal,dotted] (1.11,1.50)--(2.22,1.92);\draw[line width=0.5mm,line width=0.5mm,teal,dotted] (2.22,1.92)--(3.33,2.17);\draw[line width=0.5mm,line width=0.5mm,teal,dotted] (3.33,2.17)--(4.44,2.50);\draw[line width=0.5mm,line width=0.5mm,teal,dotted] (4.44,2.50)--(5.56,2.67);\draw[line width=0.5mm,line width=0.5mm,teal,dotted] (5.56,2.67)--(6.67,2.75);\draw[line width=0.5mm,line width=0.5mm,teal,dotted] (6.67,2.75)--(7.78,2.92);\draw[line width=0.5mm,line width=0.5mm,teal,dotted] (7.78,2.92)--(8.89,3.08);\draw[line width=0.5mm,line width=0.5mm,teal,dotted] (8.89,3.08)--(10.00,3.33);
\draw[line width=0.5mm, blue, dashdotted] (8.5, 1.00)--(9.0, 1.00);
\draw[line width=0.5mm, blue] (9.0, 1.00) node[anchor=west] {SBD};
\draw[line width=0.5mm, magenta, dashed] (8.5, 1.50)--(9.0, 1.50);
\draw[line width=0.5mm, magenta] (9.0, 1.50) node[anchor=west] {SIF};
\draw[line width=0.5mm, teal, dotted] (8.5, 2.00)--(9.0, 2.00);
\draw[line width=0.5mm, teal] (9.0, 2.00) node[anchor=west] {DIF};
\end{tikzpicture}
            \label{fig:performance-runtime}
         \end{subfigure}
        \label{fig:performance-profile}   
    \end{minipage}
\end{figure}

\begin{figure}
    \begin{minipage}{.99\textwidth}
        \centering
        \caption{Comparison between DIF and SIF, where the $y$ axis presents the number of instances with a value (solution time or optimality gap) smaller than or  equal to the reference value on the $x$ axis.}
        \label{fig:c1}
        \vspace{0.3cm}
         \begin{subfigure}[b]{0.49\textwidth}
            \centering
            \caption{Instances solved to optimality by both methods.}
            \begin{tikzpicture}[scale=.7, every node/.style={scale=.7}]
\draw[line width=0.5mm,thick,->] (0,0) -- (10.5,0);
\draw[line width=0.5mm,thick,->] (0,0) -- (0,5.5);
\draw (9.5,0.5) node[anchor=mid] {runtime (min)};
\draw (0,6) node[anchor=mid] {\# instances};
\draw (0,-0.5) node[anchor=mid] {$0$};
\draw (1,-0.5) node[anchor=mid] {$6$};
\draw (2,-0.5) node[anchor=mid] {$12$};
\draw (3,-0.5) node[anchor=mid] {$18$};
\draw (4,-0.5) node[anchor=mid] {$24$};
\draw (5,-0.5) node[anchor=mid] {$30$};
\draw (6,-0.5) node[anchor=mid] {$36$};
\draw (7,-0.5) node[anchor=mid] {$42$};
\draw (8,-0.5) node[anchor=mid] {$48$};
\draw (9,-0.5) node[anchor=mid] {$54$};
\draw (10,-0.5) node[anchor=mid] {$60$};
\draw (-0.5,0) node[anchor=mid] {$0$};
\draw (-0.5,1) node[anchor=mid] {$55$};
\draw (-0.5,2) node[anchor=mid] {$111$};
\draw (-0.5,3) node[anchor=mid] {$166$};
\draw (-0.5,4) node[anchor=mid] {$222$};
\draw (-0.5,5) node[anchor=mid] {$277$};
\draw[line width=0.5mm,line width=0.5mm,magenta,dashed] (0.00,0.00)--(0.50,3.35);\draw[line width=0.5mm,line width=0.5mm,magenta,dashed] (0.50,3.35)--(1.00,3.75);\draw[line width=0.5mm,line width=0.5mm,magenta,dashed] (1.00,3.75)--(1.50,3.95);\draw[line width=0.5mm,line width=0.5mm,magenta,dashed] (1.50,3.95)--(2.00,4.10);\draw[line width=0.5mm,line width=0.5mm,magenta,dashed] (2.00,4.10)--(2.50,4.15);\draw[line width=0.5mm,line width=0.5mm,magenta,dashed] (2.50,4.15)--(3.00,4.25);\draw[line width=0.5mm,line width=0.5mm,magenta,dashed] (3.00,4.25)--(3.50,4.35);\draw[line width=0.5mm,line width=0.5mm,magenta,dashed] (3.50,4.35)--(4.00,4.50);\draw[line width=0.5mm,line width=0.5mm,magenta,dashed] (4.00,4.50)--(4.50,4.60);\draw[line width=0.5mm,line width=0.5mm,magenta,dashed] (4.50,4.60)--(5.00,4.60);\draw[line width=0.5mm,line width=0.5mm,magenta,dashed] (5.00,4.60)--(5.50,4.65);\draw[line width=0.5mm,line width=0.5mm,magenta,dashed] (5.50,4.65)--(6.00,4.65);\draw[line width=0.5mm,line width=0.5mm,magenta,dashed] (6.00,4.65)--(6.50,4.70);\draw[line width=0.5mm,line width=0.5mm,magenta,dashed] (6.50,4.70)--(7.00,4.70);\draw[line width=0.5mm,line width=0.5mm,magenta,dashed] (7.00,4.70)--(7.50,4.75);\draw[line width=0.5mm,line width=0.5mm,magenta,dashed] (7.50,4.75)--(8.00,4.80);\draw[line width=0.5mm,line width=0.5mm,magenta,dashed] (8.00,4.80)--(8.50,4.85);\draw[line width=0.5mm,line width=0.5mm,magenta,dashed] (8.50,4.85)--(9.00,4.85);\draw[line width=0.5mm,line width=0.5mm,magenta,dashed] (9.00,4.85)--(9.50,4.90);\draw[line width=0.5mm,line width=0.5mm,magenta,dashed] (9.50,4.90)--(10.00,5.00);\draw[line width=0.5mm,line width=0.5mm,teal,dotted] (0.00,0.00)--(0.50,3.60);\draw[line width=0.5mm,line width=0.5mm,teal,dotted] (0.50,3.60)--(1.00,4.00);\draw[line width=0.5mm,line width=0.5mm,teal,dotted] (1.00,4.00)--(1.50,4.20);\draw[line width=0.5mm,line width=0.5mm,teal,dotted] (1.50,4.20)--(2.00,4.35);\draw[line width=0.5mm,line width=0.5mm,teal,dotted] (2.00,4.35)--(2.50,4.50);\draw[line width=0.5mm,line width=0.5mm,teal,dotted] (2.50,4.50)--(3.00,4.60);\draw[line width=0.5mm,line width=0.5mm,teal,dotted] (3.00,4.60)--(3.50,4.65);\draw[line width=0.5mm,line width=0.5mm,teal,dotted] (3.50,4.65)--(4.00,4.70);\draw[line width=0.5mm,line width=0.5mm,teal,dotted] (4.00,4.70)--(4.50,4.80);\draw[line width=0.5mm,line width=0.5mm,teal,dotted] (4.50,4.80)--(5.00,4.80);\draw[line width=0.5mm,line width=0.5mm,teal,dotted] (5.00,4.80)--(5.50,4.85);\draw[line width=0.5mm,line width=0.5mm,teal,dotted] (5.50,4.85)--(6.00,4.85);\draw[line width=0.5mm,line width=0.5mm,teal,dotted] (6.00,4.85)--(6.50,4.85);\draw[line width=0.5mm,line width=0.5mm,teal,dotted] (6.50,4.85)--(7.00,4.90);\draw[line width=0.5mm,line width=0.5mm,teal,dotted] (7.00,4.90)--(7.50,4.90);\draw[line width=0.5mm,line width=0.5mm,teal,dotted] (7.50,4.90)--(8.00,4.90);\draw[line width=0.5mm,line width=0.5mm,teal,dotted] (8.00,4.90)--(8.50,4.90);\draw[line width=0.5mm,line width=0.5mm,teal,dotted] (8.50,4.90)--(9.00,4.95);\draw[line width=0.5mm,line width=0.5mm,teal,dotted] (9.00,4.95)--(9.50,4.95);\draw[line width=0.5mm,line width=0.5mm,teal,dotted] (9.50,4.95)--(10.00,5.00);
\draw[line width=0.5mm, magenta, dashed] (8.5, 3.00)--(9.0, 3.00);
\draw[line width=0.5mm, magenta] (9.0, 3.00) node[anchor=west] {SIF};
\draw[line width=0.5mm, teal, dotted] (8.5, 3.50)--(9.0, 3.50);
\draw[line width=0.5mm, teal] (9.0, 3.50) node[anchor=west] {DIF};
\end{tikzpicture}
        \label{fig:runtimes-c1}
         \end{subfigure}
        \hfill
         \begin{subfigure}[b]{0.49\textwidth}
            \centering
            \caption{Instances not solved to optimality by both methods.}
            \begin{tikzpicture}[scale=.7, every node/.style={scale=.7}]
\draw[line width=0.5mm,thick,->] (0,0) -- (10.5,0);
\draw[line width=0.5mm,thick,->] (0,0) -- (0,5.5);
\draw (9.5,0.5) node[anchor=mid] {optimality gap (\%)};
\draw (0,6) node[anchor=mid] {\# instances};
\draw (0,-0.5) node[anchor=mid] {$0$};
\draw (1,-0.5) node[anchor=mid] {$10$};
\draw (2,-0.5) node[anchor=mid] {$20$};
\draw (3,-0.5) node[anchor=mid] {$30$};
\draw (4,-0.5) node[anchor=mid] {$40$};
\draw (5,-0.5) node[anchor=mid] {$50$};
\draw (6,-0.5) node[anchor=mid] {$60$};
\draw (7,-0.5) node[anchor=mid] {$70$};
\draw (8,-0.5) node[anchor=mid] {$80$};
\draw (9,-0.5) node[anchor=mid] {$90$};
\draw (10,-0.5) node[anchor=mid] {$100$};
\draw (-0.5,0) node[anchor=mid] {$0$};
\draw (-0.5,1) node[anchor=mid] {$74$};
\draw (-0.5,2) node[anchor=mid] {$148$};
\draw (-0.5,3) node[anchor=mid] {$223$};
\draw (-0.5,4) node[anchor=mid] {$297$};
\draw (-0.5,5) node[anchor=mid] {$371$};
\draw[line width=0.5mm,line width=0.5mm,magenta,dashed] (0.00,0.05)--(0.50,1.65);\draw[line width=0.5mm,line width=0.5mm,magenta,dashed] (0.50,1.65)--(1.00,2.80);\draw[line width=0.5mm,line width=0.5mm,magenta,dashed] (1.00,2.80)--(1.50,3.65);\draw[line width=0.5mm,line width=0.5mm,magenta,dashed] (1.50,3.65)--(2.00,3.90);\draw[line width=0.5mm,line width=0.5mm,magenta,dashed] (2.00,3.90)--(2.50,4.00);\draw[line width=0.5mm,line width=0.5mm,magenta,dashed] (2.50,4.00)--(3.00,4.15);\draw[line width=0.5mm,line width=0.5mm,magenta,dashed] (3.00,4.15)--(3.50,4.15);\draw[line width=0.5mm,line width=0.5mm,magenta,dashed] (3.50,4.15)--(4.00,4.20);\draw[line width=0.5mm,line width=0.5mm,magenta,dashed] (4.00,4.20)--(4.50,4.20);\draw[line width=0.5mm,line width=0.5mm,magenta,dashed] (4.50,4.20)--(5.00,4.20);\draw[line width=0.5mm,line width=0.5mm,magenta,dashed] (5.00,4.20)--(5.50,4.20);\draw[line width=0.5mm,line width=0.5mm,magenta,dashed] (5.50,4.20)--(6.00,4.20);\draw[line width=0.5mm,line width=0.5mm,magenta,dashed] (6.00,4.20)--(6.50,4.25);\draw[line width=0.5mm,line width=0.5mm,magenta,dashed] (6.50,4.25)--(7.00,4.25);\draw[line width=0.5mm,line width=0.5mm,magenta,dashed] (7.00,4.25)--(7.50,4.25);\draw[line width=0.5mm,line width=0.5mm,magenta,dashed] (7.50,4.25)--(8.00,4.25);\draw[line width=0.5mm,line width=0.5mm,magenta,dashed] (8.00,4.25)--(8.50,4.25);\draw[line width=0.5mm,line width=0.5mm,magenta,dashed] (8.50,4.25)--(9.00,4.25);\draw[line width=0.5mm,line width=0.5mm,magenta,dashed] (9.00,4.25)--(9.50,4.25);\draw[line width=0.5mm,line width=0.5mm,magenta,dashed] (9.50,4.25)--(10.00,5.00);\draw[line width=0.5mm,line width=0.5mm,teal,dotted] (0.00,0.05)--(0.50,2.05);\draw[line width=0.5mm,line width=0.5mm,teal,dotted] (0.50,2.05)--(1.00,3.30);\draw[line width=0.5mm,line width=0.5mm,teal,dotted] (1.00,3.30)--(1.50,4.25);\draw[line width=0.5mm,line width=0.5mm,teal,dotted] (1.50,4.25)--(2.00,4.55);\draw[line width=0.5mm,line width=0.5mm,teal,dotted] (2.00,4.55)--(2.50,4.55);\draw[line width=0.5mm,line width=0.5mm,teal,dotted] (2.50,4.55)--(3.00,4.60);\draw[line width=0.5mm,line width=0.5mm,teal,dotted] (3.00,4.60)--(3.50,4.60);\draw[line width=0.5mm,line width=0.5mm,teal,dotted] (3.50,4.60)--(4.00,4.60);\draw[line width=0.5mm,line width=0.5mm,teal,dotted] (4.00,4.60)--(4.50,4.60);\draw[line width=0.5mm,line width=0.5mm,teal,dotted] (4.50,4.60)--(5.00,4.60);\draw[line width=0.5mm,line width=0.5mm,teal,dotted] (5.00,4.60)--(5.50,4.60);\draw[line width=0.5mm,line width=0.5mm,teal,dotted] (5.50,4.60)--(6.00,4.60);\draw[line width=0.5mm,line width=0.5mm,teal,dotted] (6.00,4.60)--(6.50,4.60);\draw[line width=0.5mm,line width=0.5mm,teal,dotted] (6.50,4.60)--(7.00,4.60);\draw[line width=0.5mm,line width=0.5mm,teal,dotted] (7.00,4.60)--(7.50,4.60);\draw[line width=0.5mm,line width=0.5mm,teal,dotted] (7.50,4.60)--(8.00,4.60);\draw[line width=0.5mm,line width=0.5mm,teal,dotted] (8.00,4.60)--(8.50,4.60);\draw[line width=0.5mm,line width=0.5mm,teal,dotted] (8.50,4.60)--(9.00,4.60);\draw[line width=0.5mm,line width=0.5mm,teal,dotted] (9.00,4.60)--(9.50,4.60);\draw[line width=0.5mm,line width=0.5mm,teal,dotted] (9.50,4.60)--(10.00,5.00);
\draw[line width=0.5mm, magenta, dashed] (8.5, 3.00)--(9.0, 3.00);
\draw[line width=0.5mm, magenta] (9.0, 3.00) node[anchor=west] {SIF};
\draw[line width=0.5mm, teal, dotted] (8.5, 3.50)--(9.0, 3.50);
\draw[line width=0.5mm, teal] (9.0, 3.50) node[anchor=west] {DIF};
\end{tikzpicture}
            \label{fig:optgaps-c1}
         \end{subfigure}        
    \end{minipage}
    \begin{minipage}{.99\textwidth}
        \centering
        \caption{Comparison between SBD and DIF, where the $y$ axis presents the number of instances with a value (solution time or optimality gap) smaller than or  equal to the reference value on the $x$ axis.}
        \label{fig:c2}
        \vspace{0.3cm}
         \begin{subfigure}[b]{0.49\textwidth}
            \centering
            \caption{Instances solved to optimality by both methods.}
            \begin{tikzpicture}[scale=.7, every node/.style={scale=.7}]
\draw[line width=0.5mm,thick,->] (0,0) -- (10.5,0);
\draw[line width=0.5mm,thick,->] (0,0) -- (0,5.5);
\draw (9.5,0.5) node[anchor=mid] {runtime (min)};
\draw (0,6) node[anchor=mid] {\# instances};
\draw (0,-0.5) node[anchor=mid] {$0$};
\draw (1,-0.5) node[anchor=mid] {$6$};
\draw (2,-0.5) node[anchor=mid] {$12$};
\draw (3,-0.5) node[anchor=mid] {$18$};
\draw (4,-0.5) node[anchor=mid] {$24$};
\draw (5,-0.5) node[anchor=mid] {$30$};
\draw (6,-0.5) node[anchor=mid] {$36$};
\draw (7,-0.5) node[anchor=mid] {$42$};
\draw (8,-0.5) node[anchor=mid] {$48$};
\draw (9,-0.5) node[anchor=mid] {$54$};
\draw (10,-0.5) node[anchor=mid] {$60$};
\draw (-0.5,0) node[anchor=mid] {$0$};
\draw (-0.5,1) node[anchor=mid] {$59$};
\draw (-0.5,2) node[anchor=mid] {$117$};
\draw (-0.5,3) node[anchor=mid] {$176$};
\draw (-0.5,4) node[anchor=mid] {$234$};
\draw (-0.5,5) node[anchor=mid] {$293$};
\draw[line width=0.5mm,line width=0.5mm,teal,dotted] (0.00,0.00)--(0.50,3.40);\draw[line width=0.5mm,line width=0.5mm,teal,dotted] (0.50,3.40)--(1.00,3.80);\draw[line width=0.5mm,line width=0.5mm,teal,dotted] (1.00,3.80)--(1.50,4.05);\draw[line width=0.5mm,line width=0.5mm,teal,dotted] (1.50,4.05)--(2.00,4.25);\draw[line width=0.5mm,line width=0.5mm,teal,dotted] (2.00,4.25)--(2.50,4.40);\draw[line width=0.5mm,line width=0.5mm,teal,dotted] (2.50,4.40)--(3.00,4.45);\draw[line width=0.5mm,line width=0.5mm,teal,dotted] (3.00,4.45)--(3.50,4.50);\draw[line width=0.5mm,line width=0.5mm,teal,dotted] (3.50,4.50)--(4.00,4.60);\draw[line width=0.5mm,line width=0.5mm,teal,dotted] (4.00,4.60)--(4.50,4.65);\draw[line width=0.5mm,line width=0.5mm,teal,dotted] (4.50,4.65)--(5.00,4.70);\draw[line width=0.5mm,line width=0.5mm,teal,dotted] (5.00,4.70)--(5.50,4.70);\draw[line width=0.5mm,line width=0.5mm,teal,dotted] (5.50,4.70)--(6.00,4.75);\draw[line width=0.5mm,line width=0.5mm,teal,dotted] (6.00,4.75)--(6.50,4.80);\draw[line width=0.5mm,line width=0.5mm,teal,dotted] (6.50,4.80)--(7.00,4.80);\draw[line width=0.5mm,line width=0.5mm,teal,dotted] (7.00,4.80)--(7.50,4.85);\draw[line width=0.5mm,line width=0.5mm,teal,dotted] (7.50,4.85)--(8.00,4.85);\draw[line width=0.5mm,line width=0.5mm,teal,dotted] (8.00,4.85)--(8.50,4.85);\draw[line width=0.5mm,line width=0.5mm,teal,dotted] (8.50,4.85)--(9.00,4.90);\draw[line width=0.5mm,line width=0.5mm,teal,dotted] (9.00,4.90)--(9.50,4.90);\draw[line width=0.5mm,line width=0.5mm,teal,dotted] (9.50,4.90)--(10.00,5.00);\draw[line width=0.5mm,line width=0.5mm,blue,dashdotted] (0.00,0.00)--(0.50,4.60);\draw[line width=0.5mm,line width=0.5mm,blue,dashdotted] (0.50,4.60)--(1.00,4.70);\draw[line width=0.5mm,line width=0.5mm,blue,dashdotted] (1.00,4.70)--(1.50,4.80);\draw[line width=0.5mm,line width=0.5mm,blue,dashdotted] (1.50,4.80)--(2.00,4.85);\draw[line width=0.5mm,line width=0.5mm,blue,dashdotted] (2.00,4.85)--(2.50,4.90);\draw[line width=0.5mm,line width=0.5mm,blue,dashdotted] (2.50,4.90)--(3.00,4.90);\draw[line width=0.5mm,line width=0.5mm,blue,dashdotted] (3.00,4.90)--(3.50,4.95);\draw[line width=0.5mm,line width=0.5mm,blue,dashdotted] (3.50,4.95)--(4.00,5.00);\draw[line width=0.5mm,line width=0.5mm,blue,dashdotted] (4.00,5.00)--(4.50,5.00);\draw[line width=0.5mm,line width=0.5mm,blue,dashdotted] (4.50,5.00)--(5.00,5.00);\draw[line width=0.5mm,line width=0.5mm,blue,dashdotted] (5.00,5.00)--(5.50,5.00);\draw[line width=0.5mm,line width=0.5mm,blue,dashdotted] (5.50,5.00)--(6.00,5.00);\draw[line width=0.5mm,line width=0.5mm,blue,dashdotted] (6.00,5.00)--(6.50,5.00);\draw[line width=0.5mm,line width=0.5mm,blue,dashdotted] (6.50,5.00)--(7.00,5.00);\draw[line width=0.5mm,line width=0.5mm,blue,dashdotted] (7.00,5.00)--(7.50,5.00);\draw[line width=0.5mm,line width=0.5mm,blue,dashdotted] (7.50,5.00)--(8.00,5.00);\draw[line width=0.5mm,line width=0.5mm,blue,dashdotted] (8.00,5.00)--(8.50,5.00);\draw[line width=0.5mm,line width=0.5mm,blue,dashdotted] (8.50,5.00)--(9.00,5.00);\draw[line width=0.5mm,line width=0.5mm,blue,dashdotted] (9.00,5.00)--(9.50,5.00);\draw[line width=0.5mm,line width=0.5mm,blue,dashdotted] (9.50,5.00)--(10.00,5.00);
\draw[line width=0.5mm, teal, dotted] (8.5, 3.00)--(9.0, 3.00);
\draw[line width=0.5mm, teal] (9.0, 3.00) node[anchor=west] {DIF};
\draw[line width=0.5mm, blue, dashdotted] (8.5, 3.50)--(9.0, 3.50);
\draw[line width=0.5mm, blue] (9.0, 3.50) node[anchor=west] {SBD};
\end{tikzpicture}
            \label{fig:runtimes-c2}
         \end{subfigure}
         \hfill
         \begin{subfigure}[b]{0.49\textwidth}
            \centering
            \caption{Instances not solved to optimality by both methods.}
            \begin{tikzpicture}[scale=.7, every node/.style={scale=.7}]
\draw[line width=0.5mm,thick,->] (0,0) -- (10.5,0);
\draw[line width=0.5mm,thick,->] (0,0) -- (0,5.5);
\draw (9.5,0.5) node[anchor=mid] {optimality gap (\%)};
\draw (0,6) node[anchor=mid] {\# instances};
\draw (0,-0.5) node[anchor=mid] {$0$};
\draw (1,-0.5) node[anchor=mid] {$10$};
\draw (2,-0.5) node[anchor=mid] {$20$};
\draw (3,-0.5) node[anchor=mid] {$30$};
\draw (4,-0.5) node[anchor=mid] {$40$};
\draw (5,-0.5) node[anchor=mid] {$50$};
\draw (6,-0.5) node[anchor=mid] {$60$};
\draw (7,-0.5) node[anchor=mid] {$70$};
\draw (8,-0.5) node[anchor=mid] {$80$};
\draw (9,-0.5) node[anchor=mid] {$90$};
\draw (10,-0.5) node[anchor=mid] {$100$};
\draw (-0.5,0) node[anchor=mid] {$0$};
\draw (-0.5,1) node[anchor=mid] {$71$};
\draw (-0.5,2) node[anchor=mid] {$142$};
\draw (-0.5,3) node[anchor=mid] {$213$};
\draw (-0.5,4) node[anchor=mid] {$284$};
\draw (-0.5,5) node[anchor=mid] {$355$};
\draw[line width=0.5mm,line width=0.5mm,teal,dotted] (0.00,0.00)--(0.50,1.90);\draw[line width=0.5mm,line width=0.5mm,teal,dotted] (0.50,1.90)--(1.00,3.20);\draw[line width=0.5mm,line width=0.5mm,teal,dotted] (1.00,3.20)--(1.50,4.20);\draw[line width=0.5mm,line width=0.5mm,teal,dotted] (1.50,4.20)--(2.00,4.50);\draw[line width=0.5mm,line width=0.5mm,teal,dotted] (2.00,4.50)--(2.50,4.55);\draw[line width=0.5mm,line width=0.5mm,teal,dotted] (2.50,4.55)--(3.00,4.55);\draw[line width=0.5mm,line width=0.5mm,teal,dotted] (3.00,4.55)--(3.50,4.55);\draw[line width=0.5mm,line width=0.5mm,teal,dotted] (3.50,4.55)--(4.00,4.55);\draw[line width=0.5mm,line width=0.5mm,teal,dotted] (4.00,4.55)--(4.50,4.55);\draw[line width=0.5mm,line width=0.5mm,teal,dotted] (4.50,4.55)--(5.00,4.55);\draw[line width=0.5mm,line width=0.5mm,teal,dotted] (5.00,4.55)--(5.50,4.55);\draw[line width=0.5mm,line width=0.5mm,teal,dotted] (5.50,4.55)--(6.00,4.55);\draw[line width=0.5mm,line width=0.5mm,teal,dotted] (6.00,4.55)--(6.50,4.55);\draw[line width=0.5mm,line width=0.5mm,teal,dotted] (6.50,4.55)--(7.00,4.55);\draw[line width=0.5mm,line width=0.5mm,teal,dotted] (7.00,4.55)--(7.50,4.55);\draw[line width=0.5mm,line width=0.5mm,teal,dotted] (7.50,4.55)--(8.00,4.55);\draw[line width=0.5mm,line width=0.5mm,teal,dotted] (8.00,4.55)--(8.50,4.55);\draw[line width=0.5mm,line width=0.5mm,teal,dotted] (8.50,4.55)--(9.00,4.55);\draw[line width=0.5mm,line width=0.5mm,teal,dotted] (9.00,4.55)--(9.50,4.55);\draw[line width=0.5mm,line width=0.5mm,teal,dotted] (9.50,4.55)--(10.00,5.00);\draw[line width=0.5mm,line width=0.5mm,blue,dashdotted] (0.00,0.05)--(0.50,2.40);\draw[line width=0.5mm,line width=0.5mm,blue,dashdotted] (0.50,2.40)--(1.00,3.55);\draw[line width=0.5mm,line width=0.5mm,blue,dashdotted] (1.00,3.55)--(1.50,4.10);\draw[line width=0.5mm,line width=0.5mm,blue,dashdotted] (1.50,4.10)--(2.00,4.55);\draw[line width=0.5mm,line width=0.5mm,blue,dashdotted] (2.00,4.55)--(2.50,4.80);\draw[line width=0.5mm,line width=0.5mm,blue,dashdotted] (2.50,4.80)--(3.00,4.85);\draw[line width=0.5mm,line width=0.5mm,blue,dashdotted] (3.00,4.85)--(3.50,4.90);\draw[line width=0.5mm,line width=0.5mm,blue,dashdotted] (3.50,4.90)--(4.00,4.95);\draw[line width=0.5mm,line width=0.5mm,blue,dashdotted] (4.00,4.95)--(4.50,5.00);\draw[line width=0.5mm,line width=0.5mm,blue,dashdotted] (4.50,5.00)--(5.00,5.00);\draw[line width=0.5mm,line width=0.5mm,blue,dashdotted] (5.00,5.00)--(5.50,5.00);\draw[line width=0.5mm,line width=0.5mm,blue,dashdotted] (5.50,5.00)--(6.00,5.00);\draw[line width=0.5mm,line width=0.5mm,blue,dashdotted] (6.00,5.00)--(6.50,5.00);\draw[line width=0.5mm,line width=0.5mm,blue,dashdotted] (6.50,5.00)--(7.00,5.00);\draw[line width=0.5mm,line width=0.5mm,blue,dashdotted] (7.00,5.00)--(7.50,5.00);\draw[line width=0.5mm,line width=0.5mm,blue,dashdotted] (7.50,5.00)--(8.00,5.00);\draw[line width=0.5mm,line width=0.5mm,blue,dashdotted] (8.00,5.00)--(8.50,5.00);\draw[line width=0.5mm,line width=0.5mm,blue,dashdotted] (8.50,5.00)--(9.00,5.00);\draw[line width=0.5mm,line width=0.5mm,blue,dashdotted] (9.00,5.00)--(9.50,5.00);\draw[line width=0.5mm,line width=0.5mm,blue,dashdotted] (9.50,5.00)--(10.00,5.00);
\draw[line width=0.5mm, teal, dotted] (8.5, 3.00)--(9.0, 3.00);
\draw[line width=0.5mm, teal] (9.0, 3.00) node[anchor=west] {DIF};
\draw[line width=0.5mm, blue, dashdotted] (8.5, 3.50)--(9.0, 3.50);
\draw[line width=0.5mm, blue] (9.0, 3.50) node[anchor=west] {SBD};
\end{tikzpicture}
            \label{fig:optgaps-c2}
         \end{subfigure}        
    \end{minipage}
    \begin{minipage}{.99\textwidth}
        \centering
        \caption{Comparison between ABD and SBD, where the $y$ axis presents the number of instances with a value (solution time or optimality gap) smaller than or  equal to the reference value on the $x$ axis.}
        \label{fig:c3}
        \vspace{0.3cm}
         \begin{subfigure}[b]{0.49\textwidth}
            \centering
            \caption{Instances solved to optimality by both methods.}
            \begin{tikzpicture}[scale=.7, every node/.style={scale=.7}]
\draw[line width=0.5mm,thick,->] (0,0) -- (10.5,0);
\draw[line width=0.5mm,thick,->] (0,0) -- (0,5.5);
\draw (9.5,0.5) node[anchor=mid] {runtime (min)};
\draw (0,6) node[anchor=mid] {\# instances};
\draw (0,-0.5) node[anchor=mid] {$0$};
\draw (1,-0.5) node[anchor=mid] {$6$};
\draw (2,-0.5) node[anchor=mid] {$12$};
\draw (3,-0.5) node[anchor=mid] {$18$};
\draw (4,-0.5) node[anchor=mid] {$24$};
\draw (5,-0.5) node[anchor=mid] {$30$};
\draw (6,-0.5) node[anchor=mid] {$36$};
\draw (7,-0.5) node[anchor=mid] {$42$};
\draw (8,-0.5) node[anchor=mid] {$48$};
\draw (9,-0.5) node[anchor=mid] {$54$};
\draw (10,-0.5) node[anchor=mid] {$60$};
\draw (-0.5,0) node[anchor=mid] {$0$};
\draw (-0.5,1) node[anchor=mid] {$40$};
\draw (-0.5,2) node[anchor=mid] {$79$};
\draw (-0.5,3) node[anchor=mid] {$119$};
\draw (-0.5,4) node[anchor=mid] {$158$};
\draw (-0.5,5) node[anchor=mid] {$198$};
\draw[line width=0.5mm,line width=0.5mm,blue,dashdotted] (0.00,0.00)--(0.50,4.15);\draw[line width=0.5mm,line width=0.5mm,blue,dashdotted] (0.50,4.15)--(1.00,4.35);\draw[line width=0.5mm,line width=0.5mm,blue,dashdotted] (1.00,4.35)--(1.50,4.50);\draw[line width=0.5mm,line width=0.5mm,blue,dashdotted] (1.50,4.50)--(2.00,4.60);\draw[line width=0.5mm,line width=0.5mm,blue,dashdotted] (2.00,4.60)--(2.50,4.70);\draw[line width=0.5mm,line width=0.5mm,blue,dashdotted] (2.50,4.70)--(3.00,4.75);\draw[line width=0.5mm,line width=0.5mm,blue,dashdotted] (3.00,4.75)--(3.50,4.85);\draw[line width=0.5mm,line width=0.5mm,blue,dashdotted] (3.50,4.85)--(4.00,4.90);\draw[line width=0.5mm,line width=0.5mm,blue,dashdotted] (4.00,4.90)--(4.50,4.90);\draw[line width=0.5mm,line width=0.5mm,blue,dashdotted] (4.50,4.90)--(5.00,4.95);\draw[line width=0.5mm,line width=0.5mm,blue,dashdotted] (5.00,4.95)--(5.50,5.00);\draw[line width=0.5mm,line width=0.5mm,blue,dashdotted] (5.50,5.00)--(6.00,5.00);\draw[line width=0.5mm,line width=0.5mm,blue,dashdotted] (6.00,5.00)--(6.50,5.00);\draw[line width=0.5mm,line width=0.5mm,blue,dashdotted] (6.50,5.00)--(7.00,5.00);\draw[line width=0.5mm,line width=0.5mm,blue,dashdotted] (7.00,5.00)--(7.50,5.00);\draw[line width=0.5mm,line width=0.5mm,blue,dashdotted] (7.50,5.00)--(8.00,5.00);\draw[line width=0.5mm,line width=0.5mm,blue,dashdotted] (8.00,5.00)--(8.50,5.00);\draw[line width=0.5mm,line width=0.5mm,blue,dashdotted] (8.50,5.00)--(9.00,5.00);\draw[line width=0.5mm,line width=0.5mm,blue,dashdotted] (9.00,5.00)--(9.50,5.00);\draw[line width=0.5mm,line width=0.5mm,blue,dashdotted] (9.50,5.00)--(10.00,5.00);\draw[line width=0.5mm,line width=0.5mm,orange,solid] (0.00,0.00)--(0.50,4.25);\draw[line width=0.5mm,line width=0.5mm,orange,solid] (0.50,4.25)--(1.00,4.50);\draw[line width=0.5mm,line width=0.5mm,orange,solid] (1.00,4.50)--(1.50,4.65);\draw[line width=0.5mm,line width=0.5mm,orange,solid] (1.50,4.65)--(2.00,4.75);\draw[line width=0.5mm,line width=0.5mm,orange,solid] (2.00,4.75)--(2.50,4.85);\draw[line width=0.5mm,line width=0.5mm,orange,solid] (2.50,4.85)--(3.00,4.85);\draw[line width=0.5mm,line width=0.5mm,orange,solid] (3.00,4.85)--(3.50,4.85);\draw[line width=0.5mm,line width=0.5mm,orange,solid] (3.50,4.85)--(4.00,4.85);\draw[line width=0.5mm,line width=0.5mm,orange,solid] (4.00,4.85)--(4.50,4.90);\draw[line width=0.5mm,line width=0.5mm,orange,solid] (4.50,4.90)--(5.00,4.90);\draw[line width=0.5mm,line width=0.5mm,orange,solid] (5.00,4.90)--(5.50,5.00);\draw[line width=0.5mm,line width=0.5mm,orange,solid] (5.50,5.00)--(6.00,5.00);\draw[line width=0.5mm,line width=0.5mm,orange,solid] (6.00,5.00)--(6.50,5.00);\draw[line width=0.5mm,line width=0.5mm,orange,solid] (6.50,5.00)--(7.00,5.00);\draw[line width=0.5mm,line width=0.5mm,orange,solid] (7.00,5.00)--(7.50,5.00);\draw[line width=0.5mm,line width=0.5mm,orange,solid] (7.50,5.00)--(8.00,5.00);\draw[line width=0.5mm,line width=0.5mm,orange,solid] (8.00,5.00)--(8.50,5.00);\draw[line width=0.5mm,line width=0.5mm,orange,solid] (8.50,5.00)--(9.00,5.00);\draw[line width=0.5mm,line width=0.5mm,orange,solid] (9.00,5.00)--(9.50,5.00);\draw[line width=0.5mm,line width=0.5mm,orange,solid] (9.50,5.00)--(10.00,5.00);
\draw[line width=0.5mm, blue, dashdotted] (8.5, 3.00)--(9.0, 3.00);
\draw[line width=0.5mm, blue] (9.0, 3.00) node[anchor=west] {SBD};
\draw[line width=0.5mm, orange, solid] (8.5, 3.50)--(9.0, 3.50);
\draw[line width=0.5mm, orange] (9.0, 3.50) node[anchor=west] {ABD};
\end{tikzpicture}
            \label{fig:runtimes-c3}
         \end{subfigure}
         \hfill
         \begin{subfigure}[b]{0.49\textwidth}
            \centering
            \caption{Instances not solved to optimality by both methods.}
            \begin{tikzpicture}[scale=.7, every node/.style={scale=.7}]
\draw[line width=0.5mm,thick,->] (0,0) -- (10.5,0);
\draw[line width=0.5mm,thick,->] (0,0) -- (0,5.5);
\draw (9.5,0.5) node[anchor=mid] {optimality gap (\%)};
\draw (0,6) node[anchor=mid] {\# instances};
\draw (0,-0.5) node[anchor=mid] {$0$};
\draw (1,-0.5) node[anchor=mid] {$10$};
\draw (2,-0.5) node[anchor=mid] {$20$};
\draw (3,-0.5) node[anchor=mid] {$30$};
\draw (4,-0.5) node[anchor=mid] {$40$};
\draw (5,-0.5) node[anchor=mid] {$50$};
\draw (6,-0.5) node[anchor=mid] {$60$};
\draw (7,-0.5) node[anchor=mid] {$70$};
\draw (8,-0.5) node[anchor=mid] {$80$};
\draw (9,-0.5) node[anchor=mid] {$90$};
\draw (10,-0.5) node[anchor=mid] {$100$};
\draw (-0.5,0) node[anchor=mid] {$0$};
\draw (-0.5,1) node[anchor=mid] {$4$};
\draw (-0.5,2) node[anchor=mid] {$7$};
\draw (-0.5,3) node[anchor=mid] {$11$};
\draw (-0.5,4) node[anchor=mid] {$14$};
\draw (-0.5,5) node[anchor=mid] {$18$};
\draw[line width=0.5mm,line width=0.5mm,blue,dashdotted] (0.00,0.00)--(0.50,1.10);\draw[line width=0.5mm,line width=0.5mm,blue,dashdotted] (0.50,1.10)--(1.00,4.15);\draw[line width=0.5mm,line width=0.5mm,blue,dashdotted] (1.00,4.15)--(1.50,5.00);\draw[line width=0.5mm,line width=0.5mm,blue,dashdotted] (1.50,5.00)--(2.00,5.00);\draw[line width=0.5mm,line width=0.5mm,blue,dashdotted] (2.00,5.00)--(2.50,5.00);\draw[line width=0.5mm,line width=0.5mm,blue,dashdotted] (2.50,5.00)--(3.00,5.00);\draw[line width=0.5mm,line width=0.5mm,blue,dashdotted] (3.00,5.00)--(3.50,5.00);\draw[line width=0.5mm,line width=0.5mm,blue,dashdotted] (3.50,5.00)--(4.00,5.00);\draw[line width=0.5mm,line width=0.5mm,blue,dashdotted] (4.00,5.00)--(4.50,5.00);\draw[line width=0.5mm,line width=0.5mm,blue,dashdotted] (4.50,5.00)--(5.00,5.00);\draw[line width=0.5mm,line width=0.5mm,blue,dashdotted] (5.00,5.00)--(5.50,5.00);\draw[line width=0.5mm,line width=0.5mm,blue,dashdotted] (5.50,5.00)--(6.00,5.00);\draw[line width=0.5mm,line width=0.5mm,blue,dashdotted] (6.00,5.00)--(6.50,5.00);\draw[line width=0.5mm,line width=0.5mm,blue,dashdotted] (6.50,5.00)--(7.00,5.00);\draw[line width=0.5mm,line width=0.5mm,blue,dashdotted] (7.00,5.00)--(7.50,5.00);\draw[line width=0.5mm,line width=0.5mm,blue,dashdotted] (7.50,5.00)--(8.00,5.00);\draw[line width=0.5mm,line width=0.5mm,blue,dashdotted] (8.00,5.00)--(8.50,5.00);\draw[line width=0.5mm,line width=0.5mm,blue,dashdotted] (8.50,5.00)--(9.00,5.00);\draw[line width=0.5mm,line width=0.5mm,blue,dashdotted] (9.00,5.00)--(9.50,5.00);\draw[line width=0.5mm,line width=0.5mm,blue,dashdotted] (9.50,5.00)--(10.00,5.00);\draw[line width=0.5mm,line width=0.5mm,orange,solid] (0.00,0.00)--(0.50,1.35);\draw[line width=0.5mm,line width=0.5mm,orange,solid] (0.50,1.35)--(1.00,4.40);\draw[line width=0.5mm,line width=0.5mm,orange,solid] (1.00,4.40)--(1.50,5.00);\draw[line width=0.5mm,line width=0.5mm,orange,solid] (1.50,5.00)--(2.00,5.00);\draw[line width=0.5mm,line width=0.5mm,orange,solid] (2.00,5.00)--(2.50,5.00);\draw[line width=0.5mm,line width=0.5mm,orange,solid] (2.50,5.00)--(3.00,5.00);\draw[line width=0.5mm,line width=0.5mm,orange,solid] (3.00,5.00)--(3.50,5.00);\draw[line width=0.5mm,line width=0.5mm,orange,solid] (3.50,5.00)--(4.00,5.00);\draw[line width=0.5mm,line width=0.5mm,orange,solid] (4.00,5.00)--(4.50,5.00);\draw[line width=0.5mm,line width=0.5mm,orange,solid] (4.50,5.00)--(5.00,5.00);\draw[line width=0.5mm,line width=0.5mm,orange,solid] (5.00,5.00)--(5.50,5.00);\draw[line width=0.5mm,line width=0.5mm,orange,solid] (5.50,5.00)--(6.00,5.00);\draw[line width=0.5mm,line width=0.5mm,orange,solid] (6.00,5.00)--(6.50,5.00);\draw[line width=0.5mm,line width=0.5mm,orange,solid] (6.50,5.00)--(7.00,5.00);\draw[line width=0.5mm,line width=0.5mm,orange,solid] (7.00,5.00)--(7.50,5.00);\draw[line width=0.5mm,line width=0.5mm,orange,solid] (7.50,5.00)--(8.00,5.00);\draw[line width=0.5mm,line width=0.5mm,orange,solid] (8.00,5.00)--(8.50,5.00);\draw[line width=0.5mm,line width=0.5mm,orange,solid] (8.50,5.00)--(9.00,5.00);\draw[line width=0.5mm,line width=0.5mm,orange,solid] (9.00,5.00)--(9.50,5.00);\draw[line width=0.5mm,line width=0.5mm,orange,solid] (9.50,5.00)--(10.00,5.00);
\draw[line width=0.5mm, blue, dashdotted] (8.5, 3.00)--(9.0, 3.00);
\draw[line width=0.5mm, blue] (9.0, 3.00) node[anchor=west] {SBD};
\draw[line width=0.5mm, orange, solid] (8.5, 3.50)--(9.0, 3.50);
\draw[line width=0.5mm, orange] (9.0, 3.50) node[anchor=west] {ABD};
\end{tikzpicture}
            \label{fig:optgaps-c3}
         \end{subfigure}        
    \end{minipage}
\end{figure}

\subsubsection{Comparison Between DIF and SIF.}
Theorem~\ref{thm:tightness} guarantees that the DIF provides a tighter continuous relaxation than the SIF, which should help off-the-shelf solvers to find solutions and prove their optimality faster.
However, this is not always the case nowadays due to other built-in techniques that might impact the branch-and-bound tree in an unexpected (or randomized) manner.
We therefore evaluate how well Gurobi solves both formulations.
Figure~\ref{fig:runtimes-c1} presents solution times for instances solved to optimality by the DIF and the SIF, whereas Figure~\ref{fig:optgaps-c1} presents (proven) optimality gaps for instances not solved to optimality by at least one of them.
We also report solution times and (proven) optimality gaps averaged over different instance attributes in Online Appendix~F.
%
%
The DIF clearly has faster solution times than the SIF for instances solved to optimality by both formulations, being approximately $53\%$ faster on average ($4.48$ minutes versus $6.48$ minutes).
Moreover, the DIF proves smaller optimality gaps than the SIF within the time limit for instances not solved to optimality by at least one of the formulations, which are $57\%$ smaller on average ($1.97 \%$ versus $6.13 \%$).
A closer look into the integrality gaps of instances solved to optimality show that the DIF has an average integrality gap three times smaller than one of the SIF.
These results confirm that tighter continuous relaxation bounds provided by the DIF, evidenced by the average integrality gap, are advantageous to Gurobi when solving the \ourproblem{} in practice.

\subsubsection{Comparison Between SBD and DIF.}
The Benders Decomposition proposed in Section~\ref{sub:benders} allows us to avoid handling many variables $x^{\ell{}t}_{ij}$ in the DIF, but the branch-and-Benders-cut implementation may require numerous optimality cuts to converge to the optimal solution in practice.
In other words, the proposed Benders Decomposition may actually have a worse performance than the DIF if the disadvantages of having many variables $x^{\ell{}t}_{ij}$ do not outweigh the disadvantages of needing numerous optimality cuts.
We therefore assess how much better the proposed Benders decomposition performs relative to the DIF.
Figures~\ref{fig:runtimes-c2} and~\ref{fig:optgaps-c2} have analogous meaning to Figures~\ref{fig:runtimes-c1} and~\ref{fig:optgaps-c1}, respectively.
We also report solution times and (proven) optimality gaps averaged over different instance attributes in Online Appendix~F.
%
%
On average, the SBD is approximately five times faster than the DIF ($1.11$ minutes versus $5.72$), and proves optimality gaps approximately two times smaller ($7.42 \%$ versus $14.73\%$).
These results suggest that the proposed Benders Decomposition should be preferred over the DIF for the \ourproblem{}.

\subsubsection{Comparison Between ABD and SBD.}
Algorithm~\ref{alg:analytical} describes an analytical procedure to compute optimality cuts.
Such a procedure avoids the need to solve linear programs with Gurobi and allows us to generate optimality cuts with a similar structure among themselves.
There is no theoretical guarantee, however, that optimality cuts obtained through the analytical procedure are better or worse than those obtained by solving the linear program with Gurobi.
We now compare the performance of these two implementations of the branch-and-Benders-cut in practice.
Since the analytical procedure only works for instances of the 1-\ourproblem{}, we filter our benchmark to consider only $216$ instances with a single facility.
Figures~\ref{fig:runtimes-c3} and~\ref{fig:optgaps-c3} have analogous meaning to Figures~\ref{fig:runtimes-c1} and~\ref{fig:optgaps-c1}, respectively.
We also report solution times and optimality gaps averaged over different instance attributes in Online Appendix~F.
%
%
On average, the ABD has $35\%$ faster solution times than SBD ($1.85$ minutes versus $2.50$ minutes) and  {manages to prove optimality gaps about $5\%$ smaller ($6.07 \%$ versus $6.43 \%$).
These results indicate that the analytical procedure decreases the computational burden of computing optimality cuts, but unfortunately not by a large margin.

\subsection{Managerial Insights}
\label{sub:managerial-insights}

We now draw managerial insights from the benchmark proposed for the computational experiments.

\subsubsection{Impact of Heuristic Decisions.}
As mentioned in Section~\ref{sub:heuristics}, the provider might implicitly or explicitly ignore cumulative customer demand, thus deciding myopically where to place facilities over the planning horizon.
In this sense, we now investigate the impact of myopic and, more generally, heuristic decisions.
More specifically, we evaluate whether the heuristics  (\ie, DBH, FGH and BGH) provide sufficiently high-quality solutions for $370$ instances of the benchmark with a known optimal solution.
To this end, we define the opportunity gap for each heuristic as $\frac{Z^{\star} - Z^{\prime}}{Z^{\star}}$, where $Z^{\star}$ is the optimal objective value obtained through an exact method and $Z^{\prime}$ is the objective value of the heuristic at hand.
Intuitively, large (respectively, small) opportunity gaps indicate that the heuristic finds low-quality (respectively, high-quality) location policies.
In particular, small opportunity gaps indicate that the heuristic could be employed whenever the provider cannot apply exact methods (\eg, for large-scale instances or when there are no off-the-shelf solvers available).
As a sanity check, we also evaluate a randomly generated solution denoted as RND.
Figure~\ref{fig:heuristics} presents the performance profile of the heuristics, where the $y$ axis presents the number of instances with an opportunity gap smaller than or equal to the reference value on the $x$ axis for each heuristic.
We report opportunity gaps averaged over different instance attributes in the Online Appendix~F.

\begin{figure}
    \centering
    \caption{Performance profile of the heuristics, where the $y$ axis presents the number of instances with an opportunity gap smaller than or  equal to the reference value on the $x$ axis for each heuristic.}
    \begin{tikzpicture}[scale=.8, every node/.style={scale=.8}]
\draw[line width=0.5mm,thick,->] (0,0) -- (10.5,0);
\draw[line width=0.5mm,thick,->] (0,0) -- (0,10.5);
\draw (9.5,0.5) node[anchor=mid] {opportunity gap (\%)};
\draw (0,11) node[anchor=mid] {\# instances};
\draw (0,-0.5) node[anchor=mid] {$0$};
\draw (1,-0.5) node[anchor=mid] {$10$};
\draw (2,-0.5) node[anchor=mid] {$20$};
\draw (3,-0.5) node[anchor=mid] {$30$};
\draw (4,-0.5) node[anchor=mid] {$40$};
\draw (5,-0.5) node[anchor=mid] {$50$};
\draw (6,-0.5) node[anchor=mid] {$60$};
\draw (7,-0.5) node[anchor=mid] {$70$};
\draw (8,-0.5) node[anchor=mid] {$80$};
\draw (9,-0.5) node[anchor=mid] {$90$};
\draw (10,-0.5) node[anchor=mid] {$100$};
\draw (-0.5,0) node[anchor=mid] {$0$};
\draw (-0.5,1) node[anchor=mid] {$37$};
\draw (-0.5,2) node[anchor=mid] {$74$};
\draw (-0.5,3) node[anchor=mid] {$111$};
\draw (-0.5,4) node[anchor=mid] {$148$};
\draw (-0.5,5) node[anchor=mid] {$185$};
\draw (-0.5,6) node[anchor=mid] {$222$};
\draw (-0.5,7) node[anchor=mid] {$259$};
\draw (-0.5,8) node[anchor=mid] {$296$};
\draw (-0.5,9) node[anchor=mid] {$333$};
\draw (-0.5,10) node[anchor=mid] {$370$};
\draw[line width=0.5mm,line width=0.5mm,red,dashed] (0.00,0.00)--(0.50,0.20);\draw[line width=0.5mm,line width=0.5mm,red,dashed] (0.50,0.20)--(1.00,0.90);\draw[line width=0.5mm,line width=0.5mm,red,dashed] (1.00,0.90)--(1.50,2.00);\draw[line width=0.5mm,line width=0.5mm,red,dashed] (1.50,2.00)--(2.00,3.10);\draw[line width=0.5mm,line width=0.5mm,red,dashed] (2.00,3.10)--(2.50,3.80);\draw[line width=0.5mm,line width=0.5mm,red,dashed] (2.50,3.80)--(3.00,4.50);\draw[line width=0.5mm,line width=0.5mm,red,dashed] (3.00,4.50)--(3.50,5.20);\draw[line width=0.5mm,line width=0.5mm,red,dashed] (3.50,5.20)--(4.00,5.70);\draw[line width=0.5mm,line width=0.5mm,red,dashed] (4.00,5.70)--(4.50,5.90);\draw[line width=0.5mm,line width=0.5mm,red,dashed] (4.50,5.90)--(5.00,6.30);\draw[line width=0.5mm,line width=0.5mm,red,dashed] (5.00,6.30)--(5.50,6.70);\draw[line width=0.5mm,line width=0.5mm,red,dashed] (5.50,6.70)--(6.00,7.00);\draw[line width=0.5mm,line width=0.5mm,red,dashed] (6.00,7.00)--(6.50,7.80);\draw[line width=0.5mm,line width=0.5mm,red,dashed] (6.50,7.80)--(7.00,8.50);\draw[line width=0.5mm,line width=0.5mm,red,dashed] (7.00,8.50)--(7.50,9.40);\draw[line width=0.5mm,line width=0.5mm,red,dashed] (7.50,9.40)--(8.00,10.00);\draw[line width=0.5mm,line width=0.5mm,red,dashed] (8.00,10.00)--(8.50,10.00);\draw[line width=0.5mm,line width=0.5mm,red,dashed] (8.50,10.00)--(9.00,10.00);\draw[line width=0.5mm,line width=0.5mm,red,dashed] (9.00,10.00)--(9.50,10.00);\draw[line width=0.5mm,line width=0.5mm,red,dashed] (9.50,10.00)--(10.00,10.00);\draw[line width=0.5mm,line width=0.5mm,gray,dotted] (0.00,0.00)--(0.50,0.00);\draw[line width=0.5mm,line width=0.5mm,gray,dotted] (0.50,0.00)--(1.00,0.40);\draw[line width=0.5mm,line width=0.5mm,gray,dotted] (1.00,0.40)--(1.50,1.20);\draw[line width=0.5mm,line width=0.5mm,gray,dotted] (1.50,1.20)--(2.00,2.40);\draw[line width=0.5mm,line width=0.5mm,gray,dotted] (2.00,2.40)--(2.50,4.00);\draw[line width=0.5mm,line width=0.5mm,gray,dotted] (2.50,4.00)--(3.00,5.80);\draw[line width=0.5mm,line width=0.5mm,gray,dotted] (3.00,5.80)--(3.50,7.10);\draw[line width=0.5mm,line width=0.5mm,gray,dotted] (3.50,7.10)--(4.00,8.20);\draw[line width=0.5mm,line width=0.5mm,gray,dotted] (4.00,8.20)--(4.50,9.10);\draw[line width=0.5mm,line width=0.5mm,gray,dotted] (4.50,9.10)--(5.00,9.50);\draw[line width=0.5mm,line width=0.5mm,gray,dotted] (5.00,9.50)--(5.50,9.80);\draw[line width=0.5mm,line width=0.5mm,gray,dotted] (5.50,9.80)--(6.00,9.90);\draw[line width=0.5mm,line width=0.5mm,gray,dotted] (6.00,9.90)--(6.50,10.00);\draw[line width=0.5mm,line width=0.5mm,gray,dotted] (6.50,10.00)--(7.00,10.00);\draw[line width=0.5mm,line width=0.5mm,gray,dotted] (7.00,10.00)--(7.50,10.00);\draw[line width=0.5mm,line width=0.5mm,gray,dotted] (7.50,10.00)--(8.00,10.00);\draw[line width=0.5mm,line width=0.5mm,gray,dotted] (8.00,10.00)--(8.50,10.00);\draw[line width=0.5mm,line width=0.5mm,gray,dotted] (8.50,10.00)--(9.00,10.00);\draw[line width=0.5mm,line width=0.5mm,gray,dotted] (9.00,10.00)--(9.50,10.00);\draw[line width=0.5mm,line width=0.5mm,gray,dotted] (9.50,10.00)--(10.00,10.00);\draw[line width=0.5mm,line width=0.5mm,olive,dashdotted] (0.00,0.00)--(0.50,5.10);\draw[line width=0.5mm,line width=0.5mm,olive,dashdotted] (0.50,5.10)--(1.00,8.10);\draw[line width=0.5mm,line width=0.5mm,olive,dashdotted] (1.00,8.10)--(1.50,9.50);\draw[line width=0.5mm,line width=0.5mm,olive,dashdotted] (1.50,9.50)--(2.00,9.90);\draw[line width=0.5mm,line width=0.5mm,olive,dashdotted] (2.00,9.90)--(2.50,10.00);\draw[line width=0.5mm,line width=0.5mm,olive,dashdotted] (2.50,10.00)--(3.00,10.00);\draw[line width=0.5mm,line width=0.5mm,olive,dashdotted] (3.00,10.00)--(3.50,10.00);\draw[line width=0.5mm,line width=0.5mm,olive,dashdotted] (3.50,10.00)--(4.00,10.00);\draw[line width=0.5mm,line width=0.5mm,olive,dashdotted] (4.00,10.00)--(4.50,10.00);\draw[line width=0.5mm,line width=0.5mm,olive,dashdotted] (4.50,10.00)--(5.00,10.00);\draw[line width=0.5mm,line width=0.5mm,olive,dashdotted] (5.00,10.00)--(5.50,10.00);\draw[line width=0.5mm,line width=0.5mm,olive,dashdotted] (5.50,10.00)--(6.00,10.00);\draw[line width=0.5mm,line width=0.5mm,olive,dashdotted] (6.00,10.00)--(6.50,10.00);\draw[line width=0.5mm,line width=0.5mm,olive,dashdotted] (6.50,10.00)--(7.00,10.00);\draw[line width=0.5mm,line width=0.5mm,olive,dashdotted] (7.00,10.00)--(7.50,10.00);\draw[line width=0.5mm,line width=0.5mm,olive,dashdotted] (7.50,10.00)--(8.00,10.00);\draw[line width=0.5mm,line width=0.5mm,olive,dashdotted] (8.00,10.00)--(8.50,10.00);\draw[line width=0.5mm,line width=0.5mm,olive,dashdotted] (8.50,10.00)--(9.00,10.00);\draw[line width=0.5mm,line width=0.5mm,olive,dashdotted] (9.00,10.00)--(9.50,10.00);\draw[line width=0.5mm,line width=0.5mm,olive,dashdotted] (9.50,10.00)--(10.00,10.00);\draw[line width=0.5mm,line width=0.5mm,orange,solid] (0.00,2.00)--(0.50,9.10);\draw[line width=0.5mm,line width=0.5mm,orange,solid] (0.50,9.10)--(1.00,10.00);\draw[line width=0.5mm,line width=0.5mm,orange,solid] (1.00,10.00)--(1.50,10.00);\draw[line width=0.5mm,line width=0.5mm,orange,solid] (1.50,10.00)--(2.00,10.00);\draw[line width=0.5mm,line width=0.5mm,orange,solid] (2.00,10.00)--(2.50,10.00);\draw[line width=0.5mm,line width=0.5mm,orange,solid] (2.50,10.00)--(3.00,10.00);\draw[line width=0.5mm,line width=0.5mm,orange,solid] (3.00,10.00)--(3.50,10.00);\draw[line width=0.5mm,line width=0.5mm,orange,solid] (3.50,10.00)--(4.00,10.00);\draw[line width=0.5mm,line width=0.5mm,orange,solid] (4.00,10.00)--(4.50,10.00);\draw[line width=0.5mm,line width=0.5mm,orange,solid] (4.50,10.00)--(5.00,10.00);\draw[line width=0.5mm,line width=0.5mm,orange,solid] (5.00,10.00)--(5.50,10.00);\draw[line width=0.5mm,line width=0.5mm,orange,solid] (5.50,10.00)--(6.00,10.00);\draw[line width=0.5mm,line width=0.5mm,orange,solid] (6.00,10.00)--(6.50,10.00);\draw[line width=0.5mm,line width=0.5mm,orange,solid] (6.50,10.00)--(7.00,10.00);\draw[line width=0.5mm,line width=0.5mm,orange,solid] (7.00,10.00)--(7.50,10.00);\draw[line width=0.5mm,line width=0.5mm,orange,solid] (7.50,10.00)--(8.00,10.00);\draw[line width=0.5mm,line width=0.5mm,orange,solid] (8.00,10.00)--(8.50,10.00);\draw[line width=0.5mm,line width=0.5mm,orange,solid] (8.50,10.00)--(9.00,10.00);\draw[line width=0.5mm,line width=0.5mm,orange,solid] (9.00,10.00)--(9.50,10.00);\draw[line width=0.5mm,line width=0.5mm,orange,solid] (9.50,10.00)--(10.00,10.00);
\draw[line width=0.5mm, red, dashed] (8.5, 3.00)--(9.0, 3.00);
\draw[line width=0.5mm, red] (9.0, 3.00) node[anchor=west] {DBH};
\draw[line width=0.5mm, gray, dotted] (8.5, 3.50)--(9.0, 3.50);
\draw[line width=0.5mm, gray] (9.0, 3.50) node[anchor=west] {RND};
\draw[line width=0.5mm, olive, dashdotted] (8.5, 4.00)--(9.0, 4.00);
\draw[line width=0.5mm, olive] (9.0, 4.00) node[anchor=west] {FGH};
\draw[line width=0.5mm, orange, solid] (8.5, 4.50)--(9.0, 4.50);
\draw[line width=0.5mm, orange] (9.0, 4.50) node[anchor=west] {BGH};
\end{tikzpicture}
    \label{fig:heuristics}
\end{figure}

We first look at the myopic heuristics (\ie, DBH and FGH).
On the one hand, the DBH has the worst performance in terms of opportunity gap, even worse than RND (\ie, our sanity check).
On the other hand, the FGH performs considerably better than the sanity check, with an average opportunity gap of $6.07\%$.
Nevertheless, ignoring cumulative customer demand within the optimization framework leads to a considerable lost opportunity.
We now turn to the BGH, inspired by our theoretical results.
We highlight that the BGH is not completely myopic because it does not neglect future effects of current location decisions when computing a location policy.
The BGH performs surprisingly much better than expected, even for instances with different rewards, having an average opportunity gap of $1.74\%$.
From a theoretical point of view, this result indicates that tighter approximation guarantees may be achievable for the 1-\ourproblem{} with identical rewards -- this is reasonable, since the proof of Theorem~\ref{thm:approximation} does not provide a tight upper bound to the optimal objective value. 
From a practical point of view, this result underlines the relevance of the BGH to obtain high-quality solutions for large-scale instances, where exact methods may not be applicable.
Although the BGH performs considerably well for the benchmark, the development of exact methods remains important because the \ourproblem{} cannot be approximated within a constant factor, as stated by Theorem~\ref{thm:nphard-single}.
In other words, there can be instances where the BGH performs arbitrarily bad. Although these instances are not part of our benchmark, they may appear in practice when solving the \ourproblem{}.

\begin{figure}
    \begin{minipage}{.99\textwidth}
         \centering
         \caption{Number of customers per number of captures for the instance with identical rewards, long rankings, constant spawning demands, and no penalties.} 
        \vspace{0.3cm}
         \begin{subfigure}[b]{0.32\textwidth}
            \centering
            \caption{Optimal solution for 1 facility.}
            \begin{tikzpicture}[scale=.4, every node/.style={scale=.4}]
\draw[line width=0.5mm,thick,->] (0,0) -- (7.5,0);
\draw[line width=0.5mm,thick,->] (0,0) -- (0,5.5);
\draw (7.5,-0.5) node[anchor=mid] {\# captures};
\draw (0,6) node[anchor=mid] {\# customers};
\draw (1,-0.5) node[anchor=mid] {$0$};
\draw (2,-0.5) node[anchor=mid] {$1$};
\draw (3,-0.5) node[anchor=mid] {$2$};
\draw (4,-0.5) node[anchor=mid] {$3$};
\draw (5,-0.5) node[anchor=mid] {$4$};
\draw (6,-0.5) node[anchor=mid] {$5$};
\draw (-0.5,0) node[anchor=mid] {$0$};
\draw (-0.5,1) node[anchor=mid] {$10$};
\draw (-0.5,2) node[anchor=mid] {$20$};
\draw (-0.5,3) node[anchor=mid] {$30$};
\draw (-0.5,4) node[anchor=mid] {$40$};
\draw (-0.5,5) node[anchor=mid] {$50$};
\draw[line width=0.5mm,gray,fill=gray] (0.60,0) rectangle ++(0.8,0.9);\draw[line width=0.5mm,gray,fill=gray] (1.60,0) rectangle ++(0.8,3.4);\draw[line width=0.5mm,gray,fill=gray] (2.60,0) rectangle ++(0.8,0.6);\draw[line width=0.5mm,gray,fill=gray] (3.60,0) rectangle ++(0.8,0.1);\draw[line width=0.5mm,gray,fill=gray] (4.60,0) rectangle ++(0.8,0.0);\draw[line width=0.5mm,gray,fill=gray] (5.60,0) rectangle ++(0.8,0.0);
\end{tikzpicture}
            
         \end{subfigure}
        \hfill
         \begin{subfigure}[b]{0.32\textwidth}
            \centering
            \caption{Optimal solution for 3 facilities.}
            \begin{tikzpicture}[scale=.4, every node/.style={scale=.4}]
\draw[line width=0.5mm,thick,->] (0,0) -- (7.5,0);
\draw[line width=0.5mm,thick,->] (0,0) -- (0,5.5);
\draw (7.5,-0.5) node[anchor=mid] {\# captures};
\draw (0,6) node[anchor=mid] {\# customers};
\draw (1,-0.5) node[anchor=mid] {$0$};
\draw (2,-0.5) node[anchor=mid] {$1$};
\draw (3,-0.5) node[anchor=mid] {$2$};
\draw (4,-0.5) node[anchor=mid] {$3$};
\draw (5,-0.5) node[anchor=mid] {$4$};
\draw (6,-0.5) node[anchor=mid] {$5$};
\draw (-0.5,0) node[anchor=mid] {$0$};
\draw (-0.5,1) node[anchor=mid] {$10$};
\draw (-0.5,2) node[anchor=mid] {$20$};
\draw (-0.5,3) node[anchor=mid] {$30$};
\draw (-0.5,4) node[anchor=mid] {$40$};
\draw (-0.5,5) node[anchor=mid] {$50$};
\draw[line width=0.5mm,gray,fill=gray] (0.60,0) rectangle ++(0.8,0.0);\draw[line width=0.5mm,gray,fill=gray] (1.60,0) rectangle ++(0.8,2.1);\draw[line width=0.5mm,gray,fill=gray] (2.60,0) rectangle ++(0.8,2.4);\draw[line width=0.5mm,gray,fill=gray] (3.60,0) rectangle ++(0.8,0.5);\draw[line width=0.5mm,gray,fill=gray] (4.60,0) rectangle ++(0.8,0.0);\draw[line width=0.5mm,gray,fill=gray] (5.60,0) rectangle ++(0.8,0.0);
\end{tikzpicture}
            
         \end{subfigure}
        \hfill
        \begin{subfigure}[b]{0.32\textwidth}
            \centering
            \caption{Optimal solution for 5 facilities.}
            \begin{tikzpicture}[scale=.4, every node/.style={scale=.4}]
\draw[line width=0.5mm,thick,->] (0,0) -- (7.5,0);
\draw[line width=0.5mm,thick,->] (0,0) -- (0,5.5);
\draw (7.5,-0.5) node[anchor=mid] {\# captures};
\draw (0,6) node[anchor=mid] {\# customers};
\draw (1,-0.5) node[anchor=mid] {$0$};
\draw (2,-0.5) node[anchor=mid] {$1$};
\draw (3,-0.5) node[anchor=mid] {$2$};
\draw (4,-0.5) node[anchor=mid] {$3$};
\draw (5,-0.5) node[anchor=mid] {$4$};
\draw (6,-0.5) node[anchor=mid] {$5$};
\draw (-0.5,0) node[anchor=mid] {$0$};
\draw (-0.5,1) node[anchor=mid] {$10$};
\draw (-0.5,2) node[anchor=mid] {$20$};
\draw (-0.5,3) node[anchor=mid] {$30$};
\draw (-0.5,4) node[anchor=mid] {$40$};
\draw (-0.5,5) node[anchor=mid] {$50$};
\draw[line width=0.5mm,gray,fill=gray] (0.60,0) rectangle ++(0.8,0.0);\draw[line width=0.5mm,gray,fill=gray] (1.60,0) rectangle ++(0.8,1.0);\draw[line width=0.5mm,gray,fill=gray] (2.60,0) rectangle ++(0.8,1.1);\draw[line width=0.5mm,gray,fill=gray] (3.60,0) rectangle ++(0.8,1.2);\draw[line width=0.5mm,gray,fill=gray] (4.60,0) rectangle ++(0.8,1.3);\draw[line width=0.5mm,gray,fill=gray] (5.60,0) rectangle ++(0.8,0.4);
\end{tikzpicture}
            
         \end{subfigure}
        \label{fig:histogram-identical-large-fixed-0}
    \end{minipage}
    \begin{minipage}{.99\textwidth}
        \centering
         \caption{Number of customers per number of captures for the instance with different rewards, long rankings, constant spawning demand, and no penalties.} 
        \vspace{0.3cm}
         \begin{subfigure}[b]{0.32\textwidth}
            \centering
            \caption{Optimal solution for 1 facility.}
            \begin{tikzpicture}[scale=.4, every node/.style={scale=.4}]
\draw[line width=0.5mm,thick,->] (0,0) -- (7.5,0);
\draw[line width=0.5mm,thick,->] (0,0) -- (0,5.5);
\draw (7.5,-0.5) node[anchor=mid] {\# captures};
\draw (0,6) node[anchor=mid] {\# customers};
\draw (1,-0.5) node[anchor=mid] {$0$};
\draw (2,-0.5) node[anchor=mid] {$1$};
\draw (3,-0.5) node[anchor=mid] {$2$};
\draw (4,-0.5) node[anchor=mid] {$3$};
\draw (5,-0.5) node[anchor=mid] {$4$};
\draw (6,-0.5) node[anchor=mid] {$5$};
\draw (-0.5,0) node[anchor=mid] {$0$};
\draw (-0.5,1) node[anchor=mid] {$10$};
\draw (-0.5,2) node[anchor=mid] {$20$};
\draw (-0.5,3) node[anchor=mid] {$30$};
\draw (-0.5,4) node[anchor=mid] {$40$};
\draw (-0.5,5) node[anchor=mid] {$50$};
\draw[line width=0.5mm,gray,fill=gray] (0.60,0) rectangle ++(0.8,2.0);\draw[line width=0.5mm,gray,fill=gray] (1.60,0) rectangle ++(0.8,3.0);\draw[line width=0.5mm,gray,fill=gray] (2.60,0) rectangle ++(0.8,0.0);\draw[line width=0.5mm,gray,fill=gray] (3.60,0) rectangle ++(0.8,0.0);\draw[line width=0.5mm,gray,fill=gray] (4.60,0) rectangle ++(0.8,0.0);\draw[line width=0.5mm,gray,fill=gray] (5.60,0) rectangle ++(0.8,0.0);
\end{tikzpicture}
            
         \end{subfigure}
         \begin{subfigure}[b]{0.32\textwidth}
            \centering
            \caption{Optimal solution for 3 facilities.}
            \begin{tikzpicture}[scale=.4, every node/.style={scale=.4}]
\draw[line width=0.5mm,thick,->] (0,0) -- (7.5,0);
\draw[line width=0.5mm,thick,->] (0,0) -- (0,5.5);
\draw (7.5,-0.5) node[anchor=mid] {\# captures};
\draw (0,6) node[anchor=mid] {\# customers};
\draw (1,-0.5) node[anchor=mid] {$0$};
\draw (2,-0.5) node[anchor=mid] {$1$};
\draw (3,-0.5) node[anchor=mid] {$2$};
\draw (4,-0.5) node[anchor=mid] {$3$};
\draw (5,-0.5) node[anchor=mid] {$4$};
\draw (6,-0.5) node[anchor=mid] {$5$};
\draw (-0.5,0) node[anchor=mid] {$0$};
\draw (-0.5,1) node[anchor=mid] {$10$};
\draw (-0.5,2) node[anchor=mid] {$20$};
\draw (-0.5,3) node[anchor=mid] {$30$};
\draw (-0.5,4) node[anchor=mid] {$40$};
\draw (-0.5,5) node[anchor=mid] {$50$};
\draw[line width=0.5mm,gray,fill=gray] (0.60,0) rectangle ++(0.8,0.5);\draw[line width=0.5mm,gray,fill=gray] (1.60,0) rectangle ++(0.8,4.0);\draw[line width=0.5mm,gray,fill=gray] (2.60,0) rectangle ++(0.8,0.5);\draw[line width=0.5mm,gray,fill=gray] (3.60,0) rectangle ++(0.8,0.0);\draw[line width=0.5mm,gray,fill=gray] (4.60,0) rectangle ++(0.8,0.0);\draw[line width=0.5mm,gray,fill=gray] (5.60,0) rectangle ++(0.8,0.0);
\end{tikzpicture}
            
         \end{subfigure}
        \begin{subfigure}[b]{0.32\textwidth}
            \centering
            \caption{Optimal solution for 5 facilities.}
            \begin{tikzpicture}[scale=.4, every node/.style={scale=.4}]
\draw[line width=0.5mm,thick,->] (0,0) -- (7.5,0);
\draw[line width=0.5mm,thick,->] (0,0) -- (0,5.5);
\draw (7.5,-0.5) node[anchor=mid] {\# captures};
\draw (0,6) node[anchor=mid] {\# customers};
\draw (1,-0.5) node[anchor=mid] {$0$};
\draw (2,-0.5) node[anchor=mid] {$1$};
\draw (3,-0.5) node[anchor=mid] {$2$};
\draw (4,-0.5) node[anchor=mid] {$3$};
\draw (5,-0.5) node[anchor=mid] {$4$};
\draw (6,-0.5) node[anchor=mid] {$5$};
\draw (-0.5,0) node[anchor=mid] {$0$};
\draw (-0.5,1) node[anchor=mid] {$10$};
\draw (-0.5,2) node[anchor=mid] {$20$};
\draw (-0.5,3) node[anchor=mid] {$30$};
\draw (-0.5,4) node[anchor=mid] {$40$};
\draw (-0.5,5) node[anchor=mid] {$50$};
\draw[line width=0.5mm,gray,fill=gray] (0.60,0) rectangle ++(0.8,0.1);\draw[line width=0.5mm,gray,fill=gray] (1.60,0) rectangle ++(0.8,2.8);\draw[line width=0.5mm,gray,fill=gray] (2.60,0) rectangle ++(0.8,1.7);\draw[line width=0.5mm,gray,fill=gray] (3.60,0) rectangle ++(0.8,0.3);\draw[line width=0.5mm,gray,fill=gray] (4.60,0) rectangle ++(0.8,0.1);\draw[line width=0.5mm,gray,fill=gray] (5.60,0) rectangle ++(0.8,0.0);
\end{tikzpicture}
            
         \end{subfigure}
        \label{fig:histogram-inversely-large-fixed-0}
    \end{minipage}
    \begin{minipage}{.99\textwidth}
       \centering
       \caption{Number of customers per number of captures for the instance with identical rewards, long rankings, constant spawning demands, and penalties.} 
      \vspace{0.3cm}
       \begin{subfigure}[b]{0.32\textwidth}
          \centering
          \caption{Optimal solution for 1 facility.}
          \begin{tikzpicture}[scale=.4, every node/.style={scale=.4}]
\draw[line width=0.5mm,thick,->] (0,0) -- (7.5,0);
\draw[line width=0.5mm,thick,->] (0,0) -- (0,5.5);
\draw (7.5,-0.5) node[anchor=mid] {\# captures};
\draw (0,6) node[anchor=mid] {\# customers};
\draw (1,-0.5) node[anchor=mid] {$0$};
\draw (2,-0.5) node[anchor=mid] {$1$};
\draw (3,-0.5) node[anchor=mid] {$2$};
\draw (4,-0.5) node[anchor=mid] {$3$};
\draw (5,-0.5) node[anchor=mid] {$4$};
\draw (6,-0.5) node[anchor=mid] {$5$};
\draw (-0.5,0) node[anchor=mid] {$0$};
\draw (-0.5,1) node[anchor=mid] {$10$};
\draw (-0.5,2) node[anchor=mid] {$20$};
\draw (-0.5,3) node[anchor=mid] {$30$};
\draw (-0.5,4) node[anchor=mid] {$40$};
\draw (-0.5,5) node[anchor=mid] {$50$};
\draw[line width=0.5mm,gray,fill=gray] (0.60,0) rectangle ++(0.8,1.4);\draw[line width=0.5mm,gray,fill=gray] (1.60,0) rectangle ++(0.8,2.0);\draw[line width=0.5mm,gray,fill=gray] (2.60,0) rectangle ++(0.8,1.3);\draw[line width=0.5mm,gray,fill=gray] (3.60,0) rectangle ++(0.8,0.3);\draw[line width=0.5mm,gray,fill=gray] (4.60,0) rectangle ++(0.8,0.0);\draw[line width=0.5mm,gray,fill=gray] (5.60,0) rectangle ++(0.8,0.0);
\end{tikzpicture}
          
       \end{subfigure}
      \hfill
       \begin{subfigure}[b]{0.32\textwidth}
          \centering
          \caption{Optimal solution for 3 facilities.}
          \begin{tikzpicture}[scale=.4, every node/.style={scale=.4}]
\draw[line width=0.5mm,thick,->] (0,0) -- (7.5,0);
\draw[line width=0.5mm,thick,->] (0,0) -- (0,5.5);
\draw (7.5,-0.5) node[anchor=mid] {\# captures};
\draw (0,6) node[anchor=mid] {\# customers};
\draw (1,-0.5) node[anchor=mid] {$0$};
\draw (2,-0.5) node[anchor=mid] {$1$};
\draw (3,-0.5) node[anchor=mid] {$2$};
\draw (4,-0.5) node[anchor=mid] {$3$};
\draw (5,-0.5) node[anchor=mid] {$4$};
\draw (6,-0.5) node[anchor=mid] {$5$};
\draw (-0.5,0) node[anchor=mid] {$0$};
\draw (-0.5,1) node[anchor=mid] {$10$};
\draw (-0.5,2) node[anchor=mid] {$20$};
\draw (-0.5,3) node[anchor=mid] {$30$};
\draw (-0.5,4) node[anchor=mid] {$40$};
\draw (-0.5,5) node[anchor=mid] {$50$};
\draw[line width=0.5mm,gray,fill=gray] (0.60,0) rectangle ++(0.8,0.1);\draw[line width=0.5mm,gray,fill=gray] (1.60,0) rectangle ++(0.8,1.6);\draw[line width=0.5mm,gray,fill=gray] (2.60,0) rectangle ++(0.8,0.3);\draw[line width=0.5mm,gray,fill=gray] (3.60,0) rectangle ++(0.8,1.1);\draw[line width=0.5mm,gray,fill=gray] (4.60,0) rectangle ++(0.8,1.2);\draw[line width=0.5mm,gray,fill=gray] (5.60,0) rectangle ++(0.8,0.7);
\end{tikzpicture}
          
       \end{subfigure}
      \hfill
      \begin{subfigure}[b]{0.32\textwidth}
          \centering
          \caption{Optimal solution for 5 facilities.}
          \begin{tikzpicture}[scale=.4, every node/.style={scale=.4}]
\draw[line width=0.5mm,thick,->] (0,0) -- (7.5,0);
\draw[line width=0.5mm,thick,->] (0,0) -- (0,5.5);
\draw (7.5,-0.5) node[anchor=mid] {\# captures};
\draw (0,6) node[anchor=mid] {\# customers};
\draw (1,-0.5) node[anchor=mid] {$0$};
\draw (2,-0.5) node[anchor=mid] {$1$};
\draw (3,-0.5) node[anchor=mid] {$2$};
\draw (4,-0.5) node[anchor=mid] {$3$};
\draw (5,-0.5) node[anchor=mid] {$4$};
\draw (6,-0.5) node[anchor=mid] {$5$};
\draw (-0.5,0) node[anchor=mid] {$0$};
\draw (-0.5,1) node[anchor=mid] {$10$};
\draw (-0.5,2) node[anchor=mid] {$20$};
\draw (-0.5,3) node[anchor=mid] {$30$};
\draw (-0.5,4) node[anchor=mid] {$40$};
\draw (-0.5,5) node[anchor=mid] {$50$};
\draw[line width=0.5mm,gray,fill=gray] (0.60,0) rectangle ++(0.8,0.0);\draw[line width=0.5mm,gray,fill=gray] (1.60,0) rectangle ++(0.8,0.9);\draw[line width=0.5mm,gray,fill=gray] (2.60,0) rectangle ++(0.8,0.0);\draw[line width=0.5mm,gray,fill=gray] (3.60,0) rectangle ++(0.8,0.0);\draw[line width=0.5mm,gray,fill=gray] (4.60,0) rectangle ++(0.8,1.6);\draw[line width=0.5mm,gray,fill=gray] (5.60,0) rectangle ++(0.8,2.5);
\end{tikzpicture}
          
       \end{subfigure}
      \label{fig:histogram-identical-large-fixed-50}
    \end{minipage}
    \begin{minipage}{.99\textwidth}
      \centering
       \caption{Number of customers per number of captures for the instance with different rewards, long rankings, constant spawning demands, and penalties.} 
      \vspace{0.3cm}
       \begin{subfigure}[b]{0.32\textwidth}
          \centering
          \caption{Optimal solution for 1 facility.}
          \begin{tikzpicture}[scale=.4, every node/.style={scale=.4}]
\draw[line width=0.5mm,thick,->] (0,0) -- (7.5,0);
\draw[line width=0.5mm,thick,->] (0,0) -- (0,5.5);
\draw (7.5,-0.5) node[anchor=mid] {\# captures};
\draw (0,6) node[anchor=mid] {\# customers};
\draw (1,-0.5) node[anchor=mid] {$0$};
\draw (2,-0.5) node[anchor=mid] {$1$};
\draw (3,-0.5) node[anchor=mid] {$2$};
\draw (4,-0.5) node[anchor=mid] {$3$};
\draw (5,-0.5) node[anchor=mid] {$4$};
\draw (6,-0.5) node[anchor=mid] {$5$};
\draw (-0.5,0) node[anchor=mid] {$0$};
\draw (-0.5,1) node[anchor=mid] {$10$};
\draw (-0.5,2) node[anchor=mid] {$20$};
\draw (-0.5,3) node[anchor=mid] {$30$};
\draw (-0.5,4) node[anchor=mid] {$40$};
\draw (-0.5,5) node[anchor=mid] {$50$};
\draw[line width=0.5mm,gray,fill=gray] (0.60,0) rectangle ++(0.8,2.7);\draw[line width=0.5mm,gray,fill=gray] (1.60,0) rectangle ++(0.8,1.0);\draw[line width=0.5mm,gray,fill=gray] (2.60,0) rectangle ++(0.8,0.0);\draw[line width=0.5mm,gray,fill=gray] (3.60,0) rectangle ++(0.8,0.0);\draw[line width=0.5mm,gray,fill=gray] (4.60,0) rectangle ++(0.8,1.0);\draw[line width=0.5mm,gray,fill=gray] (5.60,0) rectangle ++(0.8,0.3);
\end{tikzpicture}
          
       \end{subfigure}
       \begin{subfigure}[b]{0.32\textwidth}
          \centering
          \caption{Optimal solution for 3 facilities.}
          \begin{tikzpicture}[scale=.4, every node/.style={scale=.4}]
\draw[line width=0.5mm,thick,->] (0,0) -- (7.5,0);
\draw[line width=0.5mm,thick,->] (0,0) -- (0,5.5);
\draw (7.5,-0.5) node[anchor=mid] {\# captures};
\draw (0,6) node[anchor=mid] {\# customers};
\draw (1,-0.5) node[anchor=mid] {$0$};
\draw (2,-0.5) node[anchor=mid] {$1$};
\draw (3,-0.5) node[anchor=mid] {$2$};
\draw (4,-0.5) node[anchor=mid] {$3$};
\draw (5,-0.5) node[anchor=mid] {$4$};
\draw (6,-0.5) node[anchor=mid] {$5$};
\draw (-0.5,0) node[anchor=mid] {$0$};
\draw (-0.5,1) node[anchor=mid] {$10$};
\draw (-0.5,2) node[anchor=mid] {$20$};
\draw (-0.5,3) node[anchor=mid] {$30$};
\draw (-0.5,4) node[anchor=mid] {$40$};
\draw (-0.5,5) node[anchor=mid] {$50$};
\draw[line width=0.5mm,gray,fill=gray] (0.60,0) rectangle ++(0.8,0.3);\draw[line width=0.5mm,gray,fill=gray] (1.60,0) rectangle ++(0.8,1.5);\draw[line width=0.5mm,gray,fill=gray] (2.60,0) rectangle ++(0.8,0.2);\draw[line width=0.5mm,gray,fill=gray] (3.60,0) rectangle ++(0.8,0.6);\draw[line width=0.5mm,gray,fill=gray] (4.60,0) rectangle ++(0.8,1.2);\draw[line width=0.5mm,gray,fill=gray] (5.60,0) rectangle ++(0.8,1.2);
\end{tikzpicture}
          
       \end{subfigure}
      \begin{subfigure}[b]{0.32\textwidth}
          \centering
          \caption{Optimal solution for 5 facilities.}
          \begin{tikzpicture}[scale=.4, every node/.style={scale=.4}]
\draw[line width=0.5mm,thick,->] (0,0) -- (7.5,0);
\draw[line width=0.5mm,thick,->] (0,0) -- (0,5.5);
\draw (7.5,-0.5) node[anchor=mid] {\# captures};
\draw (0,6) node[anchor=mid] {\# customers};
\draw (1,-0.5) node[anchor=mid] {$0$};
\draw (2,-0.5) node[anchor=mid] {$1$};
\draw (3,-0.5) node[anchor=mid] {$2$};
\draw (4,-0.5) node[anchor=mid] {$3$};
\draw (5,-0.5) node[anchor=mid] {$4$};
\draw (6,-0.5) node[anchor=mid] {$5$};
\draw (-0.5,0) node[anchor=mid] {$0$};
\draw (-0.5,1) node[anchor=mid] {$10$};
\draw (-0.5,2) node[anchor=mid] {$20$};
\draw (-0.5,3) node[anchor=mid] {$30$};
\draw (-0.5,4) node[anchor=mid] {$40$};
\draw (-0.5,5) node[anchor=mid] {$50$};
\draw[line width=0.5mm,gray,fill=gray] (0.60,0) rectangle ++(0.8,0.1);\draw[line width=0.5mm,gray,fill=gray] (1.60,0) rectangle ++(0.8,0.8);\draw[line width=0.5mm,gray,fill=gray] (2.60,0) rectangle ++(0.8,0.0);\draw[line width=0.5mm,gray,fill=gray] (3.60,0) rectangle ++(0.8,0.0);\draw[line width=0.5mm,gray,fill=gray] (4.60,0) rectangle ++(0.8,1.3);\draw[line width=0.5mm,gray,fill=gray] (5.60,0) rectangle ++(0.8,2.8);
\end{tikzpicture}
          
       \end{subfigure}
      \label{fig:histogram-inversely-large-fixed-50}
    \end{minipage}
\end{figure}

\subsubsection{Structure of Optimal Solutions.}

We now investigate the structure of optimal solutions for different instance attributes.
More specifically, we analyze the quality of service perceived by customers, quantified as the number of captures throughout the planning horizon.
To this end, we consider $24$ instances of the benchmark with $|\mathcal{I}| = 50$ locations, $|\mathcal{J}| = 50$ customers, and $|\mathcal{T}| = 5$ periods, for which we are able to obtain the optimal solution with our solution methods.
Figure~\ref{fig:histogram-identical-large-fixed-0} (respectively, Figure~\ref{fig:histogram-inversely-large-fixed-0}) presents the number of customers per number of captures throughout the planning horizon for an instance with long rankings, constant spawning demands, and identical (respectively, different) rewards.
We omit these figures for instances with sparse spawning demands and short rankings because they do not provide additional insights.

First, we can see that optimal solutions composed of larger number of facilities tend to capture customers more often throughout the planning horizon, as one could naturally expect.
In this sense, increasing the number of facilities is a trivial way to increase the quality of service perceived by customers.
Second, we note that optimal solutions of the instance with identical rewards serve customers more often than those of the instance with different rewards.
This result indicates that the provider prefers to capture each customer through some  key combination of location and time period in the optimal solution of instances with different rewards.
We remark that this behaviour stems from the NP-hardness of the \ourproblem{}, as the proof of Theorem~\ref{thm:nphard-single} relies on creating a \ourproblem{} instance where customers should only be captured once through a key combination of location and time period.

Although not an issue for the provider, customers may experience a low quality of service because they are only visited a couple of times in the optimal solutions whenever rewards are, in fact, different.
One strategy to prevent the provider from visiting customers only at these key time periods throughout the planning horizon is to penalize unmet demand.
We follow the penalty structure explained in Section~\ref{sec:ProblemExtensions} and set $p_{j} = 50, \forall j \in \mathcal{J}$ to discourage the provider from visiting customers too rarely throughout the planning horizon.
Figures \ref{fig:histogram-identical-large-fixed-50} and \ref{fig:histogram-inversely-large-fixed-50} have analogous meaning to Figures \ref{fig:histogram-identical-large-fixed-0} and Figure~\ref{fig:histogram-inversely-large-fixed-0}, respectively.
Intuitively, the introduction of penalties force the provider to serve customers more frequently throughout the planning horizon, since spawning demands not served as soon as they appear barely contribute to the total reward and, consequently, total profit of the provider.
In this sense, if applicable, government agencies can make use of these penalties to regulate service providers, with the utmost goal of maintaining a better service level for customers.

\section{Conclusion}
\label{sec:conclusion}

We investigate a novel multi-period deterministic location problem, where the decision maker relocates facilities over time to capture cumulative customer demand.
This paper makes practical, methodological and theoretical contributions. On the practical front, we model demand behaviour that, despite its practical relevance, has not received much attention in the literature. Indeed, in our computational experiments, ignoring cumulative demand (or, alternatively, making myopic decisions) has resulted in an average loss of 40\% (respectively, 6\%) of the optimal reward, leaving considerable room for improvement.
We also draw managerial insights on the structure of optimal solutions under cumulative customer demand, particularly when it comes to the service quality perceived by customers throughout the planning horizon. 

From a methodological perspective, we model this problem as mixed-integer programming formulations and show that one of them \textit{(i)} provides a tighter continuous relaxation and \textit{(ii)} allows for more complex demand behaviour.
On average, the tighter formulation is solved to optimality $53\%$ faster than the alternative formulation, and proves optimality gaps $57\%$ smaller within the same time limit.
We then devise an exact Benders decomposition for the tighter formulation, which is five times faster, on average, than solving the tighter formulation directly, and propose an analytical procedure to generate optimality cuts for the special case with a single facility, which yields $35\%$ faster solution times than obtaining optimality cuts through Gurobi.

On the theoretical front, we focus on the special case with a single facility, and identify which problem characteristics reduce or increase its computational complexity.
Naturally, these negative theoretical results extend to the general case.
We also present a 2-approximate algorithm for the special case with a single facility and identical rewards.
Although this algorithm is a heuristic for the general case, it still finds reasonably high-quality solutions for our benchmark (on average, within $2\%$ of the optimal solution), thus being an theoretically-inspired interesting alternative for tackling large-scale instances.

Given the potential relevance of modelling cumulative demand in other application contexts, we hope that the here provided modelling techniques, solution methods and theoretical insights will be useful to more realistically model and solve such planning problems.
Future work includes studying our planning problem  in a duopoly, where the provider competes over customers with a competitor that can provide the same service.
This setting may be faced by multiple real-world applications, and may considerably impact the structure of the optimal location policy implemented by the service provider.

\section*{Acknowledgements} \label{sec:acknowledgments}

This work was funded by the FRQNT Doctoral Scholarship No. B2X-328911, the FRQ-IVADO Research Chair in
Data Science for Combinatorial Game Theory, and the NSERC Grants No. 2017-05224 and 2024-04051. This research was also enabled
in part by support provided by \href{https://www.calculquebec.ca}{Calcul Québec} and the \href{https://alliancecan.ca/}{Digital Research Alliance of Canada}.

\bibliography{references}

\begin{thebibliography}{}

\bibitem[Alizadeh et~al., 2021]{alizadehMultiPeriodMaximalCovering2021}
Alizadeh, R., Nishi, T., Bagherinejad, J., and Bashiri, M. (2021).
\newblock Multi-{{Period Maximal Covering Location Problem}} with {{Capacitated
  Facilities}} and {{Modules}} for {{Natural Disaster Relief Services}}.
\newblock {\em Applied Sciences}, 11(1):397.

\bibitem[Archetti et~al., 2007]{archettiBranchCutAlgorithmVendorManaged2007}
Archetti, C., Bertazzi, L., Laporte, G., and Speranza, M.~G. (2007).
\newblock A {{Branch-and-Cut Algorithm}} for a {{Vendor-Managed
  Inventory-Routing Problem}}.
\newblock {\em Transportation Science}, 41(3):382--391.

\bibitem[Ausiello et~al., 1980]{ausiello1980structure}
Ausiello, G., D'Atri, A., and Protasi, M. (1980).
\newblock Structure preserving reductions among convex optimization problems.
\newblock {\em Journal of Computer and System Sciences}, 21(1):136--153.

\bibitem[Ballou, 1968]{ballouDynamicWarehouseLocation1968}
Ballou, R.~H. (1968).
\newblock Dynamic {{Warehouse Location Analysis}}.
\newblock {\em Journal of Marketing Research}, 5(3):271--276.

\bibitem[Baron et~al., 2011]{baron2011facility}
Baron, O., Milner, J., and Naseraldin, H. (2011).
\newblock Facility location: {{A}} robust optimization approach.
\newblock {\em Production and Operations Management}, 20(5):772--785.

\bibitem[Basciftci et~al., 2021]{basciftciDistributionallyRobustFacility2021}
Basciftci, B., Ahmed, S., and Shen, S. (2021).
\newblock Distributionally robust facility location problem under
  decision-dependent stochastic demand.
\newblock {\em European Journal of Operational Research}, 292(2):548--561.

\bibitem[Benders, 1962]{bendersPartitioningProceduresSolving1962}
Benders, J.~F. (1962).
\newblock Partitioning procedures for solving mixed-variables programming
  problems.
\newblock {\em Numerische mathematik}, 4(1):238--252.

\bibitem[Berbeglia et~al., 2022]{berbegliaComparativeEmpiricalStudy2022}
Berbeglia, G., Garassino, A., and Vulcano, G. (2022).
\newblock A {{Comparative Empirical Study}} of {{Discrete Choice Models}} in
  {{Retail Operations}}.
\newblock {\em Management Science}, 68(6):4005--4023.

\bibitem[Bierlaire, 1998]{bierlaireDiscreteChoiceModels1998}
Bierlaire, M. (1998).
\newblock Discrete choice models.
\newblock In Labb{\'e}, M., Laporte, G., Tanczos, K., and Toint, P., editors,
  {\em Operations Research and Decision Aid Methodologies in Traffic and
  Transportation Management}, pages 203--227. Springer.

\bibitem[B{\"u}sing et~al., 2021]{busingRobustStrategicPlanning2021}
B{\"u}sing, C., Comis, M., Schmidt, E., and Streicher, M. (2021).
\newblock Robust strategic planning for mobile medical units with steerable and
  unsteerable demands.
\newblock {\em European Journal of Operational Research}, 295(1):34--50.

\bibitem[C{\'a}novas et~al., 2007]{canovasStrengthenedFormulationSimple2007}
C{\'a}novas, L., Garc{\'i}a, S., Labb{\'e}, M., and Mar{\'i}n, A. (2007).
\newblock A strengthened formulation for the simple plant location problem with
  order.
\newblock {\em Operations Research Letters}, 35(2):141--150.

\bibitem[Clothiers, 2024]{retail1}
Clothiers, M. (2024).
\newblock Tour schedule.
\newblock \url{https://www.maxwellsclothiers.com/tour-schedule}.
\newblock Accessed: 30 April 2024.

\bibitem[Cordeau et~al., 2019]{cordeauBendersDecompositionVery2019}
Cordeau, J.-F., Furini, F., and Ljubi{\'c}, I. (2019).
\newblock Benders decomposition for very large scale partial set covering and
  maximal covering location problems.
\newblock {\em European Journal of Operational Research}, 275(3):882--896.

\bibitem[Cornu{\'e}jols et~al., 1983]{cornuejols1983uncapicitated}
Cornu{\'e}jols, G., Nemhauser, G., and Wolsey, L. (1983).
\newblock The uncapacitated facility location problem.
\newblock Technical report, {Cornell University Operations Research and
  Industrial Engineering}.

\bibitem[Daneshvar et~al., 2023]{daneshvarTwostageStochasticPostdisaster2023}
Daneshvar, M., Jena, S.~D., and Rei, W. (2023).
\newblock A two-stage stochastic post-disaster humanitarian supply chain
  network design problem.
\newblock {\em Computers \& Industrial Engineering}, 183:109459.

\bibitem[Dubinski, 2021]{healthcare2}
Dubinski, K. (2021).
\newblock This mobile medical clinic has helped over 500 people in london
  ontario.
\newblock
  \url{https://www.cbc.ca/news/canada/london/this-mobile-medical-clinic-has-helped-over-500-people-in-london-ont-1.6112957}.
\newblock Accessed: 30 April 2024.

\bibitem[Farias et~al., 2013]{fariasNonparametricApproachModeling2013}
Farias, V.~F., Jagabathula, S., and Shah, D. (2013).
\newblock A {{Nonparametric Approach}} to {{Modeling Choice}} with {{Limited
  Data}}.
\newblock {\em Management Science}, 59(2):305--322.

\bibitem[Garey and Johnson, 1979]{garey1979computers}
Garey, M.~R. and Johnson, D.~S. (1979).
\newblock {\em Computers and Intractability: {{A}} Guide to the Theory of
  {{NP-Completeness}}}.
\newblock W. H. Freeman, New York.

\bibitem[Gunawardane, 1982]{gunawardaneDynamicVersionsSet1982}
Gunawardane, G. (1982).
\newblock Dynamic versions of set covering type public facility location
  problems.
\newblock {\em European Journal of Operational Research}, 10(2):190--195.

\bibitem[Hanjoul and Peeters, 1987]{hanjoulFacilityLocationProblem1987}
Hanjoul, P. and Peeters, D. (1987).
\newblock A facility location problem with clients' preference orderings.
\newblock {\em Regional Science and Urban Economics}, 17(3):451--473.

\bibitem[Hastad, 1996]{hastadCliqueHardApproximate1996}
Hastad, J. (1996).
\newblock Clique is hard to approximate within n/sup 1-{$\varepsilon$}/.
\newblock In Shor, P.~W., editor, {\em Proceedings of 37th {{Conference}} on
  {{Foundations}} of {{Computer Science}}}, pages 627--636, Burlington, VT,
  USA. IEEE Comput. Soc. Press.

\bibitem[Hazan et~al., 2006]{hazanComplexityApproximatingKset2006}
Hazan, E., Safra, S., and Schwartz, O. (2006).
\newblock On the complexity of approximating k-set packing.
\newblock {\em computational complexity}, 15(1):20--39.

\bibitem[Hu et~al., 2025]{huRobustFacilityLocation2025}
Hu, H., Tang, J., and Tian, T. (2025).
\newblock Robust facility location and protection under facility disruptions
  with decision-dependent uncertainty.
\newblock {\em International Journal of Production Economics}, 282:109558.

\bibitem[Jena et~al., 2015]{jenaDynamicFacilityLocation2015}
Jena, S.~D., Cordeau, J.-F., and Gendron, B. (2015).
\newblock Dynamic {{Facility Location}} with {{Generalized Modular
  Capacities}}.
\newblock {\em Transportation Science}, 49(3):484--499.

\bibitem[Jena et~al., 2020]{jenaPartiallyRankedChoice2020}
Jena, S.~D., Lodi, A., Palmer, H., and Sole, C. (2020).
\newblock A {{Partially Ranked Choice Model}} for {{Large-Scale Data-Driven
  Assortment Optimization}}.
\newblock {\em INFORMS Journal on Optimization}, 2(4):297--319.

\bibitem[Karp, 1972]{karp1972reducibility}
Karp, R.~M. (1972).
\newblock Reducibility among {{Combinatorial Problems}}.
\newblock In Miller, R.~E., Thatcher, J.~W., and Bohlinger, J.~D., editors,
  {\em Complexity of {{Computer Computations}}}, pages 85--103. Springer.

\bibitem[Kuhn, 1955]{kuhn1955hungarian}
Kuhn, H.~W. (1955).
\newblock The {{Hungarian}} method for the assignment problem.
\newblock {\em Naval Research Logistics Quarterly}, 2(1-2):83--97.

\bibitem[Lamontagne et~al., 2023]{lamontagneOptimisingElectricVehicle2023}
Lamontagne, S., Carvalho, M., Frejinger, E., Gendron, B., Anjos, M.~F., and
  Atallah, R. (2023).
\newblock Optimising {{Electric Vehicle Charging Station Placement Using
  Advanced Discrete Choice Models}}.
\newblock {\em INFORMS Journal on Computing}, 35(5):1195--1213.

\bibitem[Laporte et~al., 2019]{laporteLocationScience2019}
Laporte, G., Nickel, S., and {Saldanha da Gama}, F. (2019).
\newblock {\em Location {{Science}}}.
\newblock Springer.

\bibitem[Lin et~al., 2022]{linLocatingFacilitiesCompetition2022}
Lin, Y.~H., Tian, Q., and Zhao, Y. (2022).
\newblock Locating facilities under competition and market expansion:
  {{Formulation}}, optimization, and implications.
\newblock {\em Production and Operations Management}, 31(7):3021--3042.

\bibitem[Magnanti and Wong, 1981]{magnantiAcceleratingBendersDecomposition1981}
Magnanti, T.~L. and Wong, R.~T. (1981).
\newblock Accelerating {{Benders Decomposition}}: {{Algorithmic Enhancement}}
  and {{Model Selection Criteria}}.
\newblock {\em Operations Research}, 29(3):464--484.

\bibitem[Malladi and Muthuraman, 2024]{malladiFacilityLocationProblem2024}
Malladi, V. and Muthuraman, K. (2024).
\newblock Facility {{Location Problem}}: {{Modeling Joint Disruptions Using
  Subordination}}.
\newblock {\em Transportation Science}, 58(5):1016--1032.

\bibitem[Mar{\'i}n et~al., 2018]{marinMultiperiodStochasticCovering2018a}
Mar{\'i}n, A., {Mart{\'i}nez-Merino}, L.~I., {Rodr{\'i}guez-Ch{\'i}a}, A.~M.,
  and {Saldanha-da-Gama}, F. (2018).
\newblock Multi-period stochastic covering location problems: {{Modeling}}
  framework and solution approach.
\newblock {\em European Journal of Operational Research}, 268(2):432--449.

\bibitem[Nica and Moraru, 2020]{nica2020diaspora}
Nica, F. and Moraru, M. (2020).
\newblock Diaspora policies, consular services and social protection for
  {{Romanian}} citizens abroad.
\newblock In Lafleur, J.-M. and Vintila, D., editors, {\em Migration and
  {{Social Protection}} in {{Europe}} and {{Beyond}} ({{Volume}} 2) {{Comparing
  Consular Services}} and {{Diaspora Policies}}}, pages 409--425. Springer.

\bibitem[Nickel and {Saldanha-da-Gama}, 2019]{nickel2019multi}
Nickel, S. and {Saldanha-da-Gama}, F. (2019).
\newblock Multi-period facility location.
\newblock In Laporte, G., Nickel, S., and {Saldanha-da-Gama}, F., editors, {\em
  Location Science}, pages 303--326. Springer.

\bibitem[Nzioka, 2024]{consular1}
Nzioka, T. (2024).
\newblock Kenya rolls out phase three of diaspora mobile consular services.
\newblock
  \url{https://www.the-star.co.ke/news/realtime/2024-03-18-kenya-rolls-out-phase-three-of-diaspora-mobile-consular-services/}.
\newblock Accessed: 30 April 2024.

\bibitem[Qi et~al., 2017]{qiMobileFacilityRouting2017}
Qi, M., Cheng, C., Wang, X., and Rao, W. (2017).
\newblock Mobile {{Facility Routing Problem}} with {{Service-Time-related
  Demand}}.
\newblock In {\em 2017 {{International Conference}} on {{Service Systems}} and
  {{Service Management}}}, pages 1--6, Dalian, China. IEEE.

\bibitem[Qi et~al., 2022]{qiSequentialCompetitiveFacility2022}
Qi, M., Jiang, R., and Shen, S. (2022).
\newblock Sequential {{Competitive Facility Location}}: {{Exact}} and
  {{Approximate Algorithms}}.
\newblock {\em Operations Research}, 72(1):300--316.

\bibitem[Rahmaniani et~al., 2017]{rahmanianiBendersDecompositionAlgorithm2017}
Rahmaniani, R., Crainic, T.~G., Gendreau, M., and Rei, W. (2017).
\newblock The {{Benders}} decomposition algorithm: {{A}} literature review.
\newblock {\em European Journal of Operational Research}, 259(3):801--817.

\bibitem[Rosenbaum et~al., 2021]{rosenbaumBenefitsPitfallsContemporary2021}
Rosenbaum, M.~S., Edwards, K., and Ramirez, G.~C. (2021).
\newblock The benefits and pitfalls of contemporary pop-up shops.
\newblock {\em Business Horizons}, 64(1):93--106.

\bibitem[Sahyouni et~al., 2007]{sahyouniFacilityLocationModel2007}
Sahyouni, K., Savaskan, R.~C., and Daskin, M.~S. (2007).
\newblock A {{Facility Location Model}} for {{Bidirectional Flows}}.
\newblock {\em Transportation Science}, 41(4):484--499.

\bibitem[Schuurman and Woeginger, 2008]{schuurman2001approximation}
Schuurman, P. and Woeginger, G.~J. (2008).
\newblock Approximation schemes - a tutorial.
\newblock {\em Working paper}.

\bibitem[Sweeney and Tatham, 1976]{sweeneyImprovedLongRunModel1976}
Sweeney, D.~J. and Tatham, R.~L. (1976).
\newblock An {{Improved Long-Run Model}} for {{Multiple Warehouse Location}}.
\newblock {\em Management Science}, 22(7):748--758.

\bibitem[Van~Roy and Erlenkotter, 1982]{vanroyDualBasedProcedureDynamic1982}
Van~Roy, T.~J. and Erlenkotter, D. (1982).
\newblock A {{Dual-Based Procedure}} for {{Dynamic Facility Location}}.
\newblock {\em Management Science}, 28(10):1091--1105.

\bibitem[{van Ryzin} and Vulcano, 2015]{vanryzinMarketDiscoveryAlgorithm2015}
{van Ryzin}, G. and Vulcano, G. (2015).
\newblock A {{Market Discovery Algorithm}} to {{Estimate}} a {{General Class}}
  of {{Nonparametric Choice Models}}.
\newblock {\em Management Science}, 61(2):281--300.

\bibitem[Vatsa and Jayaswal, 2021]{vatsaCapacitatedMultiperiodMaximal2021}
Vatsa, A.~K. and Jayaswal, S. (2021).
\newblock Capacitated multi-period maximal covering location problem with
  server uncertainty.
\newblock {\em European Journal of Operational Research}, 289(3):1107--1126.

\bibitem[Wesolowsky, 1973]{wesolowskyDynamicFacilityLocation1973}
Wesolowsky, G.~O. (1973).
\newblock Dynamic {{Facility Location}}.
\newblock {\em Management Science}, 19(11):1241--1248.

\bibitem[Wesolowsky and Truscott,
  1975]{wesolowskyMultiperiodLocationAllocationProblem1975}
Wesolowsky, G.~O. and Truscott, W.~G. (1975).
\newblock The {{Multiperiod Location-Allocation Problem}} with {{Relocation}}
  of {{Facilities}}.
\newblock {\em Management Science}, 22(1):57--65.

\bibitem[Zhang et~al., 2021]{zhangBendersDecompositionApproach2021}
Zhang, Z., Luo, Z., Baldacci, R., and Lim, A. (2021).
\newblock A {{Benders Decomposition Approach}} for the {{Multivehicle
  Production Routing Problem}} with {{Order-up-to-Level Policy}}.
\newblock {\em Transportation Science}, 55(1):160--178.

\end{thebibliography}

\appendix

\section{Rank-Based Choice Model}
\label{apx:choicemodel}

In this appendix, we explain in detail how we can readily account for a rank-based choice model \citep[see, \eg,][]{fariasNonparametricApproachModeling2013, vanryzinMarketDiscoveryAlgorithm2015, jenaPartiallyRankedChoice2020} in the \ourproblem{}, which performs better than parametric choice models (\eg, the Multinomial Logit Model) in practice \citep{berbegliaComparativeEmpiricalStudy2022}.

Let $\mathcal{I}^{R}, \mathcal{J}^{R}, \mathcal{T}, \boldsymbol{r}^{R}, \boldsymbol{d}^{R}, \boldsymbol{\succ}^{R}, h^{R}$ be a \ourproblem{} instance where we assume that each customer $j$ can be represented by a set of customer profiles $\mathcal{P}_{j}$ indexed by $(j, p)$.
Each customer profile $(j, p)$ has its own ranking $\succ_{jp}$ and a realization probability $\omega_{jp} \in [0, 1]$.
%
%
We can obtain a standard instance $\mathcal{I}^{S}, \mathcal{J}^{S}, \mathcal{T}, \boldsymbol{r}^{S}, \boldsymbol{d}^{S}, \boldsymbol{\succ}^{S}, h^{S}$ of the \ourproblem{} where each pair $(j,p)$ becomes a separate customer.
To be explicit, we set the input parameters of instance $S$ based on instance $R$ as follows: $\mathcal{I}^{S} = \mathcal{I}^{R}$; $\mathcal{J}^{S} = \mathset{(j,p) \mid j \in \mathcal{J}^{R}, p \in \mathcal{P}_{j}}$; $\mathcal{T}^{S} = \mathcal{T}^{R}$; $r^{S}_{i} = r^{R}_{i}, \forall i \in \mathcal{I}^{R}$;  ${d^{t}_{(j,p)}}^{S} = \omega_{jp} {d^{t}_{j}}^{R}, \forall j \in \mathcal{J}^{R}, \forall p \in \mathcal{P}_{j}, \forall t \in \mathcal{T}^{R}$; $h^{S} = h^{R}$.

We can now apply our solution methods to instance $S$ without loss of generality.
In other words, employing rank-based choice models only increases the number of targeted customers in the \ourproblem{} instance.

\section{Linearization Details}
\label{apx:linearization}

In this appendix, we provide the detailed linearization of the SIF.

We employ an additional decision variable $w^{t}_{ij} \in \mathbb{R}^{+}$ to store the demand of customer $j$ captured by location $i$ at time period $t$, and an additional parameter $M^{t}_{j} \in \mathbb{R}^{+}$ to represent a sufficiently large constant for each time period $t$ and each customer $j$.
The linearized version of the SIF writes as follows:
\begin{subequations}
    \label{pgm:linearized-dflp-ccd}
    \begin{align}
        \max_{\boldsymbol{u}, \boldsymbol{c}, \boldsymbol{w}, \boldsymbol{x}, \boldsymbol{y}} \quad  &  \sum_{t \in \mathcal{T}} \sum_{i \in \mathcal{I}} \sum_{j \in \mathcal{J}} r_{i} w^{t}_{ij}
        && \label{eq:linearized-dflp-ccd-obj}\\
        \text{s.t.} \quad
        & \boldsymbol{y} \in \mathcal{Y} 
        && \label{eq:linearized-dflp-ccd-ct1} \\
        & x^{t}_{ij} \leq a_{ij} y^{t}_{i} 
        && \forall i \in \mathcal{I}, \forall j \in \mathcal{J}, \forall t \in \mathcal{T} \label{eq:linearized-dflp-ccd-ct2} \\
        & \sum_{k \in \mathcal{I}} x^{t}_{kj} \geq a_{ij} y^{t}_{i}
        && \forall i \in \mathcal{I}, \forall j \in \mathcal{J}, \forall t \in \mathcal{T} \label{eq:linearized-dflp-ccd-ct3} \\
        & a_{ij} y^{t}_{i} + \sum_{\substack{k \in \mathcal{I}: \\ i \succ_{j} k}} x^{t}_{kj} \leq  a_{ij} &&  \forall i \in \mathcal{I}, \forall j \in \mathcal{J}, \forall t \in \mathcal{T} \label{eq:linearized-dflp-ccd-ct4} \\
        & u^{0}_{j} = 0
        && \forall j \in \mathcal{J}
        \label{eq:linearized-dflp-ccd-ct5} \\
        & c^{t}_{j} = u^{t-1}_{j} + d^{t}_{j}
        && \forall j \in \mathcal{J}, \forall t \in \mathcal{T} \label{eq:linearized-dflp-ccd-ct6} \\
        & u^{t}_{j} = c^{t}_{j} - \sum_{i \in \mathcal{I}} w^{t}_{ij}
        && \forall j \in \mathcal{J}, \forall t \in \mathcal{T} \label{eq:linearized-dflp-ccd-ct7} \\
        & w^{t}_{ij} \leq M^{t}_{j} x^{t}_{ij}
        && \forall i \in \mathcal{I}, \forall j \in \mathcal{J}, \forall t \in \mathcal{T} \label{eq:linearized-dflp-ccd-ct8a} \\
        & w^{t}_{ij} \leq c^{t}_{j} + M^{t}_{j} (1- x^{t}_{ij})
        && \forall i \in \mathcal{I}, \forall j \in \mathcal{J}, \forall t \in \mathcal{T} \label{eq:linearized-dflp-ccd-ct8b} \\
        & w^{t}_{ij} \geq - M^{t}_{j} x^{t}_{ij}
        && \forall i \in \mathcal{I}, \forall j \in \mathcal{J}, \forall t \in \mathcal{T} \label{eq:linearized-dflp-ccd-ct8c} \\
        & w^{t}_{ij} \geq c^{t}_{j} - M^{t}_{j} (1- x^{t}_{ij})
        && \forall i \in \mathcal{I}, \forall j \in \mathcal{J}, \forall t \in \mathcal{T} \label{eq:linearized-dflp-ccd-ct8d} \\
        & u^{t}_{j} \in \mathbb{R}^{+}
        && \forall j \in \mathcal{J}, \forall t \in \mathcal{T}^{S} \label{eq:linearized-dflp-ccd-dm1} \\
        & c^{t}_{j} \in \mathbb{R}^{+}
        && \forall j \in \mathcal{J}, \forall t \in \mathcal{T} \label{eq:linearized-dflp-ccd-dm2} \\
        & w^{t}_{ij} \in \mathbb{R}^{+}
        && \forall i \in \mathcal{I}, \forall j \in \mathcal{J}, \forall t \in \mathcal{T} \label{eq:linearized-dflp-ccd-dm3} \\
        & x^{t}_{ij} \in \mathset{0,1}
        && \forall i \in \mathcal{I}, \forall j \in \mathcal{J}, \forall t \in \mathcal{T} \label{eq:linearized-dflp-ccd-dm4} \\
        & y^{t}_{i} \in \mathset{0,1}
        && \forall i \in \mathcal{I}, \forall t \in \mathcal{T}. \label{eq:linearized-dflp-ccd-dm5}
    \end{align}
\end{subequations}

Continuous variables $w^{t}_{ij}$ linearize bilinear terms $c^{t}_{j} x^{t}_{ij}$ in Objective Function (2a) and Constraints (2h).
In this sense, Objective Function \eqref{eq:linearized-dflp-ccd-obj} and Constraints \eqref{eq:linearized-dflp-ccd-ct7} are straightforward adaptations of their counterparts in the SIF.
Constraints \eqref{eq:linearized-dflp-ccd-ct8a}--\eqref{eq:linearized-dflp-ccd-ct8d} regulate the behaviour of the continuous variable $w^{t}_{ij}$.
If customer $j$ is captured by location $i$ at time period $t$ (\ie, $x^{t}_{ij} = 1$), Constraints \eqref{eq:linearized-dflp-ccd-ct8a} and \eqref{eq:linearized-dflp-ccd-ct8c} become nonrestrictive, whereas Constraints \eqref{eq:linearized-dflp-ccd-ct8b} and \eqref{eq:linearized-dflp-ccd-ct8d} become binding, thus ensuring $w^{t}_{ij} = c^{t}_{j}$.
The exact opposite happens if customer $j$ is not captured by location $i$ at time period $t$ (\ie, $x^{t}_{ij} = 0$), thus ensuring $w^{t}_{ij} = 0$.
Setting tight values for parameters $M^{t}_{j}$ is important to obtain a tighter continuous relaxation, which in turn tends to help off-the-shelf solvers employing branch-and-bound to find optimal solutions and prove optimality faster.
In the computational experiments, we set $M^{t}_{j} = D^{0t}_{j} = \sum_{\substack{s \in \mathcal{T} \\ s > 0 \\ s \leq t}} d^{s}_{j}$.

\section{Formulation Tightness}
\label{apx:tightness}

We present here the proof of Theorem~1.

\textbf{Proof.}
    Let $\overline{DIF}$ and $\overline{SIF}$ be the linear relaxations of the DIF and the SIF, respectively.
    Recall that we employ the linearized version of the SIF presented in Appendix~\ref{apx:linearization}.
    We show that the DIF provides a tighter continuous relaxation than the SIF in two steps.
    First, we prove that the DIF is at least as tight as the SIF by showing that each feasible solution $(\boldsymbol{x}, \boldsymbol{y})$ in the $\overline{DIF}$ has an equivalent feasible solution $(\boldsymbol{u}, \boldsymbol{c}, \boldsymbol{w}, \boldsymbol{x}, \boldsymbol{y})$ in the $\overline{SIF}$ with the same objective value.
    Second, we prove that the DIF is strictly tighter than the SIF by providing an example where the optimal objective value of the $\overline{DIF}$ provides a strictly better bound than the one of the $\overline{DIF}$ (\ie, closer to the objective value of the optimal integer solution).
    
    \textbf{First Step.}
    We first generate a feasible solution $(\boldsymbol{u}, \boldsymbol{c}, \boldsymbol{w}, \boldsymbol{x}, \boldsymbol{y})$ in the $\overline{SIF}$ from a feasible solution $(\boldsymbol{x}, \boldsymbol{y})$ in the $\overline{DIF}$.
    First, we set variables $y^{t}_{i}$ to the same value.
    Then, we set the remaining variables as follows:
    \begin{align}
        & w^{t}_{ij} = \sum_{\substack{\ell{} \in \mathcal{T}^{S}: \\ \ell{} < t}} D^{\ell{}t}_{j} x^{\ell{}t}_{ij} && \forall i \in \mathcal{I}, \forall j \in \mathcal{J}, \forall t \in \mathcal{T} \label{eq:transformation-w} \\
        & x^{t}_{ij} = \sum_{\substack{\ell{} \in \mathcal{T}^{S}: \\ \ell{} < t}} x^{\ell{}t}_{ij} && \forall i \in \mathcal{I}, \forall j \in \mathcal{J}, \forall t \in \mathcal{T} \label{eq:transformation-x} \\
        & c^{t}_{j} = u^{t-1}_{j} + d^{t}_{j} && \forall j \in \mathcal{J}, \forall t \in \mathcal{T} \label{eq:transformation-c} \\ 
        & u^{t}_{j} = c^{t}_{j} - \sum_{i \in \mathcal{I}} w^{t}_{ij}  && \forall j \in \mathcal{J}, \forall t \in \mathcal{T}. \label{eq:transformation-b}
    \end{align}

    This feasible solution in the $\overline{SIF}$ has the same objective value as the one in the $\overline{DIF}$, as replacing values $w^{t}_{ij}$ in Objective Function \eqref{eq:linearized-dflp-ccd-obj} results in Objective Function (1a).
    Constraints \eqref{eq:linearized-dflp-ccd-ct1} are satisfied trivially due to Constraints (1b), Constraints \eqref{eq:linearized-dflp-ccd-ct2}--\eqref{eq:linearized-dflp-ccd-ct4} are satisfied due to Constraints (1c)--(1e) and Equation~\eqref{eq:transformation-x}, and Constraints \eqref{eq:linearized-dflp-ccd-ct5}--\eqref{eq:linearized-dflp-ccd-ct7} are satisfied by construction due to Equations \eqref{eq:transformation-c}--\eqref{eq:transformation-b}.
    We now show that values $w^{t}_{ij}$, built with Equation~\eqref{eq:transformation-w}, respect Constraints \eqref{eq:linearized-dflp-ccd-ct8a}--\eqref{eq:linearized-dflp-ccd-ct8d} by showing that they respect $w^{t}_{ij} = c^{t}_{j} x^{t}_{ij}$.

    Assume for the sake of contradiction that there is a customer $j$, a location $i$, and a time period $t$ such that $w^{t}_{ij} \neq c^{t}_{j} x^{t}_{ij}$.
    If $w^{t}_{ij} < c^{t}_{j} x^{t}_{ij}$, then $\sum_{\substack{\ell{} \in \mathcal{T}^{S} : \\ \ell{} < t}} D^{\ell{}t}_{j} x^{\ell{}t}_{ij} <  c^{t}_{j} \sum_{\substack{\ell{} \in \mathcal{T}^{S} : \\ \ell{} < t}} x^{\ell{}t}_{ij}$ due to Equations~\eqref{eq:transformation-w}--\eqref{eq:transformation-x}.
    This implies $D^{\ell{}t}_{j} < c^{t}_{j} \, \forall \ell{} \in \mathcal{T}^{S} : \ell{} < t$, which is an absurd because $c^{t}_{j} \leq D^{0t}_{j}$ due to Equation~\eqref{eq:transformation-c}.
    If $w^{t}_{ij} > c^{t}_{j} x^{t}_{ij}$, then $\sum_{\substack{\ell{} \in \mathcal{T}^{S} : \\ \ell{} < t}} D^{\ell{}t}_{j} x^{\ell{}t}_{ij} >  c^{t}_{j} \sum_{\substack{\ell{} \in \mathcal{T}^{S} : \\ \ell{} < t}}  x^{\ell{}t}_{ij}$ due to Equation~\eqref{eq:transformation-w}--\eqref{eq:transformation-x}.
    This implies $D^{\ell{}t}_{j} > c^{t}_{j} \, \forall \ell{} \in \mathcal{T}^{S} : \ell{} < t$, which is an absurd because $c^{t}_{j} \geq d^{t}_{j} = D^{(t-1)t}_{j}$ due to Equation~\eqref{eq:transformation-c}.
    Therefore, Constraints~\eqref{eq:linearized-dflp-ccd-ct8a}--\eqref{eq:linearized-dflp-ccd-ct8d} are satisfied.

    \textbf{Second Step.}
    Consider a 1-\ourproblem{} instance with three locations $\mathcal{I} = \mathset{1,2,3}$, two customers $\mathcal{J} = \mathset{A,B}$, and two time periods $\mathcal{T} = \mathset{1,2}$.
    We employ letters instead of numbers to refer to customers for the sake of clarity.
    Locations $1$ and $2$ have a reward of $100$ (\ie, $r_{1} = r_{2} = 100$), whereas location $3$ has a reward of $51$ (\ie, $r_{3} = 51$).
    Customer $A$ is willing to attend locations $1$ and $3$ (\ie, $a_{1A} = a_{3A} = 1, a_{2A} = 0$), whereas customer $B$ is willing to attend locations $2$ and $3$ (\ie, $a_{2B} = a_{3B} = 1, a_{1B} = 0$).
    Both customers have $1$ unit of spawning demand at each time period (\ie, $d^{t}_{j} = 1 \, \forall j \in \mathcal{J}, \forall t \in \mathcal{T}$).

    The optimal (integer) solution is to install the temporary facility at location $1$ at time period $t = 1$ and location $2$ at time period $t = 2$, which gives an optimal (integer) objective value $Z^{\star} = 300$.
    The $\overline{DIF}$ has an optimal objective value of $Z^{1} = 300$ (\ie, the tightest bound possible), whereas the $\overline{SIF}$ has an optimal objective value of $Z^{M} = 302$ (\ie, a looser bound).
    Instead of opening location $2$ at time period $t = 2$ completely (\ie, setting $y^{2}_{2} = 1$), the $\overline{SIF}$ opens \textit{(i)} half of location $2$ ($y^{2}_{2} = \frac{1}{2}$) to capture some demand from customer $B$ with a higher reward, and \textit{(ii)} half of location $3$ ($y^{2}_{3} = \frac{1}{2}$) to capture remaining demand from customers $A$ and $B$.
    Therefore, the DIF is strictly tighter than the SIF, so Theorem~1 holds. \Halmos

\section{Complexity Proofs}
\label{apx:proofs}

In this appendix, we prove the complexity results presented in Section~4.
More specifically, Appendices~\ref{apx:nphard-multiple}--\ref{apx:approximation} present the proofs for Theorems~2--5, respectively.
We recall that we fix $\mathcal{Y} = \mathset{\sum_{i \in \mathcal{I}} y^{t}_{i} \leq h, \forall t \in \mathcal{T}}$,  $f(\boldsymbol{y}) = 0, \forall \boldsymbol{y} \in \mathcal{Y}$, and do not consider problem extensions of Section~3.5 here.

Before proceeding, we define some additional notation useful for some proofs presented in this appendix.
We remark that a location policy $\boldsymbol{y}$ for the 1-\ourproblem{} can be written as $i_{1}, \ldots, i_{T}$ without loss of generality.
We first define the set of captured customers by a location $i$ or a subset of locations $\mathcal{I}^{\prime}$.

\begin{definition}[Captured customer set]
    The set of customers captured by location $i \in \mathcal{I}$ is $\mathcal{J}(i) = \{ j \in \mathcal{J} \mid a_{ij} = 1\}$.
    Similarly, the set of customers captured by a subset of locations $\mathcal{I}^{\prime} \subseteq \mathcal{I}$ is $\mathcal{J}(\mathcal{I}^{\prime}) = \bigcup_{i \in \mathcal{I}^{\prime}} \mathcal{J}(i)$.
\end{definition}

Then, we define the time period where customer $j$ was last captured within a feasible solution $i_{1}, \ldots, i_{T}$, which allow us to compute the accumulated demand of said customer at time period $t$.

\begin{definition}[Time of previous capture]
    Let $i_{1}, \ldots, i_{T}$ be a feasible solution of a 1-\ourproblem{} instance.
    The time of previous capture of customer $j$ in the feasible solution  $i_{1}, \ldots, i_{T}$ before time period $t$ is
    $$
        \tau(j, t \mid i_{1}, \ldots, i_{T}) = \max\{0, t^{\prime} \in \mathcal{T} \mid t^{\prime} < t, j \in \mathcal{J}(i_{t^{\prime}}) \}.
   $$
\end{definition}

Finally, we define the marginal reward brought by location $i_{t}$ at time period $t$ to the total reward of the provider, as well as the total reward of the provider for a feasible solution $i_{1}, \ldots, i_{T}$, as follows.

\begin{definition}[Marginal reward function]
    \label{def:marginal-reward}
    Let $i_{1}, \ldots, i_{T}$ be a feasible solution of a 1-\ourproblem{} instance.
    The marginal reward function of location $i_{t}$ at time period $t$ in the feasible solution $i_{1}, \ldots, i_{T}$ is
    $$
    \rho(i_{t}, t \mid i_{1}, \ldots, i_{T}) =
        r_{i_{t}} \sum_{j \in \mathcal{J}(i_{t})} \sum_{\substack{s \in \mathcal{T}: \\ s > \tau(j, t \mid i_{1}, \ldots, i_{T}) \\ s \leq t}} d^{s}_{j}.
    $$
\end{definition}

\begin{definition}[Total reward function]
    Let $i_{1}, \ldots, i_{T}$ be a feasible solution of a 1-\ourproblem{} instance.
    The total reward function of the feasible solution $i_{1}, \ldots, i_{T}$ is $$z(i_{1}, \ldots, i_{T}) = \sum_{t \in \mathcal{T}} \rho(i_{t}, t \mid i_{1}, \ldots, i_{T}).$$
\end{definition}

\subsection{Proof of Theorem~2}
\label{apx:nphard-multiple}

\textbf{Proof.}
    First, we show that the \ourproblem{} is in NP.
    A decision problem is in NP if the certificate answering the decision question can be verified in polynomial time.
    For the \ourproblem{}, this means verifying if $z(\boldsymbol{y}) \geq Z$ for some location policy $\boldsymbol{y}$.
    Proposition~1 guarantees that $z(\boldsymbol{y})$ can be computed in polynomial time, and it  suffices to check whether $z(\boldsymbol{y}) \geq Z$ to answer the decision question.
    Therefore, it follows that the \ourproblem{} is in NP.

    Second, we show that the \ourproblem{} is NP-hard because it has the classical facility location problem, which is known to be NP-hard \citep[][]{cornuejols1983uncapicitated}, as a special case.
    Consider a \ourproblem{} instance where \textit{(i)} the planning horizon is composed of a single time period (\ie, $|\mathcal{T}| = 1$) and \textit{(ii)} customer rankings align with rewards at time period $t = 1$ (\ie, $i \succ_{j} k$ if and only if $r_{i} > r_{k}, \forall i \in \mathcal{I}, \forall k \in \mathcal{I}, \forall j \in \mathcal{J}$).
    This special case is an instance of  the classical facility location problem, so the \ourproblem{} is, in fact, NP-hard.

\subsection{Proof of Theorem~3}
\label{apx:nphard-single}

We first formalize the decision version of the SPP, and then present the proof.

\vspace{0.3cm}
\boxxx{
{\bf Decision version of the SPP}: \\
{\sc instance}: A finite collection of $n$ sets $\mathcal{C} = \mathset{\mathcal{C}_1, \ldots, \mathcal{C}_n}$, a finite set of $m$ elements $\mathcal{B} = \mathset{B_1, \ldots, B_m}$ appearing in $\mathcal{C}$ and a positive integer $1\leq K \leq |\mathcal{C}|$.\\
{\sc question}: Is there at least $K$ mutually disjoint sets in $\mathcal{C}$?
}
\vspace{0.3cm}

\textbf{Proof.}
    The 1-\ourproblem{} is in NP as the general \ourproblem{} is in NP (see Theorem~2).
    We show that the 1-\ourproblem{} is NP-hard by reducing the SPP to it.
    By showing that if there is a certificate satisfying the SPP decision question, then there is a certificate satisfying the 1-\ourproblem{} decision question (referred to as forward direction) and vice-versa (referred to as backward direction), it holds that the 1-\ourproblem{} is NP-hard.

    \textbf{Reduction.}
    Consider the 1-\ourproblem{} instance built from an SPP instance with $2n$ locations, $m + n$ customers, and $n$ time periods as follows. Each element $B_{j}$ generates an element-customer $j$ with spawning demands $d^{1}_{j} = 1, d^{t}_{j} = 0, \forall t \in \mathcal{T} : t \geq 2$. 
    Each set $\mathcal{C}_{i}$ generates an authentic location $i_{1}$ with a reward of $r_{i_{1}} = \frac{1- \gamma}{|\mathcal{C}_{i}|}$, where $0 < \gamma < 1$; a fictive location $i_{2}$ with a reward of $r_{i_{2}} = \epsilon$, where $\epsilon = \frac{M + 1}{2}$ and $M = \max_{\mathcal{C}_{i} \in \mathcal{C}} \mathset{\frac{|\mathcal{C}_{i}| - 1}{|\mathcal{C}_{i}|} (1 - \gamma)}$; and a set-customer $j$ with spawning demands $d^{1}_{j} = 1, d^{t}_{j} = 0, \forall t \in \mathcal{T} : t \geq 2$.
    We build preference rules $a_{ij}$ so that each authentic location $i_{1}$ captures all element-customers that belong to set $\mathcal{C}_{i}$, and each fictive location $i_{2}$ captures solely the respective set-customer.

    We remark that the reward for a time period $t$ is always between $\epsilon$ (\ie, if the provider chooses a fictive location that captures the demand unit from the respective set-customer) and $(1 - \gamma)$ (\ie, if the provider chooses an authentic location that captures the demand unit from all element-customers in the respective set).
    Note that, in a given time period $t$, the provider prefers a fictive location over an authentic one if the latter cannot capture the demand unit from all element-customers in the respective set.
    As a direct implication, if $K$ is the maximum number of disjoint sets in $\mathcal{C}$, the total reward over the planning horizon is at most $K (1 - \gamma) + (n - K) \epsilon$, where $K (1 - \gamma)$ units come from authentic locations describing mutually disjoint sets and $(n - K) \epsilon$ units come from fictive locations.
    We look therefore for an objective value $Z$ of the 1-1-\ourproblem{} instance equal to $K (1 - \gamma) + (n - K) \epsilon$.

    \textbf{Forward Direction.}
    Assume that there is a certificate $\mathcal{D} = \mathset{\mathcal{C}_{h(1)}, \ldots \mathcal{C}_{h(K)}}$ to the SPP decision question, \ie, a subcollection $\mathcal{D}$ with at least $K$ mutually disjoint sets, where function $h(t)$ maps the $t$-th set in the certificate to its index in collection $\mathcal{C}$.
    We show that this certificate implies the existence of a certificate to the 1-\ourproblem{} decision question, \ie, a location policy $\boldsymbol{y}$ for the 1-\ourproblem{} instance with an objective value of at least $Z = K (1 - \gamma) + (n - K) \epsilon$.

    In the associated solution of the 1-\ourproblem{}, the provider opens authentic locations linked to sets $\mathcal{C}_{h(t)}$ for time periods $t \in \mathset{1, \ldots,K}$, and fictive locations for the remaining time periods $t \in \mathset{K+1, \ldots, n}$.
    Since sets $\mathcal{C}_{h(1)}, \ldots \mathcal{C}_{h(K)}$ are mutually disjoint, we obtain a reward of $(1 - \gamma)$ from each time period between $1$ and $K$, and a reward of $\epsilon$ from each time period between $K + 1$ and $n$, which yields a total reward of $K (1 - \gamma) + (n - K) \epsilon$.
    Thus, the existence of a certificate $\mathcal{D}$ to the SPP decision question implies the existence of a certificate $\boldsymbol{y}$ for the 1-\ourproblem{} decision question.

    \textbf{Backward Direction.}
    Assume that there is no certificate $\mathcal{D} = \mathset{\mathcal{C}_{h(1)}, \ldots \mathcal{C}_{h(K)}}$ to the SPP decision question 
    with at least $K$ mutually disjoint sets.
    We show that this implies no location policy $\boldsymbol{y}$ for the 1-\ourproblem{} instance with an objective value of at least $Z = K (1 - \gamma) + (n - K) \epsilon$.

    We know from the forward direction that, if the SPP instance has at least $K$ mutually disjoint sets, then there is a location policy $\boldsymbol{y}$ such that the total reward of the associated 1-\ourproblem{} instance is $Z = K (1 - \gamma) + (n - K) \epsilon$.
    Moreover, from the implication at the start of this proof, if the maximum number of mutually disjoint sets for the SPP instance is at most $K^{\prime} < K$, then the optimal value of the associated 1-\ourproblem{} instance is $Z = K^{\prime} (1 - \gamma) + (n - K^{\prime}) \epsilon <  K (1 - \gamma) + (n - K) \epsilon$. Hence, if there are no $K$ mutually disjoint sets in collection $\mathcal{C}$, then there cannot be a location policy $\boldsymbol{y}$ for the associated  1-\ourproblem{} instance with an objective value of at least $Z = K (1 - \gamma) + (n - K) \epsilon$.
    
    Since both directions hold, the 1-\ourproblem{} is in fact NP-hard.
    We now show that the 1-\ourproblem{} is as inapproximable as the SPP through the so-called gap technique \citep{schuurman2001approximation}.
    Let $obj_{dSPP} (\mathcal{C}^{ANSWER})$ be the optimal value assigned to an instance $\mathcal{C}^{ANSWER}$ of the decision version of the SPP, and $obj_{CCD} (f(\mathcal{C}^{ANSWER}))$ be the optimal value assigned to the associated 1-\ourproblem{} instance $f(\mathcal{C}^{ANSWER})$, which can be determined in polynomial time through the reduction presented earlier.
    Let  $\mathcal{C}^{YES}$ and $\mathcal{C}^{NO}$ denote instances of the decision version of the SPP with the same $\mathcal{C}$ but with different target outcomes: the former resulting in $YES$ and the latter resulting in $NO$. 
    Since the SPP cannot be approximated within a factor $|\mathcal{C}|^{1-\alpha}$ for any $\alpha  > 0$ \citep{ausiello1980structure, hastadCliqueHardApproximate1996, hazanComplexityApproximatingKset2006}, it holds that $\frac{obj_{dSPP} (\mathcal{C}^{NO})}{obj_{dSPP} (\mathcal{C}^{YES})}  < |\mathcal{C}|^{1-\alpha}$ for any $\alpha > 0$.
    We now analyze the ratio $\frac{obj_{CCD} (f(\mathcal{C}^{NO}))}{obj_{CCD} (f(\mathcal{C}^{YES}))}$ for the 1-\ourproblem{} to determine if we can answer the decision version of the SPP through the 1-\ourproblem{}, where $obj_{CCD} (f(\mathcal{C})) = obj_{dSPP} (\mathcal{C}) (1 -\gamma - \epsilon) + n \epsilon$ holds by construction: 
    \begin{align*}
        \frac{obj_{CCD} (f(\mathcal{C}^{NO}))}{obj_{CCD} (f(\mathcal{C}^{YES}))} =
        \frac{obj_{dSPP} (\mathcal{C}^{NO}) (1 - \gamma - \epsilon) + n \epsilon}{obj_{dSPP} (\mathcal{C}^{YES}) (1 - \gamma - \epsilon) + n \epsilon} < \frac{|\mathcal{C}|^{1-\alpha} obj_{dSPP} (\mathcal{C}^{YES}) (1 - \gamma - \epsilon) + n \epsilon}{obj_{dSPP} (\mathcal{C}^{YES}) (1 - \gamma - \epsilon) + n \epsilon} .
    \end{align*}
    
    We can then simplify the right-hand side by employing the fact that, if $a > b$, then $g(x) = \frac{a + x}{b + x} < \frac{a}{b}$:
    \begin{align*}
        \frac{|\mathcal{C}|^{1-\alpha} obj_{dSPP} (\mathcal{C}^{YES}) (1 - \gamma - \epsilon) + n \epsilon}{obj_{dSPP} (\mathcal{C}^{YES}) (1 - \gamma - \epsilon) + n \epsilon} < \frac{|\mathcal{C}|^{1-\alpha} obj_{dSPP} (\mathcal{C}^{YES}) (1 - \gamma - \epsilon)}{obj_{dSPP} (\mathcal{C}^{YES}) (1 - \gamma - \epsilon)} < |\mathcal{C}|^{1 - \alpha} = T^{1 - \alpha}
    \end{align*}

    Therefore, the 1-\ourproblem{} also cannot be approximated within a factor $T^{1 - \alpha}$, unless P = NP. \Halmos

\subsection{Proof of Theorem~4}
\label{apx:nphard-identical}

We define the 3SAT, and then present the proof.

\vspace{0.3cm}
\boxxx{
{\bf 3SAT}: \\
{\sc instance}: A finite set of $n$ Boolean variables $\mathcal{B} = \mathset{B_1, \ldots, B_n}$ and a finite set of $m$ clauses $\mathcal{C} = \mathset{\mathcal{C}_1, \ldots, \mathcal{C}_m}$, each with exactly three variables.\\
{\sc question}: Is there a literal assignment such that all clauses are satisfied?
}
\vspace{0.3cm}

\textbf{Proof.}
The 1-\ourproblem{} with identical rewards is in NP as the general \ourproblem{} is in NP (see Theorem~2).
    We show here that the 1-\ourproblem{} with identical rewards is NP-hard by reducing the 3SAT to it. 
    By showing that if there is a certificate satisfying the 3SAT, then there is a certificate satisfying the 1-\ourproblem{} decision question (referred to as forward direction) and vice-versa (referred to as backward direction), it holds that the 1-\ourproblem{} with identical rewards is NP-hard.

    \textbf{Reduction.} By hypothesis, $r_{i}  = R, \forall i \in \mathcal{I}$ and $R \in \mathbb{Q^+}$. 
    Now, consider the 1-\ourproblem{} instance with identical rewards built from a 3SAT instance with $2n$ locations, one for each assignment of a literal to a variable (\eg, $[x_1, true]$ and $[x_1, false]$ are two different locations); $n + m$ customers, one for each variable and each clause; and $n$ time periods.
    We set $d^{1}_{j} = 1, d^{t}_{j} = 0, \forall t \in \mathcal{T} : t \geq 2$ and build preference rules $a_{ij}$ so that each location (\ie, assignment of literal to variable) captures satisfied clause-customers (\ie, customers originated from clauses) and the respective variable-customer (\ie, customers originated from the variables).
    Note that the total reward for the planning horizon is at most $R (n + m)$, which can only be achieved in $n$ time periods if each variable-customer is captured at least once (\ie, each variable has been assigned to a literal), as well as each clause-customer is captured once (\ie, each clause has been satisified through some assignment).
    In fact, we look for an objective value $Z$ of the 1-\ourproblem{} with identical rewards equal to $R (n + m)$, which guarantees that there is a literal assignment that satisfies all clauses.

    \textbf{Forward Direction.}
    Assume that there is a certificate to the 3SAT, \ie, an assignment of literals to variables in $\mathcal{B}$ such that the clauses in $\mathcal{C}$ are satisfied.
    We show that this certificate implies the existence of a certificate to the 1-\ourproblem{} decision question, \ie, a location policy $\boldsymbol{y}$ for the 1-\ourproblem{} instance with an objective value of at least $Z = R (n + m)$.

    In the associated solution of the 1-\ourproblem{}, the provider opens, at time period $t$, the location linked to the $t$-th variable and its literal.
    Since the assignment satisfies all clauses, each clause-customer and each variable-customer is captured at least once over the planning horizon, yielding a total reward of $R (n + m)$.
    Thus, the existence of a certificate to the 3SAT implies the existence of a certificate $\boldsymbol{y}$ for the 1-\ourproblem{} decision question.

    \textbf{Backward Direction.}
    Assume that there is a certificate $\boldsymbol{y}$ to the 1-\ourproblem{} decision question, \ie, a location policy $\boldsymbol{y}$ with an objective value of at least $Z = R (n + m)$.
    We show that this certificate implies the existence of a certificate to the 3SAT, \ie, an assignment of literals to variables in $\mathcal{B}$ such that the clauses in $\mathcal{C}$ are satisfied.

    The backward assumption guarantees that the total reward is at least $Z = R (n + m)$.
    In turn, the only way to obtain such a solution is by capturing each clause-customer and each variable-customer at least once throughout the planning horizon, which means that each variable is assigned to a literal while satisfying all clauses.
    Thus, the existence of a certificate $\boldsymbol{y}$ to the 1-\ourproblem{} decision question implies the existence of a certificate of the 3SAT.

    Since both directions hold, the 1-\ourproblem{} with identical rewards is in fact NP-hard. 
    In addition, since the 3SAT is strongly NP-hard \citep{garey1979computers}, the same result holds for the 1-\ourproblem{}. \Halmos

\subsection{Proof of Theorem~5}
\label{apx:approximation}

For the sake of simplicity, we abuse the notation by writing $\mathcal{J} (i^{B}_{1}, \ldots, i^{B}_{T})$ instead of $\mathcal{J} (\{i^{B}_{1}, \ldots, i^{B}_{T}\})$ for the set of customers captured by locations $i^{B}_{1}, \ldots, i^{B}_{T}$.
Based on Algorithm~1, we define the location $i^{B}_{t}$ chosen by the Backward Greedy Heuristic at time period $t$ as the one with the largest marginal reward:
\begin{align}
    i^{B}_{t} = \argmax_{i \in \mathcal{I}} \mathset{ \sum_{\substack{s \in \mathcal{T}: \\ s \leq t}} \sum_{\substack{j \in \mathcal{J}(i): \\ j \not\in \mathcal{J} (i^{B}_{t+1}, \ldots, i^{B}_{T})}} d^{s}_{j} } \label{eq:bcw-selection}
\end{align}
and the total reward of the heuristic location policy $i^{B}_{1}, \ldots, i^{B}_{T}$ as follows:
\begin{align*}
    z(\boldsymbol{y}^{B}) = 
    R \left[
    \sum_{j \in \mathcal{J} (i^{B}_{T})} \sum_{\substack{s \in \mathcal{T}: \\ s \leq T}} d^{s}_{j} + 
    \sum_{\substack{j \in \mathcal{J} (i^{B}_{T-1}): \\ j \not\in \mathcal{J} (i^{B}_{T}) }} \sum_{\substack{s \in \mathcal{T}: \\ s \leq T - 1}} d^{s}_{j} + \ldots + 
    \sum_{\substack{j \in \mathcal{J} (i^{B}_{1}): \\ j \not\in \mathcal{J} (i^{B}_{2}, \ldots, i^{B}_{T}) }} \sum_{\substack{s \in \mathcal{T}: \\ s \leq 1}} d^{s}_{j}
    \right]
\end{align*}

\textbf{Proof.}
    We analyze the difference $z(\boldsymbol{y}^{\star}) - z(\boldsymbol{y}^{B})$ to upper bound it by $z(\boldsymbol{y}^{B})$.
    We rewrite first the total reward as a function of spawning demands $d^{s}_{j}$, grouping them by time period $s$, and the identical reward $R$:
    \begin{align*}
        z(\boldsymbol{y}^{\star}) - z(\boldsymbol{y}^{B}) 
        =  R \left[
            \sum_{j \in \mathcal{J}(i^{\star}_{1}, i^{\star}_{2}, \ldots, i^{\star}_{T})} d^{1}_{j} + 
            \sum_{j \in \mathcal{J}(i^{\star}_{2}, \ldots, i^{\star}_{T})} d^{2}_{j} + \ldots +
            \sum_{j \in \mathcal{J}(i^{\star}_{T-1}, i^{\star}_{T})} d^{T-1}_{j} +
            \sum_{j \in \mathcal{J}(i^{\star}_{T})} d^{T}_{j}
        \right] 
        - \\ R \left[
            \sum_{j \in \mathcal{J}(i^{B}_{1}, i^{B}_{2}, \ldots, i^{B}_{T})} d^{1}_{j} + 
            \sum_{j \in \mathcal{J}(i^{B}_{2}, \ldots, i^{B}_{T})} d^{2}_{j} + \ldots +
            \sum_{j \in \mathcal{J}(i^{B}_{T-1}, i^{B}_{T})} d^{T-1}_{j} +
            \sum_{j \in \mathcal{J}(i^{B}_{T})} d^{T}_{j}
        \right]  = (\#^{1}).
    \end{align*}
    
    We can further upper bound this difference by taking spawning demands that appear with a positive coefficient after executing the subtraction by time periods as follows:
    \begin{align*}
        (\#^{1}) \leq 
        R \left[
            \sum_{\substack{j \in \mathcal{J}(i^{\star}_{1}, i^{\star}_{2}, \ldots, i^{\star}_{T}): \\ \substack{j \not\in \mathcal{J}(i^{B}_{1}, i^{B}_{2}, \ldots, i^{B}_{T})}}} d^{1}_{j} + 
            \sum_{\substack{j \in \mathcal{J}(i^{\star}_{2}, \ldots, i^{\star}_{T}): \\ j \not\in \mathcal{J}(i^{B}_{2}, \ldots, i^{B}_{T})}} d^{2}_{j} + \ldots +
            \sum_{\substack{j \in \mathcal{J}(i^{\star}_{T-1}, i^{\star}_{T}): \\ j \not\in \mathcal{J}(i^{B}_{T-1}, i^{B}_{T})}} d^{T-1}_{j} +
            \sum_{\substack{j \in \mathcal{J}(i^{\star}_{T}): \\ j \not\in \mathcal{J}(i^{B}_{T})}} d^{T}_{j}
        \right] = (\#^{2}).
    \end{align*}
    
    We now isolate spawning demands related to location $i^{\star}_{t}$ for each time period $t$, and show that they are upper bounded by the contribution of location $i^{B}_{t}$ to the heuristic location policy $i^{B}_{1}, \ldots, i^{B}_{T}$.
    We start with time period $T$, where $\sum_{\substack{s \in \mathcal{T}: \\ s \leq T}}\sum_{\substack{j \in \mathcal{J}(i^{\star}_{T}): \\ j \not\in \mathcal{J}(i^{B}_{T})}} d^{s}_{j} \leq \sum_{\substack{s \in \mathcal{T}: \\ s \leq T}}\sum_{j \in \mathcal{J}(i^{\star}_{T})} d^{s}_{j}$ holds trivially and $ \sum_{\substack{s \in \mathcal{T}: \\ s \leq T}}\sum_{j \in \mathcal{J}(i^{\star}_{T})} d^{s}_{j} \leq \sum_{\substack{s \in \mathcal{T}: \\ s \leq T}} \sum_{j \in \mathcal{J}(i^{B}_{T})} d^{s}_{j}$ holds thanks to Equation~\eqref{eq:bcw-selection}:
    \begin{align*}
        (\#^{2}) = R \left[
            \sum_{\substack{j \in \mathcal{J}(i^{\star}_{1}, i^{\star}_{2}, \ldots, i^{\star}_{T-1}): \\ \substack{j \not\in \mathcal{J}(i^{B}_{1}, i^{B}_{2}, \ldots, i^{B}_{T})}}} d^{1}_{j} + 
            \sum_{\substack{j \in \mathcal{J}(i^{\star}_{2}, \ldots, i^{\star}_{T-1}): \\ j \not\in \mathcal{J}(i^{B}_{2}, \ldots, i^{B}_{T})}} d^{2}_{j} + \ldots +
            \sum_{\substack{j \in \mathcal{J}(i^{\star}_{T-1}): \\ j \not\in \mathcal{J}(i^{B}_{T-1}, i^{B}_{T})}} d^{T-1}_{j} +
            \sum_{\substack{s \in \mathcal{T}: \\ s \leq T}}\sum_{\substack{j \in \mathcal{J}(i^{\star}_{T}): \\ j \not\in \mathcal{J}(i^{B}_{T})}} d^{s}_{j} 
        \right] \leq \\
        R \left[
            \sum_{\substack{j \in \mathcal{J}(i^{\star}_{1}, i^{\star}_{2}, \ldots, i^{\star}_{T-1}): \\ \substack{j \not\in \mathcal{J}(i^{B}_{1}, i^{B}_{2}, \ldots, i^{B}_{T})}}} d^{1}_{j} + 
            \sum_{\substack{j \in \mathcal{J}(i^{\star}_{2}, \ldots, i^{\star}_{T-1}): \\ j \not\in \mathcal{J}(i^{B}_{2}, \ldots, i^{B}_{T})}} d^{2}_{j} + \ldots +
            \sum_{\substack{j \in \mathcal{J}(i^{\star}_{T-1}): \\ j \not\in \mathcal{J}(i^{B}_{T-1}, i^{B}_{T})}} d^{T-1}_{j}
        \right] +
        R \left[
            \sum_{\substack{s \in \mathcal{T}: \\ s \leq T}} \sum_{j \in \mathcal{J}(i^{B}_{T})} d^{s}_{j} 
        \right] = (\#^{3}).
    \end{align*}

    We do the same for time period $T-1$, where again $\sum_{\substack{s \in \mathcal{T}: \\ s \leq T-1}} \sum_{\substack{j \in \mathcal{J}(i^{\star}_{T-1}): \\ j \not\in \mathcal{J}(i^{B}_{T-1}, i^{B}_{T})}} d^{s}_{j} \leq \sum_{\substack{s \in \mathcal{T}: \\ s \leq T-1}} \sum_{\substack{j \in \mathcal{J}(i^{\star}_{T-1}): \\ j \not\in \mathcal{J}(i^{B}_{T})}} d^{s}_{j}$ holds trivially and $ \sum_{\substack{s \in \mathcal{T}: \\ s \leq T-1}} \sum_{\substack{j \in \mathcal{J}(i^{\star}_{T-1}): \\ j \not\in \mathcal{J}(i^{B}_{T})}} d^{s}_{j}  \leq \sum_{\substack{s \in \mathcal{T}: \\ s \leq T-1}} \sum_{\substack{j \in \mathcal{J}(i^{B}_{T-1}): \\ j \not\in \mathcal{J} (i^{B}_{T})}} d^{s}_{j}$ holds thanks to Equation~\eqref{eq:bcw-selection}:
    {\tiny
    \begin{align*}
        (\#^{3}) = R \left[
            \sum_{\substack{j \in \mathcal{J}(i^{\star}_{1}, i^{\star}_{2}, \ldots, i^{\star}_{T-2}): \\ \substack{j \not\in \mathcal{J}(i^{B}_{1}, i^{B}_{2}, \ldots, i^{B}_{T})}}} d^{1}_{j} + 
            \sum_{\substack{j \in \mathcal{J}(i^{\star}_{2}, \ldots, i^{\star}_{T-2}): \\ j \not\in \mathcal{J}(i^{B}_{2}, \ldots, i^{B}_{T})}} d^{2}_{j} + \ldots +
            \sum_{\substack{s \in \mathcal{T}: \\ s \leq T-1}} \sum_{\substack{j \in \mathcal{J}(i^{\star}_{T-1}): \\ j \not\in \mathcal{J}(i^{B}_{T-1}, i^{B}_{T})}} d^{s}_{j}
        \right] +
        R \left[
            \sum_{\substack{s \in \mathcal{T}: \\ s \leq T}} \sum_{j \in \mathcal{J}(i^{B}_{T})} d^{s}_{j} 
        \right] \leq \\
        R \left[
            \sum_{\substack{j \in \mathcal{J}(i^{\star}_{1}, i^{\star}_{2}, \ldots, i^{\star}_{T-2}): \\ \substack{j \not\in \mathcal{J}(i^{B}_{1}, i^{B}_{2}, \ldots, i^{B}_{T})}}} d^{1}_{j} + 
            \sum_{\substack{j \in \mathcal{J}(i^{\star}_{2}, \ldots, i^{\star}_{T-2}): \\ j \not\in \mathcal{J}(i^{B}_{2}, \ldots, i^{B}_{T})}} d^{2}_{j} + \ldots
        \right] +
        R \left[
            \sum_{\substack{s \in \mathcal{T}: \\ s \leq T-1}} \sum_{\substack{j \in \mathcal{J}(i^{B}_{T-1}): \\ j \not\in \mathcal{J} (i^{B}_{T})}} d^{s}_{j} +
            \sum_{\substack{s \in \mathcal{T}: \\ s \leq T}} \sum_{j \in \mathcal{J}(i^{B}_{T})} d^{s}_{j} 
        \right]. \\ 
    \end{align*}
    }%
    
    Intuitively, repeating this reasoning gradually builds the total reward $z(\boldsymbol{y}^{B})$ of the heuristic location policy $i^{B}_{1}, \ldots, i^{B}_{T}$ on the right-hand side.
    Therefore, it holds that $z(\boldsymbol{y}^{\star}) - z(\boldsymbol{y}^{B}) \leq z(\boldsymbol{y}^{B})$, and the backward greedy heuristic is a 2-approximation algorithm for the 1-\ourproblem{} with identical rewards. \Halmos

\subsection{Proof of Theorem~6}
\label{apx:polynomial}

\begin{lemma}
    \label{lem:no-repetition}
    Let $i_{1}, \ldots, i_{T}$ be a feasible solution of a 1-\ourproblem{} instance with loyal customers.
    If there is a location $i^{\prime} \in \mathcal{I}$ chosen for two or more time periods of the planning horizon (\ie, $\exists \, \mathcal{T}^{\prime} = \{t_{1}^{\prime}, t_{2}^{\prime}, \ldots, t_{K-1}^{\prime}, t_{K}^{\prime}\} \subseteq \mathcal{T}$ such that $i_{t} = i^{\prime} \, \forall t \in \mathcal{T}^{\prime}$), there is an equivalent solution (\ie, with the same total reward) $\tilde{i}_{1}, \ldots, \tilde{i}_{T}$ without repetition where $\tilde{i}_{t} = \emptyset \, \forall t \in \mathcal{T}^{\prime}$, $\tilde{i}_{t} = i_{t} \, \forall t \in \mathcal{T} \backslash \mathcal{T}^{\prime}$ and $\tilde{i}_{t_{K}^{\prime}} = i^{\prime}$.
\end{lemma}

\textbf{Proof.}
    Let $i^{\prime}$ be the repeated location.
    Under loyal customers, we can split the total reward into contributions from location $i^{\prime}$ and those from other locations, denoted by $Z$, because there is no intersection between customers captured by location $i^{\prime}$ and those captured by other locations.
    This reasoning gives
    \begin{align*}
        z(i_{1}, \ldots, i_{T}) = Z + \sum_{t \in \mathcal{T}^{\prime}} \rho(i^{\prime}, t \mid i_{1}, \ldots, i_{T}) =  
        Z + r_{i^{\prime}} \sum_{j \in \mathcal{J}(i^{\prime})} \left( \sum_{\substack{s \in \mathcal{T}: \\ s > 0 \\ s \leq t_{1}^{\prime}}}  d^{s}_{j} + \sum_{\substack{s \in \mathcal{T}: \\ s > t_{1}^{\prime} \\ s \leq t_{2}^{\prime}}} d^{s}_{j} + \ldots + \sum_{\substack{s \in \mathcal{T}: \\ s > t_{K-1}^{\prime} \\ s \leq t_{K}^{\prime}}} d^{s}_{j}  \right) = \\
        Z + r_{i^{\prime}} \sum_{j \in \mathcal{J}(i^{\prime})} \sum_{\substack{s \in \mathcal{T}: \\ s \leq t_{K}^{\prime}}} d^{s}_{j} =
        z(\tilde{i}_{1}, \ldots, \tilde{i}_{T}). \hfill \square 
    \end{align*}

\begin{corollary}
    Let $i_{1}, \ldots, i_{T}$ be a feasible solution of a 1-\ourproblem{} instance with loyal customers. 
    We assume that this feasible solution has no repetition without loss of generality.
    The marginal reward function can be rewritten as $\rho(i_{t}, t) = r_{i_{t}} \sum_{j \in \mathcal{J}(i_{t})} \sum_{\substack{s \in \mathcal{T}: \\ s \leq t}} d^{s}_{j}$, and the total reward function can be rewritten as $z(i_{1}, \ldots, i_{T}) =
            \sum_{t \in \mathcal{T}} \rho(i_{t}, t)$.
\end{corollary}

\textbf{Proof.}
    This is a direct outcome of Lemma~\ref{lem:no-repetition}. \Halmos

We are now ready to present the proof of Theorem 6.

\textbf{Proof.}
    Each feasible solution  $i_{1}, \ldots, i_{T}$ of the 1-\ourproblem{} with loyal customers has an equivalent solution without repetition, which is nothing but an exact assignment of locations to time periods.
    Thus, we can solve this 1-\ourproblem{} instance by solving an assignment problem in polynomial time \citep{kuhn1955hungarian}.
    Let $E = I-T$ be the difference between the number of candidate locations and time periods.
    We can create an instance of the assignment in polynomial time by setting the weight of assigning candidate location $i$ to time period $t$ as $\rho(i, t)$.
    If $E = 0$, we do not have to conduct further adaptations.
    If $E < 0$, we need to add $|E|$ virtual candidate locations, such that their weight is $0$ to all time periods.
    Similarly, if $E > 0$, we need to add $|E|$ virtual time periods, such that their weight is $0$ to all candidate locations.
    The feasible solution can be naturally drawn from the assignment of (true) candidate locations to (true) time periods in polynomial time.
    Thus, the 1-\ourproblem{} with loyal customers is polynomially solvable. \Halmos

\section{Analytical Procedure}
\label{apx:analytical}

We present here the proof of Theorem~7.
We first write the primal subproblem $w^{P}_{j} (\boldsymbol{y})$ as follows for the 1-\ourproblem{}, where Constraints~(1c)--(1d) can be merged into one and Constraints~(1e) can be omitted:
\begin{subequations}
    \label{pgm:primal-1-dflp-ccd}
    \begin{align}
        w^{P}_{j} (\boldsymbol{y}): \quad  
        \max_{\boldsymbol{x}} \quad
        &  \sum_{t \in \mathcal{T}} \sum_{\substack{\ell{} \in \mathcal{T}^{S} : \\ \ell{} < t}} \sum_{i \in \mathcal{I}}G^{\ell{}t}_{ij}  x^{\ell{}t}_{i}
        && \label{eq:primal-1-dflp-ccd-obj}\\
        \text{s.t.} \quad
        & \sum_{\substack{\ell{} \in \mathcal{T}^{S}: \\ \ell{} < t}} x^{\ell{}t}_{i} = a_{ij}y^{t}_{i}
         && \forall i \in \mathcal{I},  \forall t \in \mathcal{T} & [\lambda^{t}_{i}] 
        \label{eq:primal-1-dflp-ccd-ct1} \\
        & \sum_{\substack{s \in \mathcal{T}^{F}: \\ s > t}} \sum_{i \in \mathcal{I}} x^{ts}_{i} - \sum_{\substack{s \in \mathcal{T}^{S}: \\ s < t}} \sum_{i \in \mathcal{I}} x^{st}_{i} = 0
        && \forall i \in \mathcal{I}, \forall t \in \mathcal{T} & [\theta^{t}] 
        \label{eq:primal-1-dflp-ccd-ct2} \\
        & \sum_{s \in \mathcal{T}^{F}} \sum_{i \in \mathcal{I}} x^{0s}_{i} = 1
        && & [\theta^{0}] 
        \label{eq:primal-1-dflp-ccd-ct3} \\
        & x^{\ell{}t}_{i} \in \mathbb{R}^{+}
        && \forall i \in \mathcal{I}, \forall \ell{} \in \mathcal{T}^{S}.  \forall t \in \mathcal{T}^{F} : \ell{} < t, \label{eq:primal-1-dflp-ccd-dm1}
    \end{align}
\end{subequations}
In this sense, the dual subproblem $w^{D}_{j} (\boldsymbol{y})$ can be rewritten as follows after some rearrangement:
\begin{subequations}
    \label{pgm:dual2-1-dflp-ccd}
    \begin{align}
        w^{D}_{j} (\boldsymbol{y}): \quad 
        \min_{\boldsymbol{\lambda}, \boldsymbol{\theta}} \quad
        &  \sum_{t \in \mathcal{T}} \sum_{i \in \mathcal{I}} a_{ij} y^{t}_{i} \lambda^{t}_{i} + \theta^{0} 
        && \label{eq:dual2-1-dflp-ccd-obj}\\
        \text{s.t.} \quad
        & \lambda^{t}_{i} \geq \max_{\ell{} \in \mathcal{T}^{S} : \ell{} < t} \mathset{G^{\ell{}t}_{ij} - \theta^{\ell{}} + \theta^{t}}
        && \forall i \in \mathcal{I}, \forall t \in \mathcal{T} \label{eq:dual2-1-dflp-ccd-ct1} \\
        & \theta^{\ell{}} \geq G^{\ell{}(T+1)}_{ij}
        && \forall i \in \mathcal{I}, \forall \ell{} \in \mathcal{T}^{S} \label{eq:dual2-1-dflp-ccd-ct2} \\
        & \lambda^{t}_{i} \in \mathbb{R}
        && \forall i \in \mathcal{I}, \forall t \in \mathcal{T} \label{eq:dual2-1-dflp-ccd-dm1} \\
        & \theta^{t} \in \mathbb{R}
        && \forall t \in \mathcal{T}^{S}. \label{eq:dual2-1-dflp-ccd-dm2}
    \end{align}
\end{subequations}
where dual variables $\lambda^{t}_{i}$ and $\theta^{t}$ are related to Constraints \eqref{eq:primal-1-dflp-ccd-ct1} and \eqref{eq:primal-1-dflp-ccd-ct2}--\eqref{eq:primal-1-dflp-ccd-ct3}, respectively.

\begin{lemma}
    \label{lem:projection}
    In an optimal solution of the dual subproblem $w^{D}_{j} (\boldsymbol{y})$, Constraints~\eqref{eq:dual2-1-dflp-ccd-ct1} may be posed with equality without changing the optimal objective value.
\end{lemma}

\textbf{Proof.}
    Let $\boldsymbol{x}^{\star}$ be the optimal solution of the primal subproblem $w^{P}_{j}(\boldsymbol{y})$ computed in polynomial time based on Proposition~1.
    If $y^{t}_{i} = 0$, variable $\lambda^{t}_{i}$ does not appear in Objective Function \eqref{eq:dual2-1-dflp-ccd-obj}, and we can satisfy Constraints \eqref{eq:dual2-1-dflp-ccd-ct1} for $i$ and $t$ with equality without changing the objective value for the term $i$ and $t$ of the sum.
    If $y^{t}_{i} = 1$, there are two cases.
    On the one hand, if $a_{ij} = 0$, variable $\lambda^{t}_{i}$ also does not appear in Objective Function \eqref{eq:dual2-1-dflp-ccd-obj}, and the previous reasoning applies.
    On the other hand, if $a_{ij} = 1$, we know that ${x^{st}_{i}}^{\star} = 1$ for some time period $s < t$ through Constraints~\eqref{eq:primal-1-dflp-ccd-ct2}. Then, by complementary slackness, Constraints~\eqref{eq:dual2-1-dflp-ccd-ct1} must be satisfied for $i$ and $t$ with equality (\ie, if ${x^{st}_{i}}^{\star} = 1$, then $\lambda^{t}_{i} = G^{st}_{ij} - \theta^{s} + \theta^{t}$).
    Therefore, Constraints~\eqref{eq:dual2-1-dflp-ccd-ct1} may be posed with equality without changing the optimal objective value.  \Halmos

Lemma~\ref{lem:projection} allow us to project out variables $\lambda^{t}_{i}$ out of the dual subproblem $w^{D}_{j}(\boldsymbol{y})$, and obtain an equivalent form for Constraints \eqref{eq:dual2-1-dflp-ccd-ct1}.
In other words, from the proof of Lemma~\ref{lem:projection}, we have:
\begin{align*}
    & \lambda^{t}_{i} \geq \max_{\ell{} \in \mathcal{T}^{S} : \ell{} < t} \mathset{G^{\ell{}t}_{ij} - \theta^{\ell{}} + \theta^{t}}
        && \forall i \in \mathcal{I}, \forall t \in \mathcal{T} \\
   \implies & G^{st}_{ij} - \theta^{s} + \theta^{t} \geq \max_{\ell{} \in \mathcal{T}^{S} : \ell{} < t, \ell \neq s} \mathset{G^{\ell{}t}_{ij} - \theta^{\ell{}} + \theta^{t}}
        && \substack{\forall i \in \mathcal{I}, \forall t \in \mathcal{T}, \forall s \in \mathcal{T}^{S} :\\ s < t, {x^{st}_{i}}^{\star} = 1} \\
    \implies & G^{st}_{ij} - \theta^{s} + \theta^{t} \geq G^{\ell{}t}_{ij} - \theta^{\ell{}} + \theta^{t} 
    && \substack{\forall i \in \mathcal{I},  \forall t \in \mathcal{T}, \forall s \in \mathcal{T}^{S}, \forall \ell{} \in \mathcal{T}^{S}: \\ s < t, {x^{st}_{i}}^{\star} = 1, \ell{} < t, \ell{} \neq s}
\end{align*}
which gives the following reduced but equivalent dual subproblem $w^{R}_{j} (\boldsymbol{y})$:
\begin{subequations}
    \label{pgm:dual3-dsflp-ccd}
    \begin{align}
        w^{R}_{j} (\boldsymbol{y}): \quad 
        \min_{\boldsymbol{\theta}} \quad
        &  \theta^{0}
        && \label{eq:dual3-dsflp-ccd-obj}\\
        \text{s.t.} \quad
        & \theta^{\ell{}} \geq G^{\ell{}t}_{ij} - G^{st}_{ij} + \theta^{s} 
        && \substack{\forall i \in \mathcal{I},  \forall t \in \mathcal{T}, \forall s \in \mathcal{T}^{S}, \forall \ell{} \in \mathcal{T}^{S}: \\ s < t, {x^{st}_{i}}^{\star} = 1, \ell{} < t, \ell{} \neq s} \label{eq:dual3-dsflp-ccd-ct1}\\
        & \theta^{\ell{}} \geq G^{\ell{}(T+1)}_{ij}
        && \forall i \in \mathcal{I}, \forall \ell{} \in \mathcal{T}^{S} \label{eq:dual3-dsflp-ccd-ct2} \\
        & \theta^{\ell{}} \in \mathbb{R}
        && \forall \ell{} \in \mathcal{T}^{S} \label{eq:dual3-dsflp-ccd-dm1}
    \end{align}
\end{subequations}
from which variables $\lambda^{t}_{i}$ can be computed based on variables $\theta^{t}$ as follows:
\begin{align}
    \lambda^{t}_{i} = \max_{\ell{} \in \mathcal{T}^{S} : \ell{} < t} \mathset{G^{\ell{}t}_{ij} - \theta^{\ell{}} + \theta^{t}}
        && \forall i \in \mathcal{I}, \forall t \in \mathcal{T}. \label{eq:magical-formula-p}
\end{align}

We can now easily derive an analytical formula to compute a feasible solution for this reduced dual subproblem $w^{R} (\boldsymbol{y})$ based on Constraints \eqref{eq:dual3-dsflp-ccd-ct1}--\eqref{eq:dual3-dsflp-ccd-ct2}:
\begin{align}
    \theta^{\ell{}} = \max_{\substack{i \in \mathcal{I}, s \in \mathcal{T}^{S}, t \in \mathcal{T}: \\ s < t, {x^{st}_{i}}^{\star} = 1 \\ \ell{} < t, \ell{} \neq s}} \mathset{G^{\ell{}t}_{ij} - G^{st}_{ij} + \theta^{s}, G^{\ell{}(T+1)}_{ij}} && \forall \ell{} \in \mathcal{T}^{S}. \label{eq:magical-formula-q}
\end{align}

The main challenge here is that Equation~\eqref{eq:magical-formula-q} requires variable $\theta^{s}$ to compute variable $\theta^{\ell{}}$.
In this sense, we must ensure that there is a feasible order to compute them independently (\ie, without recurrence).
Let $\mathcal{T}^{+} (\boldsymbol{y}) = \mathset{\ell{} \in \mathcal{T} \mid \sum_{i\in \mathcal{I}} a_{ij} y^{\ell{}}_{i} = 1} \cup \mathset{0}$ and $\mathcal{T}^{-} (\boldsymbol{y}) = \mathcal{T} \backslash \mathcal{T}^{+} (\boldsymbol{y})$ be the set of time periods where customer $j$ has been captured ($+$) or remained free ($-$), respectively, in location policy $\boldsymbol{y}$.

\begin{lemma}
    \label{lem:feasibility}
    Equation~\eqref{eq:magical-formula-q} provides feasible values for variables $\theta^{\ell{}}$ in the dual subproblem $w^{R} (\boldsymbol{y})$ if we first compute them for time periods $\ell{} \in \mathcal{T}^{+} (\boldsymbol{y})$ in decreasing order, then for time periods $\ell{} \in \mathcal{T}^{-} (\boldsymbol{y})$.
\end{lemma}

\textbf{Proof.}
    We show first, by contradiction, that there is no recurrent definition in the first pass, where we compute variable $\theta^{\ell{}}$ for time periods $\ell{} \in \mathcal{T}^{+} (\boldsymbol{y})$ in decreasing order.
    Assume, for the sake of contradiction, that Equation~\eqref{eq:magical-formula-q} needs the value of some variable $\theta^{s}$ to compute some variable $\ell{} \in \mathcal{T}^{+} (\boldsymbol{y})$ that has not been computed yet.
    This implies that \textit{(i)} ${x^{\ell{}t}_{{i_{1}}}}^{\star} = 1$ for some time period $t$ and some location ${i_{1}}$ since $\ell{} \in \mathcal{T}^{+} (\boldsymbol{y})$, and \textit{(ii)}  ${x^{st}_{{i_{2}}}}^{\star} = 1$ for some location $i_{2}$ since $\theta^{s}$ appears in Equation~\eqref{eq:magical-formula-q}, which is absurd.
    In other words, we cannot have ${x^{\ell{}t}_{{i_{1}}}}^{\star} = 1$ and ${x^{st}_{{i_{2}}}}^{\star} = 1$  without violating flow conservation of the primal subproblem $w^{P} (\boldsymbol{y})$.
    Therefore, the lemma is correct for $\ell{} \in \mathcal{T}^+(y)$ (\ie, we can determine feasible values for $q^\ell{}$ for time periods $\ell{} \in \mathcal{T}^+(y)$ by applying Equation~\eqref{eq:magical-formula-q} in decreasing order).
    We consider now the second pass, where we compute variable $\theta^{\ell{}}$ for time periods $\ell{} \in \mathcal{T}^{-} (\boldsymbol{y})$.
    Note that we only employ variables $\theta^{s}$ with $s \in \mathcal{T}^{+} (\boldsymbol{y})$ in Equation~\eqref{eq:magical-formula-q}, which were already computed in the first pass.
    Therefore, the lemma holds. \Halmos

\begin{lemma}
    \label{lem:optimality}
    Let $(\boldsymbol{\lambda}, \boldsymbol{\theta})$ be a feasible solution of the dual subproblem $w^{D}_{j} (\boldsymbol{y})$ computed with Equations~\eqref{eq:magical-formula-q} and \eqref{eq:magical-formula-p}.
    This solution is, in fact, optimal.
\end{lemma}

\textbf{Proof.}
    We show the optimality of this feasible solution by strong duality (\ie, by showing that the objective value of the dual solution $(\boldsymbol{\lambda}, \boldsymbol{\theta})$ is equal to the objective value $\sum_{t \in \mathcal{T}} \sum_{\mathcal{\substack{\ell{} \in \mathcal{T} : \\ \ell{} < t}}} G^{\ell{}t}_{ij} {x^{\ell{}t}_{i}}^{\star}$ of the optimal primal solution $\boldsymbol{x}^{\star}$).
    Lemma~\ref{lem:projection} shows that only variables $\lambda^{t}_{i}$ such that $y^{t}_{i} = 1$ and $a_{ij} = 1$ appear in Objective Function \eqref{eq:dual2-1-dflp-ccd-obj}.
    In turn, each variable $\lambda^{t}_{i}$ that appears in the objective function is equal to $G^{st}_{ij} - \theta^{s} + \theta^{t}$ for some ${x^{st}_{i}}^{\star} = 1$ due to complementary slackness.
    Plugging this information into Objective Function \eqref{eq:dual2-1-dflp-ccd-obj} gives $\sum_{t \in \mathcal{T}} \sum_{\mathcal{\substack{\ell{} \in \mathcal{T} : \\ \ell{} < t}}} G^{\ell{}t}_{ij} {x^{\ell{}t}_{i}}^{\star} + \theta^{last(\boldsymbol{y})}$, where $last(\boldsymbol{y})$ is the last time period where customer $j$ was captured in location policy $\boldsymbol{y}$.
    More precisely, variables $\theta^{s}$ and $\theta^{t}$ in the definition of each variable $\lambda^{t}_{i}$ cancel each other through the sum, with the exception of $\theta^{last(\boldsymbol{y})}$.
    Note that, by construction, Equations~\eqref{eq:magical-formula-q} always sets $\theta^{last(\boldsymbol{y})} = 0$.
    Therefore, by strong duality, this feasible solution $(\boldsymbol{\lambda}, \boldsymbol{\theta})$ is, in fact, optimal.    \Halmos

Theorem~7 naturally holds as an outcome of Lemmas~\ref{lem:projection}, \ref{lem:feasibility}, and \ref{lem:optimality}.

\section{Additional Results}
\label{apx:results}

We present here some additional results to support the computational experiments presented in Section~6.2.
In what follows, we filter our benchmark by one attribute (\eg, $|\mathcal{I} |= 50$ locations) to compute average metrics (\eg, solution times) and draw insights on which attributes induce easier or harder instances for each solution method. This strategy may yield large standard deviations because we consider instances with considerably different dimensions at the same time (\eg. instances with $|\mathcal{I} |= 50$ locations have $|\mathcal{J}| \in \mathset{1|\mathcal{I}|, 3|\mathcal{I}|, 5|\mathcal{I}|}$ customers, $|\mathcal{T}| \in \mathset{5,7,9}$ periods, and so on). However, they are still useful to provide some insights on the performance of solution methods for different instance attributes. In addition, Figures~3--6 help us support our claims about the advantages and disadvantages of each solution method without relying on large standard deviations.

\subsection{Computational Performance}

\subsubsection{Comparison Between DIF and SIF.}

Table~\ref{tab:lrz-vs-net-intgap} reports average integrality gaps and solution times, as well as their standard deviations, of SIF and DIF for instances solved to optimality by both formulations.
We compute the integrality gap as $\frac{Z^{\prime} - Z^{\star}}{Z^{\prime}}$, where $Z^{\star}$ is the optimal objective value of some mixed-integer programming formulation and $Z^{\prime}$ is the optimal objective value of its continuous relaxation.
Table~\ref{tab:lrz-vs-net-optgap} reports the number of instances solved to optimality and average optimality gaps (and their standard deviations) of SIF and DIF for instances not solved to optimality by at least one of them.
Recall that we provide optimality gaps proven by the solution  method within the time limit.
\begin{table}[!ht]
    \centering
    \scriptsize
    \caption{Average integrality gaps and solution times of SIF and DIF, as well as their standard deviations, for instances solved to optimality by both formulations.}
    \begin{tabular}{crrrrr}
        \toprule
        \multirow{2}*{Instance attributes} &  \# instances & \multicolumn{2}{c}{SIF} & \multicolumn{2}{c}{DIF} \\ \cmidrule{3-6}
        & considered & int. gap (\%)  & time (min) & int. gap (\%)  & time (min)  \\ \midrule
Benchmark&$277 \, (648)$&$6.13\pm4.00$&$6.84\pm13.13$&$1.90\pm1.97$&$4.48\pm9.49$\\
\midrule
$50$ locations&$142 \, (216)$&$6.63\pm4.29$&$6.73\pm14.77$&$1.45\pm1.71$&$2.48\pm7.28$\\
$100$ locations&$87 \, (216)$&$5.79\pm3.79$&$6.35\pm11.37$&$2.12\pm2.06$&$5.64\pm11.00$\\
$150$ locations&$48 \, (216)$&$5.26\pm3.30$&$8.08\pm10.95$&$2.85\pm2.17$&$8.29\pm10.91$\\
\midrule
$1|\mathcal{I}|$ customers&$149 \, (216)$&$6.21\pm4.23$&$5.04\pm11.44$&$1.11\pm1.57$&$2.73\pm7.80$\\
$3|\mathcal{I}|$ customers&$79 \, (216)$&$6.56\pm4.15$&$9.97\pm15.45$&$2.69\pm2.13$&$5.91\pm9.71$\\
$5|\mathcal{I}|$ customers&$49 \, (216)$&$5.19\pm2.74$&$7.28\pm13.20$&$3.04\pm1.79$&$7.48\pm12.44$\\
\midrule
5 periods&$119 \, (216)$&$5.69\pm3.98$&$6.23\pm13.84$&$2.01\pm1.97$&$3.22\pm7.65$\\
7 periods&$92 \, (216)$&$6.34\pm3.91$&$7.60\pm13.52$&$1.99\pm2.12$&$5.45\pm11.12$\\
9 periods&$66 \, (216)$&$6.63\pm4.13$&$6.90\pm11.27$&$1.57\pm1.74$&$5.39\pm9.90$\\
\midrule
1 facility&$165 \, (216)$&$4.65\pm3.14$&$2.92\pm5.43$&$2.35\pm2.22$&$4.57\pm9.55$\\
3 facilities&$54 \, (216)$&$8.39\pm4.21$&$11.08\pm17.60$&$1.46\pm1.57$&$3.52\pm8.01$\\
5 facilities&$58 \, (216)$&$8.22\pm4.12$&$14.07\pm18.65$&$1.03\pm0.93$&$5.11\pm10.60$\\
\midrule
Identical rewards&$152 \, (324)$&$7.39\pm3.67$&$7.01\pm13.27$&$2.15\pm2.12$&$4.86\pm9.44$\\
Different rewards&$125 \, (324)$&$4.60\pm3.86$&$6.64\pm13.01$&$1.60\pm1.74$&$4.01\pm9.56$\\
\midrule
Short rankings&$163 \, (324)$&$5.22\pm4.02$&$5.14\pm11.08$&$1.27\pm1.33$&$2.17\pm4.60$\\
Long rankings&$114 \, (324)$&$7.43\pm3.61$&$9.28\pm15.33$&$2.80\pm2.37$&$7.78\pm13.07$\\
\midrule
Constant demand&$115 \, (324)$&$5.56\pm4.18$&$7.40\pm13.95$&$1.65\pm1.85$&$4.03\pm8.36$\\
Sparse demand&$162 \, (324)$&$6.54\pm3.83$&$6.45\pm12.55$&$2.08\pm2.05$&$4.80\pm10.23$\\
        \bottomrule
    \end{tabular}
    \label{tab:lrz-vs-net-intgap}
\end{table}
\begin{table}[!ht]
    \centering
    \scriptsize
    \caption{Number of instances solved to optimality, average optimality gaps and their standard deviations, of SIF and DIF, for instances not solved to optimality by at least one of the formulations.}
    \begin{tabular}{crrrrr}
        \toprule
        \multirow{2}*{Instance attributes} & \# instances & \multicolumn{2}{c}{SIF} & \multicolumn{2}{c}{DIF} \\ \cmidrule{3-6}
        & considered & \# opt. & prv. opt. gap (\%) & \# opt. & prv. opt. gap (\%)  \\ \midrule
Benchmark&$371 \, (648)$&$12$&$22.17\pm33.26$&$24$&$14.10\pm25.52$\\
\midrule
$50$ locations&$74 \, (216)$&$3$&$3.36\pm2.21$&$14$&$2.43\pm1.84$\\
$100$ locations&$129 \, (216)$&$3$&$14.67\pm21.56$&$8$&$9.81\pm14.66$\\
$150$ locations&$168 \, (216)$&$6$&$36.22\pm41.16$&$2$&$22.54\pm33.62$\\
\midrule
$1|\mathcal{I}|$ customers&$67 \, (216)$&$2$&$3.55\pm2.40$&$11$&$2.15\pm1.81$\\
$3|\mathcal{I}|$ customers&$137 \, (216)$&$3$&$20.68\pm30.15$&$9$&$11.41\pm19.46$\\
$5|\mathcal{I}|$ customers&$167 \, (216)$&$7$&$30.87\pm38.75$&$4$&$21.10\pm32.04$\\
\midrule
5 periods&$97 \, (216)$&$4$&$17.64\pm27.52$&$7$&$7.40\pm5.62$\\
7 periods&$124 \, (216)$&$0$&$20.96\pm31.62$&$9$&$13.02\pm23.42$\\
9 periods&$150 \, (216)$&$8$&$26.11\pm37.48$&$8$&$19.32\pm32.95$\\
\midrule
1 facility&$51 \, (216)$&$12$&$9.37\pm19.09$&$1$&$15.86\pm28.29$\\
3 facilities&$162 \, (216)$&$0$&$25.68\pm36.06$&$18$&$12.99\pm23.01$\\
5 facilities&$158 \, (216)$&$0$&$22.72\pm33.05$&$5$&$14.67\pm27.09$\\
\midrule
Identical rewards&$172 \, (324)$&$5$&$24.53\pm35.10$&$11$&$20.17\pm32.63$\\
Different rewards&$199 \, (324)$&$7$&$20.13\pm31.54$&$13$&$8.86\pm15.41$\\
\midrule
Short rankings&$161 \, (324)$&$3$&$17.77\pm29.04$&$14$&$8.40\pm13.67$\\
Long rankings&$210 \, (324)$&$9$&$25.55\pm35.87$&$10$&$18.47\pm31.07$\\
\midrule
Constant demand&$209 \, (324)$&$9$&$26.24\pm37.19$&$18$&$13.68\pm24.64$\\
Sparse demand&$162 \, (324)$&$3$&$16.92\pm26.59$&$6$&$14.64\pm26.67$\\
        \bottomrule
    \end{tabular}
    \label{tab:lrz-vs-net-optgap}
\end{table}

\subsubsection{Comparison Between SBD and DIF.}

Table~\ref{tab:bds-d-vs-a-time} reports average solution times, as well as their standard deviations, of DIF and SBD for instances solved to optimality by both methods.
Table~\ref{tab:bds-d-vs-a-optgap} report the number of instances solved to optimality and average optimality gaps (and their standard deviations) of DIF and SBD for instances not solved to optimality by at least one of them.
Recall that we provide optimality gaps proven by the solution  method within the time limit.

\begin{table}[!ht]
    \centering
    \scriptsize
    \caption{Average solution times and their standard deviations of DIF and SBD, for instances solved to optimality by both exact methods.}
    \begin{tabular}{crrr}
        \toprule
        \multirow{2}*{Instance attributes} & \# instances & \multicolumn{1}{c}{DIF} & \multicolumn{1}{c}{SBD} \\ \cmidrule{3-4}
        & considered & time (min)  & time (min) \\ \midrule
Benchmark&$293 \, (648)$&$5.72\pm11.58$&$1.11\pm2.91$\\
\midrule
$50$ locations&$151 \, (216)$&$4.15\pm10.60$&$1.15\pm3.30$\\
$100$ locations&$92 \, (216)$&$6.09\pm11.80$&$0.97\pm2.76$\\
$150$ locations&$50 \, (216)$&$9.77\pm13.13$&$1.26\pm1.72$\\
\midrule
$1|\mathcal{I}|$ customers&$152 \, (216)$&$3.01\pm8.16$&$0.79\pm2.82$\\
$3|\mathcal{I}|$ customers&$88 \, (216)$&$8.26\pm13.44$&$1.48\pm3.21$\\
$5|\mathcal{I}|$ customers&$53 \, (216)$&$9.27\pm14.58$&$1.43\pm2.56$\\
\midrule
5 periods&$126 \, (216)$&$4.36\pm9.85$&$1.15\pm2.69$\\
7 periods&$97 \, (216)$&$6.89\pm13.41$&$1.21\pm3.27$\\
9 periods&$70 \, (216)$&$6.53\pm11.66$&$0.91\pm2.79$\\
\midrule
1 facility&$166 \, (216)$&$4.88\pm10.32$&$0.58\pm1.08$\\
3 facilities&$69 \, (216)$&$7.91\pm14.73$&$1.75\pm3.96$\\
5 facilities&$58 \, (216)$&$5.50\pm10.58$&$1.87\pm4.41$\\
\midrule
Identical rewards&$163 \, (324)$&$6.12\pm11.28$&$0.46\pm0.87$\\
Different rewards&$130 \, (324)$&$5.22\pm11.98$&$1.93\pm4.13$\\
\midrule
Short rankings&$173 \, (324)$&$3.53\pm8.71$&$0.91\pm2.36$\\
Long rankings&$120 \, (324)$&$8.87\pm14.24$&$1.41\pm3.55$\\
\midrule
Constant demand&$126 \, (324)$&$5.55\pm10.87$&$1.08\pm2.49$\\
Sparse demand&$167 \, (324)$&$5.84\pm12.13$&$1.13\pm3.20$\\
        \bottomrule
    \end{tabular}
    \label{tab:bds-d-vs-a-time}
\end{table}
\begin{table}[!ht]
    \centering
    \scriptsize
    \caption{Number of instances solved to optimality, average optimality gaps and their standard deviations, of DIF and ABD, for instances not solved to optimality by at least one of these exact methods.}
    \begin{tabular}{crrrrrrrrr}
        \toprule
        \multirow{2}*{Instance attributes} & \# instances  & \multicolumn{2}{c}{DIF} &\multicolumn{2}{c}{SBD} \\ \cmidrule{3-6}
         & considered & \# opt. & prv. opt. gap (\%) & \# opt. & prv. opt. gap (\%)  \\ \midrule
Benchmark&$355 \, (648)$&$8$&$14.73\pm25.91$&$69$&$7.42\pm7.90$\\
\midrule
$50$ locations&$65 \, (216)$&$5$&$2.76\pm1.71$&$14$&$2.45\pm2.48$\\
$100$ locations&$124 \, (216)$&$3$&$10.20\pm14.81$&$17$&$6.88\pm5.93$\\
$150$ locations&$166 \, (216)$&$0$&$22.81\pm33.73$&$38$&$9.77\pm9.48$\\
\midrule
$1|\mathcal{I}|$ customers&$64 \, (216)$&$8$&$2.25\pm1.79$&$26$&$3.85\pm6.01$\\
$3|\mathcal{I}|$ customers&$128 \, (216)$&$0$&$12.21\pm19.89$&$15$&$8.57\pm8.75$\\
$5|\mathcal{I}|$ customers&$163 \, (216)$&$0$&$21.62\pm32.26$&$28$&$7.92\pm7.49$\\
\midrule
5 periods&$90 \, (216)$&$0$&$7.97\pm5.43$&$16$&$6.54\pm6.15$\\
7 periods&$119 \, (216)$&$4$&$13.57\pm23.76$&$21$&$7.34\pm7.37$\\
9 periods&$146 \, (216)$&$4$&$19.85\pm33.25$&$32$&$8.03\pm9.17$\\
\midrule
1 facility&$50 \, (216)$&$0$&$16.18\pm28.48$&$33$&$2.32\pm3.71$\\
3 facilities&$147 \, (216)$&$3$&$14.31\pm23.77$&$16$&$8.53\pm7.79$\\
5 facilities&$158 \, (216)$&$5$&$14.67\pm27.09$&$20$&$8.00\pm8.38$\\
\midrule
Identical rewards&$161 \, (324)$&$0$&$21.54\pm33.29$&$51$&$4.31\pm4.76$\\
Different rewards&$194 \, (324)$&$8$&$9.08\pm15.54$&$18$&$10.00\pm9.00$\\
\midrule
Short rankings&$151 \, (324)$&$4$&$8.95\pm13.94$&$19$&$7.98\pm8.03$\\
Long rankings&$204 \, (324)$&$4$&$19.01\pm31.36$&$50$&$7.01\pm7.80$\\
\midrule
Constant demand&$198 \, (324)$&$7$&$14.44\pm25.10$&$34$&$8.43\pm8.34$\\
Sparse demand&$157 \, (324)$&$1$&$15.10\pm26.96$&$35$&$6.15\pm7.14$\\
        \bottomrule
    \end{tabular}
    \label{tab:bds-d-vs-a-optgap}
\end{table}

\subsubsection{Comparison Between ABD and SBD.}

Table~\ref{tab:bda-time} reports average solution times, as well as their standard deviations, of SBD and ABD for instances solved to optimality by both methods.
Table~\ref{tab:bda-optgap} report the number of instances solved to optimality and average optimality gaps (and their standard deviations) of SBD and ABD for instances not solved to optimality by at least one of them.
Recall that we provide optimality gaps proven by the solution  method within the time limit.

\begin{table}[!ht]
    \centering
     \scriptsize
    \caption{Average solution times and their standard deviations of SBD and ABD, for instances solved to optimality by both exact methods.}
    \begin{tabular}{crrr}
        \toprule
        \multirow{2}*{Instance attributes} & \# instances & \multicolumn{1}{c}{SBD} & \multicolumn{1}{c}{ABD} \\ \cmidrule{3-4}
        & considered & time (min)  & time (min) \\ \midrule
Benchmark&$198 \, (216)$&$2.50\pm5.55$&$1.85\pm5.08$\\
\midrule
$50$ locations&$72 \, (72)$&$0.75\pm3.46$&$0.40\pm1.62$\\
$100$ locations&$66 \, (72)$&$2.51\pm5.54$&$2.21\pm5.83$\\
$150$ locations&$60 \, (72)$&$4.59\pm6.83$&$3.19\pm6.41$\\
\midrule
$1|\mathcal{I}|$ customers&$71 \, (72)$&$1.02\pm3.93$&$0.83\pm4.03$\\
$3|\mathcal{I}|$ customers&$65 \, (72)$&$2.66\pm5.43$&$1.56\pm3.81$\\
$5|\mathcal{I}|$ customers&$62 \, (72)$&$4.03\pm6.78$&$3.32\pm6.79$\\
\midrule
5 periods&$72 \, (72)$&$0.72\pm1.69$&$0.27\pm0.79$\\
7 periods&$69 \, (72)$&$2.53\pm4.94$&$2.01\pm5.29$\\
9 periods&$57 \, (72)$&$4.71\pm8.14$&$3.66\pm7.04$\\
\midrule
1 facility&$198 \, (216)$&$2.50\pm5.55$&$1.85\pm5.08$\\
\midrule
Identical rewards&$100 \, (108)$&$1.75\pm3.89$&$1.67\pm4.56$\\
Different rewards&$98 \, (108)$&$3.27\pm6.78$&$2.03\pm5.58$\\
\midrule
Short rankings&$107 \, (108)$&$1.31\pm3.60$&$0.81\pm2.52$\\
Long rankings&$91 \, (108)$&$3.89\pm6.97$&$3.07\pm6.80$\\
\midrule
Constant demand&$95 \, (108)$&$2.99\pm5.93$&$2.25\pm5.27$\\
Sparse demand&$103 \, (108)$&$2.05\pm5.17$&$1.48\pm4.90$\\
        \bottomrule
    \end{tabular}
    \label{tab:bda-time}
\end{table}
\begin{table}[!ht]
    \centering
    \scriptsize
    \caption{Number of instances solved to optimality, average optimality gaps and their standard deviations, of SBD and ABD, for instances not solved to optimality by at least one of these exact methods.}
    \begin{tabular}{crrrrrrrrr}
        \toprule
        \multirow{2}*{Instance attributes} & \# instances  & \multicolumn{2}{c}{SBD} &\multicolumn{2}{c}{ABD} \\ \cmidrule{3-6}
         & considered & \# opt. & prv. opt. gap (\%) & \# opt. & prv. opt. gap (\%)  \\ \midrule
Benchmark&$18 \, (216)$&$1$&$6.43\pm3.44$&$0$&$6.07\pm3.12$\\
\midrule
$50$ locations&$0 \, (72)$&$0$&$-$&$0$&$-$\\
$100$ locations&$6 \, (72)$&$1$&$5.45\pm2.78$&$0$&$5.15\pm2.67$\\
$150$ locations&$12 \, (72)$&$0$&$6.92\pm3.73$&$0$&$6.53\pm3.34$\\
\midrule
$1|\mathcal{I}|$ customers&$1 \, (72)$&$1$&$0.01\pm0.00$&$0$&$0.22\pm0.00$\\
$3|\mathcal{I}|$ customers&$7 \, (72)$&$0$&$6.50\pm2.92$&$0$&$5.83\pm2.71$\\
$5|\mathcal{I}|$ customers&$10 \, (72)$&$0$&$7.02\pm3.41$&$0$&$6.82\pm3.01$\\
\midrule
5 periods&$0 \, (72)$&$0$&$-$&$0$&$-$\\
7 periods&$3 \, (72)$&$0$&$4.47\pm1.75$&$0$&$5.73\pm1.58$\\
9 periods&$15 \, (72)$&$1$&$6.82\pm3.59$&$0$&$6.14\pm3.38$\\
\midrule
1 facility&$18 \, (216)$&$1$&$6.43\pm3.44$&$0$&$6.07\pm3.12$\\
\midrule
Identical rewards&$8 \, (108)$&$0$&$6.46\pm3.29$&$0$&$6.45\pm2.82$\\
Different rewards&$10 \, (108)$&$1$&$6.40\pm3.73$&$0$&$5.77\pm3.46$\\
\midrule
Short rankings&$1 \, (108)$&$0$&$2.23\pm0.00$&$0$&$0.57\pm0.00$\\
Long rankings&$17 \, (108)$&$1$&$6.67\pm3.37$&$0$&$6.40\pm2.89$\\
\midrule
Constant demand&$13 \, (108)$&$1$&$6.67\pm3.85$&$0$&$6.37\pm3.39$\\
Sparse demand&$5 \, (108)$&$0$&$5.79\pm2.25$&$0$&$5.30\pm2.43$\\
        \bottomrule
    \end{tabular}
    \label{tab:bda-optgap}
\end{table}

\subsection{Managerial Insights}

\subsubsection{Impact of Heuristic Decisions.}

Table~\ref{tab:heuristics} reports average opportunity gaps, as well as their standard deviations, of heuristics for instances with a known optimal solution. 
Recall that we define the opportunity gap for each heuristic as $\frac{Z^{\star} - Z^{\prime}}{Z^{\star}}$, where $Z^{\star}$ is the optimal objective value obtained through an exact method and $Z^{\prime}$ is the objective value of the heuristic at hand.
\begin{table}[!ht]
    \centering
    \scriptsize
    \caption{Average opportunity gaps and their standard deviations of DBH, RND, FGH and BGH, for instances with a known optimal solution.}
    \begin{tabular}{crrrrr}
        \toprule
        \multirow{2}*{Instance attributes} & \# instances & DBH & RND & FGH & BGH \\
        & considered & opp. gap (\%) & opp. gap (\%) & opp. gap (\%) & opp. gap (\%) \\ \midrule
Benchmark&$370 \, (648)$&$38.68\pm23.98$&$28.58\pm11.67$&$6.07\pm4.14$&$1.74\pm1.90$\\
\midrule
$50$ locations&$170 \, (216)$&$38.97\pm22.87$&$27.17\pm12.39$&$7.06\pm4.60$&$1.62\pm1.84$\\
$100$ locations&$112 \, (216)$&$39.22\pm24.75$&$29.98\pm11.93$&$5.87\pm3.92$&$2.10\pm2.14$\\
$150$ locations&$88 \, (216)$&$37.43\pm25.29$&$29.53\pm9.54$&$4.43\pm2.75$&$1.53\pm1.66$\\
\midrule
$1|\mathcal{I}|$ customers&$186 \, (216)$&$39.66\pm22.58$&$33.55\pm12.88$&$7.41\pm4.96$&$2.23\pm2.31$\\
$3|\mathcal{I}|$ customers&$103 \, (216)$&$38.40\pm24.46$&$24.09\pm8.25$&$5.36\pm2.76$&$1.27\pm1.31$\\
$5|\mathcal{I}|$ customers&$81 \, (216)$&$36.79\pm26.57$&$22.88\pm6.47$&$3.91\pm1.76$&$1.21\pm1.02$\\
\midrule
5 periods&$142 \, (216)$&$36.01\pm22.01$&$29.64\pm11.38$&$5.59\pm3.76$&$1.70\pm1.97$\\
7 periods&$122 \, (216)$&$39.48\pm24.35$&$28.48\pm11.82$&$6.27\pm4.22$&$1.66\pm1.75$\\
9 periods&$106 \, (216)$&$41.35\pm25.89$&$27.27\pm11.87$&$6.50\pm4.50$&$1.89\pm1.99$\\
\midrule
1 facility&$199 \, (216)$&$45.19\pm25.01$&$29.35\pm11.79$&$5.37\pm3.30$&$1.34\pm1.45$\\
3 facilities&$88 \, (216)$&$35.91\pm20.31$&$31.13\pm10.06$&$7.88\pm4.35$&$2.40\pm2.24$\\
5 facilities&$83 \, (216)$&$26.03\pm19.17$&$24.03\pm11.88$&$5.84\pm5.11$&$2.00\pm2.23$\\
\midrule
Identical rewards&$214 \, (324)$&$35.04\pm22.59$&$29.50\pm12.22$&$5.29\pm3.03$&$1.27\pm1.44$\\
Different rewards&$156 \, (324)$&$43.68\pm24.99$&$27.32\pm10.79$&$7.15\pm5.12$&$2.39\pm2.25$\\
\midrule
Short rankings&$196 \, (324)$&$43.08\pm23.19$&$30.62\pm12.06$&$6.33\pm3.58$&$1.82\pm2.05$\\
Long rankings&$174 \, (324)$&$33.73\pm23.96$&$26.28\pm10.81$&$5.78\pm4.70$&$1.65\pm1.73$\\
\midrule
Constant demand&$167 \, (324)$&$60.99\pm14.51$&$25.22\pm11.36$&$5.43\pm3.93$&$0.63\pm0.84$\\
Sparse demand&$203 \, (324)$&$20.33\pm11.31$&$31.34\pm11.22$&$6.60\pm4.25$&$2.65\pm2.05$\\
        \bottomrule
    \end{tabular}
    \label{tab:heuristics}
\end{table}

\end{document}